%% file: ComparisonOfNervesMay2022.tex
\documentclass[10pt]{amsart}
\input{PreambleSep2020.tex}
\title[Model independence of $(\infty,2)$-categorical nerves]{Model independence of $(\infty,2)$-categorical nerves}

\author{Lyne Moser}
\address{Max Planck Institute for Mathematics, Bonn, Germany}
\email{moser@mpim-bonn.mpg.de}

\author{Viktoriya Ozornova}
\address{Max Planck Institute for Mathematics, Bonn, Germany}
\email{viktoriya.ozornova@mpim-bonn.mpg.de}

\author{Martina Rovelli}
\address{Department of Mathematics and Statistics,
University of Massachusetts, 
Amherst,
USA
}
\email{rovelli@math.umass.edu}


\begin{document}
\begin{abstract}
For most models of $(\infty,2)$-categories an embedding of the $\infty$-category of $2$-categories into that of $(\infty,2)$-categories has been constructed in the form of a nerve construction of some flavor. We prove that all those nerve embeddings induce equivalent functors, modulo change of model. We also show that all the nerve embeddings realize the $\infty$-category of $2$-categories as the sub-$\infty$-category of $(\infty,2)$-categories that are local with respect to a certain class of maps.
\end{abstract}

\maketitle
\tableofcontents

\section*{Introduction}

%%% Infinity,2-categories
It has become apparent that many phenomena of interest, such as the cobordism hypothesis, can only be properly formalized using the language of higher categories, often in the form of $(\infty,n)$-categories for $n\ge0$, and this paper is concerned with $(\infty,2)$-categories. The structure of an $(\infty,2)$-category could be summarized as a weakening of the structure present in a traditional $2$-category. It consists of objects, $1$- and $2$-morphisms that compose suitably, as well as higher weakly invertible morphisms in dimension higher than $2$ that serve as witnesses for relations between lower dimensional morphisms.

%%% Model comparisons
Many mathematical objects have been proposed to formalize $(\infty,2)$-categories, each model presenting its own advantages and disadvantages. These include Barwick's $2$-fold complete Segal spaces \cite{BarwickThesis}, Verity's saturated $2$-complicial sets \cite{VerityComplicialI,VeritySlides,EmilyNotes,or,RiehlVerityBook}, Lurie's $\infty$-bicategories \cite{LurieGoodwillie}, Rezk's complete Segal $\Theta_2$-spaces \cite{rezkTheta}, Ara's $2$-quasi-categories \cite{ara}, and $2$-comical sets \cite{CKM,DKM}, as well as categories strictly enriched over a model of $(\infty,1)$-categories \cite{LurieGoodwillie,br1,br2}.
In the past few years, the proof that all models are equivalent was completed, combining work by  Lurie \cite{LurieGoodwillie}, Bergner--Rezk \cite{br1,br2},  Ara \cite{ara}, Gagna--Harpaz--Lanari \cite{GHL}, Campion--Doherty--Kapulkin--Maehara \cite{CKM,DKM}.

It is often the case that the same construction gets implemented independently into two or more models. It is then necessary to verify that they indeed encode the same construction, modulo a change of model given by a direct comparison or a zigzag of such. In this paper, we specifically address the compatibility of several embeddings of the homotopy theory of $2$-categories into that of $(\infty,2)$-categories that have been constructed for different models.\footnote{In the past, many ways to associate to any $2$-category a classifying \emph{space} -- namely an $(\infty,0)$-category, as opposed to an $(\infty,2)$-category -- have been provided by Street \cite{StreetOrientedSimplexes}, Duskin \cite{duskin}, Bullejos--Cegarra \cite{BullejosCegarra}, Lack--Paoli \cite{LackPaoli2Nerves}--and the equivalence of such constructions as spaces is proven in \cite{CCG}.} 

%%% Embedding of 2-categories
By design, the idea of an $(\infty,2)$-category is supposed to weaken and generalize the notion of a strict $2$-category. In particular, it is expected that any question about the homotopy theory of $2$-categories should be equivalently addressable in the world of $2$-categories or in that of $(\infty,2)$-categories. This requirement, which is even partially axiomatized in the abstract setup by Barwick--Schommer-Pries \cite{BarwickSchommerPries}, could be phrased by expecting an embedding of the homotopy theory of $2$-categories into that of $(\infty,2)$-categories. Beside for providing a consistency check, the embedding of $2$-categories into $(\infty,2)$-categories is crucial in that several structural components of $(\infty,2)$-categories, such as pasting schemes, are parametrized by strict $2$-categories.

%%% State of the art
The analog question for $(\infty,0)$-categories (a.k.a.\ $\infty$-groupoids) and $(\infty,1)$-categories (a.k.a.\ $\infty$-categories) is equally valid although easier to address and by now fairly understood.
In essentially all models for $(\infty,0)$- and $(\infty,1)$-categories one can easily identify or find in the literature a simple nerve construction for $0$-categories (a.k.a.\ sets) and $1$-categories and prove that this nerve construction realizes an embedding of homotopy theories into $(\infty,0)$- and $(\infty,1)$-categories, respectively. 

For $(\infty,2)$-categories the situation is more subtle. For instance, when one works with model categories, one technical difficulty is the fact that most models don't admit a homotopical nerve embedding that is at once fully faithful at the pointset level and a right Quillen functor at the model categorical level.
In the recent years a well-behaved embedding has also been constructed in most models of $(\infty,2)$-categories presented by model categories in the form of a homotopical functor that is homotopically fully faithful, which is in addition either right Quillen or fully faithful (but generally not both). This was achieved by the second and third author \cite{Nerves2Cat} for $2$-complicial sets, by Campbell \cite{CampbellHoCoherent} for $2$-quasi-categories, by Gagna--Harpaz--Lanari for scaled simplicial sets \cite{GHL}, and by the first author \cite{MoserNerve} for $2$-fold complete Segal spaces. For $(\infty,2)$-categories presented by categories enriched over a model of $(\infty,1)$-categories, this can be done by base-change along a suitable $1$-dimensional nerve.

The first result of this paper, proven as \cref{MainTheorem} is to check that all the mentioned nerve constructions (along with a few more that we add) are compatible with each other via the known model comparisons.

\begin{thmalph}
The aforementioned nerve embeddings of $2$-categories into $(\infty,2)$-categories constructed in different model categories are compatible with each other via known equivalences of models.
\end{thmalph}

At the level of $\infty$-categories, as part of a more general machinery Gepner--Haugseng \cite{GH} identified that the $\infty$-category of $2$-categories can be understood as a localization of the $\infty$-category of $(\infty,2)$-categories. More precisely, $2$-categories are exactly the $(\infty,2)$-categories that are local with respect to the $2$-fold $2$-point suspension of the inclusion of a point into a positive-dimensional sphere.

We prove as \cref{MainTheorem2} that all the considered nerve embeddings induce at the level of $\infty$-categories precisely the inclusion of $2$-categories as local objects amongst $(\infty,2)$-categories with respect to the class of maps from the previous paragraph.

\begin{thmalph}
The aforementioned  nerve embeddings of $2$-categories into $(\infty,2)$-categories constructed in different model categories implement the embedding of $2$-categories as local $(\infty,2)$-categories.
\end{thmalph}

While overall expected, the compatibility of the nerve constructions in different models is a fundamental verification for the consistency of the theory, and a necessary ingredient in phrasing model independently many statements originally proven in a specific model.

To mention one example, in the paper \cite{HORR} the second and third author proved with Hackney and Riehl an $(\infty,2)$-dimensional pasting theorem for $(\infty,2)$-categories modeled by categories enriched over quasi-categories, and it is there explained how the compatibility of nerves which is the subject of the current paper is necessary to conclude that the pasting theorem holds in all other models.

The compatibility of nerves is expected to play a similar role in other circumstances, for instance in work in progress by the first and third author with Rasekh with the goal of developing a model independent theory of weighted limits valued in $(\infty,2)$-categories.

Beside the novel result, we are also taking this project as an opportunity to write an accessible paper that surveys over the different nerve constructions and how they relate to each other, helping a non-specialist navigate the complex literature of $(\infty,2)$-categories.

 \addtocontents{toc}{\protect\setcounter{tocdepth}{1}}
\subsection*{Acknowledgements}
We are thankful to Rune Haugseng and Lennart Meier for valuable conversations, and to tslil clingman for their help with \LaTeX{} and TikZ.
 This work was completed while the authors
visited the 
 Instituto de Matemáticas de UNAM in Cuernavaca for the program \emph{Higher categories -- Part 2}, supported by the National Science Foundation under Grant No. DMS-1928930. The third author is grateful for support from the National Science
Foundation under Grant No. DMS-2203915.

\section{Organization of the paper}

\subsection{Model categorical framework}
In this paper, we use the language of model categories to formalize the $\infty$-categories of $(\infty,2)$-categories presented by different models.  We refer the reader to e.g.~\cite{hovey,hirschhorn} for the basic definitions from model category theory. We also assume familiarity with the basics of $\infty$-categories in the form of quasi-categories, see e.g.~\cite{htt}. Here, we only briefly recall the key facts needed to interpret the model categorical statements as statements about homotopy theories and $\infty$-categories.
\begin{itemize}[leftmargin=*]
    \item Any model category $\cM$ has an underlying $\infty$-category $[\cM]_\infty$. Explicitly, the $\infty$-category $[\cM]_\infty$ is obtained as the homotopy coherent nerve of a fibrant replacement of the \emph{Hammock localization} of $\cM$; see  e.g.~\cite{DwyerKanCalculating,DwyerKanFunction} or \cref{NerveComparisonAppendix}
    for more details. For this specific model of $[\cM]_\infty$, the set of objects is the same as the sets of objects of $\cM$;
    \item Any homotopical functor $F\colon\cM\to\cM'$, i.e. a functor that preserves weak equivalences, induces a functor of $\infty$-categories $[F]_\infty\colon[\cM]_{\infty}\to[\cM']_{\infty}$. It can be computed on objects as $[F]_\infty(X)=F(X)$.
   \item Any right (resp.~left) Quillen functor $F\colon\cM\to\cM'$ induces right (resp.~left) adjoint functor
   of $\infty$-categories $[F]_\infty\colon[\cM]_\infty\to[\cM']_\infty$, as proven in \cite[Thm\ 2.1]{MazelGeeQuillenAdj}.
   It can be computed on objects as $[F]_\infty(X)\simeq F(X^{\textrm{fib}})$ (resp.~$[F]_\infty(X)\simeq F(X^{\textrm{cof}})$). Here, $X^{\textrm{fib}}$ (resp.~$X^{\textrm{cof}}$) denotes a fibrant (resp.~cofibrant) replacement of $X$ in~$\cM$.
   \item Any right (resp.~left) Quillen embedding\footnote{By a \emph{right Quillen embedding} we mean a right Quillen functor in which the derived counit of any fibrant object is a weak equivalence. This is the right Quillen functor occurring in what is known in the literature as a \emph{Quillen reflection} or \emph{homotopy reflection} introduced in \cite[\S6.3]{JoyalVolumeII}.
   A \emph{left Quillen embedding} is defined dually.}  $F\colon\cM\to\cM'$ induces a fully faithful right adjoint (resp.~left adjoint) of $\infty$-categories $[F]_\infty\colon[\cM]_\infty\to[\cM']_\infty$.
   \item Any left (resp.~right) Quillen equivalence $F\colon\cM\to\cM'$ induces an equivalence of $\infty$-categories $[F]_\infty\colon[\cM]_\infty\to[\cM']_\infty$, as a consequence of what discussed in \cite[\S A.2]{MazelGeeQuillenAdj} and \cite[\S1.3.4]{LurieHA}.
   In particular, a zigzag of Quillen equivalences induces an equivalence of $\infty$-categories.
   \item If a functor $F\colon\cM\to\cM'$ is such that it induces a functor $[F]_\infty\colon[\cM]_\infty\to[\cM']_\infty$ in more than one way, for instance it is both left and right Quillen, or it is both right Quillen and homotopical, the resulting functors are canonically equivalent.
   \item If functors $F\colon\cM\to\cM'$ and $F'\colon\cM'\to\cM''$ and their composite $F'\circ F\colon\cM\to\cM''$ induce functors of $\infty$-categories $[F]_\infty\colon[\cM]_\infty\to[\cM']_\infty$, $[F']_\infty\colon[\cM']_\infty\to [\cM'']_\infty$, and ${[F']_\infty\circ [F]_\infty\colon[\cM]_\infty\to[\cM'']_\infty}$ each computed using any of the rules described above, then there is a canonical equivalence $[F'\circ F]_\infty\simeq [F']_\infty\circ [F]_\infty$.
\end{itemize}

 \subsection{Models of \texorpdfstring{$(\infty,2)$}{(infinity,2)}-categories}
 \label{allmodels}
We briefly recall also the main different approaches to modeling $(\infty,2)$-categories that will be relevant for the paper. For each of the approaches, it is possible to realize the homotopy theory of $(\infty,2)$-categories by means of a model structure
in which the $(\infty,2)$-categories are precisely the fibrant objects. 

\begin{enumerate}[leftmargin=*,label=(\alph*)]
    \item \emph{Globular models:} based on presheaves over Joyal's disk category $\Theta_2$ \cite{JoyalDisks} or variants of it. They include Ara's $2$-quasi-categories \cite{ara} and Rezk's complete Segal $\Theta_2$-spaces \cite{rezkTheta}. The supporting model structures $\psh{\Theta_2}_{(\infty,2)}$ and $\spsh{\Theta_2}_{(\infty,2)}$ will be recalled in more detail in \cref{AraModelStructure,RezkModelStructure}.
    \item \emph{Bisimplicial models:} based on presheaves over $\Delta\times\Delta$. They include Barwick's $2$-fold complete Segal spaces \cite{BarwickThesis} and Bergner--Rezk's Segal precategories \cite{br1}. The supporting model structures $\spsh{(\Delta\times\Delta)}_{(\infty,2)}$ and $P\cat(\spsh{\Delta})_{(\infty,2)}$
will be recalled in more detail in \cref{BarwickModelStructure,PrecatSubsection}.
    \item \emph{Enriched models:} based on categories \emph{strictly} enriched over a model of $(\infty,1)$-ca\-te\-gories. They include categories enriched over Joyal's quasi-categories \cite{joyalnotes}, over Rezk's complete Segal spaces \cite{rezk} and over Lurie's marked simplicial sets \cite{htt}. The supporting model structures $\vcat{\sset_{(\infty,1)}}$, $\vcat{\spsh{\Delta}_{(\infty,1)}}$ and $\vcat{\sset^+_{(\infty,1)}}$ will be recalled in more detail in \cref{EnrichedCatMS}.
    \item \emph{Simplicial models:} based on presheaves over variants of the simplex category~$\Delta$. They include Verity's saturated $2$-complicial sets \cite{VeritySlides,or,RiehlVerityBook}, Lurie's $\infty$-bicategories \cite{LurieGoodwillie}, and saturated $2$-precomplicial sets \cite{or} by the second and third author. The supporting models structures $\msset_{(\infty,2)}$, $\sset^{sc}_{(\infty,2)}$ and $\psh{t\Delta}_{(\infty,2)}$ will be recalled in more detail in \cref{LurieModelStructure,ORModelStructure,VerityModelStructure}.
\item \emph{Cubical models:} based on presheaves over a suitable category of cubes. The main incarnation is given by Doherty--Kapulkin--Maehara's $2$-comical sets \cite{DKM}, supported by the model structure $mc\set_{(\infty,2)}$, which is a variant of a previous version by Campion--Kapulkin--Maehara \cite{CKM}.
\end{enumerate}

We know that these models of $(\infty,2)$-categories have equivalent homotopy theories because the supporting model structures are connected by the following zigzags of Quillen equivalences\footnote{Three of these Quillen equivalences are denoted $R$ in the original sources. To distinguish them in this paper, we are using $R$, $\mathrm{Refl}$ and $\mathfrak R$.}.
\begin{tz}
\node[](1) {$\psh{\Theta_2}_{(\infty,2)}$}; 
\node[below of=1](2) {$\spsh{\Theta_2}_{(\infty,2)}$}; 
\node[below of=2](3) {$\spsh{(\Delta\times \Delta)}_{(\infty,2)}$}; 
\node[below of=3](4) {$P\cat(\spsh{\Delta})_{(\infty,2)}$};
\node[right of=4,xshift=2cm,yshift=-2.3pt](5) {$\vcat{\spsh{\Delta}_{(\infty,1)}}$};  
\node[right of=5,xshift=1.5cm,yshift=1.7pt](6) {$\vcat{\sset_{(\infty,1)}}$};
\node[right of=6,xshift=1.5cm,yshift=-1.7pt](7) {$\vcat{\sset^+_{(\infty,1)}}$};
\node[above of=7,yshift=1pt](8) {$\sset^{sc}_{(\infty,2)}$}; 
\node[above of=8](9) {$\msset_{(\infty,2)}$};
\node[above of=6,yshift=1.5cm,yshift=-1pt](8') {$mc\set_{(\infty,2)}$};
\node[above of=9](10) {$\psh{t\Delta}_{(\infty,2)}$};

\draw[->] (2) to node[left,la]{$(-)_{\bullet,0}$} node[right,la]{\cite{ara}} (1);
\draw[->] (3) to node[left,la]{$d_*$} node[right,la]{\cite{br2}} (2);
\draw[->] (3) to node[left,la]{$R$} node[right,la]{\cite{br2}} (4);
\draw[->] ($(5.west)+(0,2.3pt)$) to node[above,la]{\cite{br1}} node[below,la]{$\BRequi$} (4);
\draw[->] ($(5.east)+(0,2.3pt)$) to node[below,la]{$((-)_{\bullet,0})_*$} ($(6.west)+(0,0.6pt)$);
\draw[->] ($(7.west)+(0,2.3pt)$) to node[below,la]{$U_*$} ($(6.east)+(0,0.6pt)$);
\draw[->] (7) to node[right,la]{$\schoco$} node[left,la]{\cite{LurieGoodwillie}} (8);
\draw[->] (9) to node[right,la]{$U$} node[left,la]{\cite{GHL}} (8);
\draw[->] (10) to node[right,la]{$\mathrm{Refl}$} node[left,la]{\cite{or}} (9);
\draw[->] (8') to node[above,la]{$T$} node[below,la]{\cite{DKM}} (9);
\end{tz} 

 \subsection{Models of \texorpdfstring{$(\infty,2)$}{(infinity,2)}-categorical nerves}
 
 \label{allnerves}
The canonical homotopy theory of strict $2$-categories is presented by the following model structure, due to Lack.
 \begin{thm}[\cite{lack1,lack2}]
 \label{lackMS}
The category $\twocat$ of small $2$-categories and $2$-functors admits a model structure in which
\begin{itemize}[leftmargin=*]
    \item all $2$-categories are fibrant,
    \item the weak equivalences are precisely the biequivalences, and 
    \item the trivial fibrations are precisely the $2$-functors that are surjective on objects, full on $1$-morphisms, and fully faithful on $2$-morphisms.
\end{itemize}
 \end{thm}
 
 Several nerve constructions for $2$-categories valued in a model of $(\infty,2)$-categories have been constructed in the form of a right Quillen embedding.\footnote{Given the numerous nerve constructions considered in this paper, for the sake of exposition we decided to change some of their notations to something more evocative of which model they refer to. We point out when these constructions are recalled what is the notation used in the original sources.}

\begin{enumerate}[leftmargin=0.75cm]
\item[(a)] \emph{Nerve into $2$-quasi-categories:} A functor
\[\NLein\colon \twocat\to\spsh{\Theta_2}_{(\infty,2)}\]
was first considered by Leinster \cite[Def.\ J]{LeinsterSurvey} and later shown by Campbell \cite[Rmk 5.16,~Thm 5.10]{CampbellHoCoherent} to be a right Quillen embedding. This nerve and its properties will be recalled in \cref{LeinsterNerve,LeinsterEmbedding}.
\item[(b)] \emph{Nerve into $2$-fold complete Segal spaces:} A functor
\[\NMos\colon \twocat\to\spsh{(\Delta\times\Delta)}_{(\infty,2)}\]
was constructed and shown to be a right Quillen embedding by the first author in \cite[\S5.1, Thms~6.1.1,~6.1.3]{MoserNerve}. This nerve and its properties will be recalled in \cref{MoserNerve,MoserEmbedding}.
\item[(c)] \emph{Nerve into categories enriched over quasi-categories:} A functor
\[{\bf N}_*\colon\twocat\to\vcat{\sset_{(\infty,1)}}\]
obtained by base-change along the usual nerve functor, is used e.g.~in \cite[\S1.4.2]{RiehlVerityBook}, and can be shown to be a right Quillen embedding. This nerve and its properties will be recalled in \cref{EnrichedNerves}.
\item[(d)] \emph{Nerve into $\infty$-bicategories:} A functor
\[\NGHL\colon\twocat\to\msset^{sc}_{(\infty,2)}\]
was considered by Harpaz--Nuiten--Prasma in \cite[\S2]{HarpazNuitenPrasmaCohomology} and shown to be a right Quillen embedding by Gagna--Harpaz--Lanari in \cite[Prop.\ 8.2, 8.3]{GHL}.
This nerve and its properties will be recalled in \cref{GHLNerve,GHLEmbedding}.
\item[(d')] \emph{Nerve into $2$-precomplicial sets:} A functor
\[\NOR\colon\twocat\to\psh{t\Delta}_{(\infty,2)}\]
was constructed and shown to be a right Quillen embedding by the second and third authors in \cite[Thm~4.12,~Cor.~~4.13]{Nerves2Cat}.
This nerve and its properties will be recalled in \cref{ORNerve,OREmbedding}. 
\end{enumerate}

\subsection{Equivalences of the nerve constructions}
The goal of this paper is to study how all those nerve constructions interact with the model comparison functors and prove the compatibility. In practice, this amounts to considering the following diagram of (model) categories,
\begin{diagram} \label{diagramnerves}
\node[](1) {$\psh{\Theta_2}_{(\infty,2)}$}; 
\node[below of=1](2) {$\spsh{\Theta_2}_{(\infty,2)}$}; 
\node[below of=2](3) {$\spsh{(\Delta\times \Delta)}_{(\infty,2)}$}; 
\node[below of=3](4) {$P\cat(\spsh{\Delta})_{(\infty,2)}$};
\node[right of=4,xshift=2cm,yshift=-2.3pt](5) {$\vcat{\spsh{\Delta}_{(\infty,1)}}$};  
\node[right of=5,xshift=1.5cm,yshift=1.7pt](6) {$\vcat{\sset_{(\infty,1)}}$}; 
\node[right of=6,xshift=1.5cm,yshift=-1.7pt](7) {$\vcat{\sset^+_{(\infty,1)}}$};
\node[above of=7,yshift=1pt](8) {$\sset^{sc}_{(\infty,2)}$}; 
\node[above of=8](9) {$\msset_{(\infty,2)}$};
\node[above of=9](10) {$\psh{t\Delta}_{(\infty,2)}$};
\node(cat) at ($(2)!0.5!(9)+(0,1pt)$) {$\twocat$}; 

\draw[->] (2) to node[left,la]{$(-)_{\bullet,0}$} (1);
\draw[->] (3) to node[left,la]{$d_*$} (2);
\draw[->] (3) to node[left,la]{$R$} (4);
\draw[->] ($(5.west)+(0,2.3pt)$) to node[below,la]{$\BRequi$} (4);
\draw[->] ($(5.east)+(0,2.3pt)$) to node[below,la]{$((-)_{\bullet,0})_*$} ($(6.west)+(0,0.6pt)$);
\draw[->] ($(7.west)+(0,2.3pt)$) to node[below,la]{$U_*$} ($(6.east)+(0,0.6pt)$);
\draw[->] (7) to node[right,la]{$\schoco$} (8);
\draw[->] (9) to node[right,la]{$U$} (8);
\draw[->] (10) to node[right,la]{$\mathrm{Refl}$} (9);

\draw[->] (cat) to node[above,la,yshift=3pt]{$\NLein$} (1); 
\draw[->] (cat) to node[below,la,yshift=-3pt]{$\NMos$} (3); 
\draw[->] (cat) to node[left,la]{$\bf N_*$} (6); 
\draw[->] (cat) to node[below,la,yshift=-3pt]{$\NGHL$} (8); 
\draw[->] (cat) to node[above,la]{$\NOR$} (10); 
\end{diagram}
built using some of the model comparison functors and the aforementioned nerve constructions, and show that all regions induce commutative diagrams at the level of underlying $\infty$-categories.

\begin{thm}
\label{MainTheorem}
The diagram of underlying $\infty$-categories induced by \eqref{diagramnerves} commutes up to equivalence.

\newcounter{region}

\begin{tz}
\node[](1) {$[\psh{\Theta_2}_{(\infty,2)}]_\infty$}; 
\node[below of=1](2) {$[\spsh{\Theta_2}_{(\infty,2)}]_\infty$}; 
\node[below of=2](3) {$[\spsh{(\Delta\times \Delta)}_{(\infty,2)}]_\infty$}; 
\node[below of=3](4) {$[P\cat(\spsh{\Delta})_{(\infty,2)}]_\infty$};
\node[right of=4,xshift=2.3cm,yshift=-2pt](5) {$[\vcat{\spsh{\Delta}_{(\infty,1)}}]_\infty$};  
\node[right of=5,xshift=2.2cm,yshift=1.2pt](6) {$[\vcat{\sset_{(\infty,1)}}]_\infty$};
\node[right of=6,xshift=1.6cm,yshift=-1.2pt](7) {$[\vcat{\sset^+_{(\infty,1)}}]_\infty$};
\node[above of=7,yshift=1pt](8) {$[\sset^{sc}_{(\infty,2)}]_\infty$}; 
\node[above of=8](9) {$[\msset_{(\infty,2)}]_\infty$};
\node[above of=9](10) {$[\psh{t\Delta}_{(\infty,2)}]_\infty$};
\node(cat) at ($(2)!0.5!(9)$) {$[\twocat]_\infty$}; 

\draw[->] (2) to node[left,la]{$[(-)_{\bullet,0}]_\infty$} (1);
\draw[->] (3) to node[left,la]{$[d_*]_\infty$} (2);
\draw[->] (3) to node[left,la]{$[R]_\infty$} (4);
\draw[->] ($(5.west)+(0,2pt)$) to node[below,la,xshift=3pt]{$[\BRequi]_\infty$} (4);
\draw[->] ($(5.east)+(0,2pt)$) to node[below,la,xshift=4pt]{$[((-)_{\bullet,0})_*]_\infty$} ($(6.west)+(0,.8pt)$);
\draw[->] ($(7.west)+(0,2pt)$) to node[below,la,xshift=4pt]{$[U_*]_\infty$} ($(6.east)+(0,.8pt)$);
\draw[->] (7) to node[right,la]{$[\schoco]_\infty$} (8);
\draw[->] (9) to node[right,la]{$[U]_\infty$} (8);
\draw[->] (10) to node[right,la]{$[\mathrm{Refl}]_\infty$} (9);

\draw[->] (cat) to node[above,la,yshift=3pt]{$[\NLein]_\infty$} (1); 
\draw[->] (cat) to node[below,la,yshift=-3pt]{$[\NMos]_\infty$} (3); 
\draw[->] (cat) to node[left,la]{$[\bf N_*]_\infty$} (6); 
\draw[->] (cat) to node[below,la,yshift=-3pt]{$[\NGHL]_\infty$} (8); 
\draw[->] (cat) to node[above,la,yshift=3pt]{$[\NOR]_\infty$} (10); 

\refstepcounter{region}
\node[scale=0.9,right of=2,xshift=.5cm] {\emph{(\theregion)}\label{ThetaDelta}};
\refstepcounter{region}
\node[scale=0.9,above of=5,yshift=-.3cm,xshift=-.5cm]{\emph{(\theregion)}\label{DeltaEnriched}};
\refstepcounter{region}
\node[scale=0.9,above of=6,yshift=-.3cm,xshift=.7cm]{\emph{(\theregion)}\label{LevJoyal}};
\refstepcounter{region}
\node[scale=0.9,left of=9,xshift=-.7cm] {\emph{(\theregion)}\label{MarkedModels}};
\end{tz}
\end{thm}

\begin{proof}[Outline of the proof]
We address the commutativity of each of the regions as follows.
 \begin{itemize}[leftmargin=*]
    \item The fact that the region (\ref{ThetaDelta}) commutes is addressed as \cref{corLeinsterVSMoser}.
    \item The fact that the region (\ref{DeltaEnriched}) commutes is addressed as a combination of \cref{corcomparisonenriched,corcomparisonPrecat}.
    \item The fact that the region (\ref{LevJoyal}) commutes is addressed as a combination of \cref{corcomparisonenriched,corComplicialEnriched}.
    \item The fact that the region (\ref{MarkedModels}) commutes is addressed as \cref{corsimplicial}.\qedhere
 \end{itemize}
\end{proof}

\subsection{Universal property of nerve embeddings}

In the following diagram of adjunctions of $\infty$-categories, \cref{MainTheorem} guarantees that the diagram involving the functors induced by the nerve construction functors commutes up to equivalence. Hence so does the one involving the left adjoints
to the functors induced by the nerve constructions.
\begin{diagram} \label{DiagramAdjunctions}
\node[](1) {$[\psh{\Theta_2}_{(\infty,2)}]_\infty$}; 
\node[below of=1](2) {$[\spsh{\Theta_2}_{(\infty,2)}]_\infty$}; 
\node[below of=2](3) {$[\spsh{(\Delta\times \Delta)}_{(\infty,2)}]_\infty$}; 
\node[below of=3](4) {$[P\cat(\spsh{\Delta})_{(\infty,2)}]_\infty$};
\node[right of=4,xshift=2.15cm,yshift=-2pt](5) {$[\vcat{\spsh{\Delta}_{(\infty,1)}}]_\infty$};  
\node[right of=5,xshift=1.55cm,yshift=1.2pt](6) {$[\vcat{\sset_{(\infty,1)}}]_\infty$}; 
\node[right of=6,xshift=1.55cm,yshift=-1.2pt](7) {$[\vcat{\sset^+_{(\infty,1)}}]_\infty$};
\node[above of=7,yshift=1pt](8) {$[\sset^{sc}_{(\infty,2)}]_\infty$}; 
\node[above of=8](9) {$[\msset_{(\infty,2)}]_\infty$};
\node[above of=9](10) {$[\psh{t\Delta}_{(\infty,2)}]_\infty$};
\node(cat) at ($(2)!0.5!(9)$) {$[\twocat]_\infty$}; 

\node at ($(1.south)!0.5!(2.north)$) {\rotatebox{270}{$\eqarrow$}};

\node at ($(2.south)!0.5!(3.north)$) {\rotatebox{270}{$\eqarrow$}};

\node at ($(3.south)!0.5!(4.north)$) {\rotatebox{270}{$\eqarrow$}};

\node at ($(4.east)!0.5!(5.west)$) {$\eqarrow$};

\node at ($(5.east)!0.5!(6.west)$) {$\eqarrow$};
\node at ($(6.east)!0.5!(7.west)$) {$\eqarrow$};

\node at ($(10.south)!0.5!(9.north)$) {\rotatebox{270}{$\longeqarrow$}};

\node at ($(9.south)!0.5!(8.north)$) {\rotatebox{270}{$\longeqarrow$}};

\node at ($(8.south)!0.5!(7.north)$) {\rotatebox{270}{$\longeqarrow$}};

\draw[->] ($(cat.west)+(0,2pt)$) to node[above,la,sloped,xshift=-2pt,yshift=-2pt]{\rotatebox{90}{$\vdash$}} node[below,la,yshift=-3pt,xshift=-3pt]{$[\NLein]_\infty$} ($(1.east)-(0,7pt)$); 
\draw[->] ($(1.east)+(0,1pt)$) to ($(cat.west)+(0,10pt)$); 

\draw[->] ($(cat.west)-(0,10pt)$) to node[above,la,sloped,xshift=-2pt,yshift=-2pt]{\rotatebox{90}{$\vdash$}} node[below,la,yshift=-3pt,xshift=9pt]{$[\NMos]_\infty$} ($(3.east)-(0,1pt)$); 
\draw[->] ($(3.east)+(0,7pt)$) to ($(cat.west)-(0,2pt)$); 

\draw[->] ($(cat.south)-(5pt,0)$) to node[above,la,sloped,xshift=-2pt,yshift=-2pt]{\rotatebox{270}{$\dashv$}} node[left,la,xshift=-2pt,yshift=-2pt]{$[\bf N_*]_\infty$} ($(6.north)+(1pt,0)$);
\draw[->] ($(6.north)+(10pt,0)$) to ($(cat.south)+(4pt,0)$);

\draw[->] ($(cat.east)-(0,10pt)$) to node[above,la,sloped,xshift=-2pt,yshift=-2pt]{\rotatebox{90}{$\vdash$}} node[below,la,yshift=-3pt,xshift=-3pt]{$[\NGHL]_\infty$} ($(8.west)-(0,1pt)$); 
\draw[->] ($(8.west)+(0,7pt)$) to ($(cat.east)-(0,2pt)$); 

\draw[->] ($(cat.east)+(0,2pt)$) to node[above,la,sloped,xshift=-2pt,yshift=-2pt]{\rotatebox{90}{$\vdash$}} node[below,la,yshift=-3pt,xshift=3pt]{$[\NOR]_\infty$} ($(10.west)-(0,7pt)$); 
\draw[->] ($(10.west)+(0,1pt)$) to ($(cat.east)+(0,10pt)$); 
\end{diagram}

Any of the $\infty$-categories underlying one of the model structures for $(\infty,2)$-categories from \cref{allmodels} can be taken to be \emph{the} $\infty$-category of $(\infty,2)$-categories $\mathscr{C}at_{(\infty,2)}$, and all others are equivalent to this one -- explicitly via the equivalences of $\infty$-categories given by the mentioned Quillen equivalences. Also, if $\mathscr{C}at_{2}$ denotes the $\infty$-category of $2$-categories, then there is an equivalence of $\infty$-categories $\mathscr{C}at_{2}\simeq [\twocat]_\infty$.
The many models of nerves discussed in \cref{allnerves}, all induce equivalent right adjoint functors between
the $\infty$-category of $(\infty,2)$-categories $\mathscr{C}at_{(\infty,2)}$ and the $\infty$-category $[\twocat]_\infty$ of $2$-categories:
\begin{equation}
\label{genericadjunction}
  \mathscr{C}at_{(\infty,2)}\simeq [\cM]_\infty\rightleftarrows[\twocat]_\infty\simeq\mathscr{C}at_{2}.  
\end{equation}
One may argue at this point that, although it was shown that those functors \emph{do the same thing}, do they actually do \emph{the right thing}?

To address this question, we first observe the compatibility of the embedding $\mathscr Cat_2\hookrightarrow\mathscr Cat_{(\infty,2)}$ with Barwick--Schommer-Pries' framework from \cite{BarwickSchommerPries}.

\begin{rmk}
In \cite[\S7]{BarwickSchommerPries}, Barwick--Schommer-Pries identified an axiomatic setup that guarantees that an $\infty$-category $\mathscr M$ (with extra structure) models correctly the theory of $(\infty,2)$-categories, satisfying in particular $\mathscr M\simeq\mathscr Cat_{(\infty,2)}$ and deserving the name of a \emph{model for $(\infty,2)$-categories}. Given a model of $(\infty,2)$-categories $\mathscr M$, the extra structure that is required is an embedding
\begin{equation}
    \label{embeddinggaunt}
    \mathscr Gaunt_2\hookrightarrow \mathscr M\simeq\mathscr Cat_{(\infty,2)}
\end{equation}
of the $\infty$-category $\mathscr Gaunt_2$ of gaunt\footnote{A $2$-category is \emph{gaunt} or \emph{rigid}  if it has no non-identity $2$-isomorphisms and no non-identity $1$-equivalences (or equivalently no non-identity $1$- and $2$-isomorphisms).}
$2$-categories into the $\infty$-category $\mathscr M$. If we take e.g.~$\mathscr M\coloneqq[\cM]_\infty$, for $\cM$ any of the model categories of $(\infty,2)$-categories from \cref{allnerves} for which a nerve construction was described, then the embedding \eqref{embeddinggaunt} can be taken to be the restriction
\[\mathscr Gaunt_2\hookrightarrow\mathscr Cat_2\simeq[\twocat]_\infty\hookrightarrow[\cM]_\infty\simeq\mathscr Cat_{(\infty,2)}\]
of the functor from \eqref{genericadjunction}, for a suitably chosen equivalence $[\twocat]_\infty\simeq\mathscr Cat_2$. This could be seen by employing \cite[Lem.~10.2]{BarwickSchommerPries}.
\end{rmk}

Next, we address how the equivalent adjunctions $\mathscr Cat_2\rightleftarrows\mathscr Cat_{(\infty,2)}$ from \eqref{genericadjunction} relate to work by Gepner--Haugseng \cite[\S6]{GH}.

\begin{rmk}
In \cite[Prop.\ 6.1.7]{GH}, Gepner--Haugseng identify a universal property that relates the $\infty$-category of $2$-categories $\mathscr Cat_{2}$
and the $\infty$-category $\mathscr Cat_{(\infty,2)}$ of $(\infty,2)$-categories. More precisely, the former can be understood as a localization of the latter with respect to the class of maps
\[\Sigma^2\Lambda\coloneqq \{\Sigma^2 \Delta[0]\hookrightarrow\Sigma^2 S^k\ |\ k>0\},\]
where $S^k$ denotes the $k$-th sphere as an object of the $\infty$-category $\mathscr S$ of spaces, and
\begin{equation}
    \label{GHsuspension}
    \Sigma^2\colon\mathscr S\to\mathscr Cat_{(\infty,2)}
\end{equation}
implements a suitable $2$-fold $2$-point suspension, constructed in \cite[Def.~4.3.21]{GH}.
From this, one deduces the existence of an adjunction
\begin{equation}
\label{GHextended}
\mathscr Cat_{(\infty,2)}\rightleftarrows\mathscr L_{\Sigma^2\Lambda}\mathscr Cat_{(\infty,2)}
\simeq\mathscr Cat_2\end{equation}
with left adjoint being reflector and right adjoint being inclusion.
\end{rmk}

Our goal in \cref{GHcharacterization} is to prove that, for compatibly chosen equivalences of $\infty$-categories, the incarnation
\begin{equation}
  \label{adjunctionfrommodelcats}  
[c_*]_{\infty}\colon \mathscr{C}at_{(\infty,2)}\simeq [\vcat{\sset_{(\infty,1)}}]_\infty\rightleftarrows [\twocat]_\infty\simeq\mathscr{C}at_{2}\colon [{\bf N}_*]_\infty\end{equation}
 induces precisely the adjunction \eqref{GHextended}. This will show that, hence, all the adjunctions of $\infty$-categories \eqref{DiagramAdjunctions} have the correct universal property and do the right thing.

\begin{thm}
\label{MainTheorem2}
The adjunctions of $\infty$-categories \eqref{GHextended} and \eqref{adjunctionfrommodelcats} are equivalent.
\end{thm}

\begin{proof}[Outline of the proof]
The proof involves three steps.
\begin{itemize}[leftmargin=*]
    \item First, in \cref{suspensioncomparison}, we will discuss why the functor between model categories
    \[\Sigma^2\colon\sset_{(\infty,0)}\to\vcat{\sset_{(\infty,1)}}\]
    from \cref{Sigma2} induces at the level of underlying $\infty$-categories the functor from \eqref{GHsuspension}, where $\sset_{(\infty,0)}$ is the Kan-Quillen model structure.
    \item Then, in \cref{QE2} we will show that the Quillen pair 
\[c_*\colon \vcat{\sset_{(\infty,1)}}\rightleftarrows\twocat\colon {\bf N}_*\]
and the left Bousfield localization adjunction
\[\mathrm{Id}\colon\vcat{\sset_{(\infty,1)}}\rightleftarrows\cL_{\Sigma^2\Lambda}(\vcat{\sset_{(\infty,1)}})\colon\mathrm{Id}\]
induce equivalent adjunctions at the level of underlying $\infty$-categories.
\item Finally, in \cref{equivalencewithGH2} we use the previous two steps, as well as other results from the literature, to establish that the functor of $\infty$-categories \eqref{adjunctionfrommodelcats} is indeed equivalent to \eqref{GHextended}, as desired.\qedhere
\end{itemize}
\end{proof}

%%%
%%%

\section{Nerves in \texorpdfstring{$\Theta_2$}{Theta-2}-models}

We devote this section to briefly recalling the main globular models of $(\infty,2)$-categories, namely those based on Joyal's cell category $\Theta_2$, and the relevant nerve constructions.

We refer the reader to \cite{JoyalDisks} for the category $\Theta_2$, which is a full subcategory of $\twocat$ \cite{BergerCellular,MakkaiZawadowski}. The generic object of $\Theta_2$ is a $2$-category of the form $\theta=[i|j_1,\ldots,j_i]$ for $i\ge0$ and $j_k\ge0$ for $k=1,\ldots,i$. For example, the $2$-category $[4|2,0,3,1]$ is the $2$-category generated by the following data.
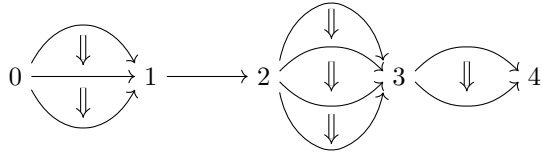
\begin{tz}
\node[](0) {$0$}; 
\node[right of=0,xshift=.3cm](1) {$1$}; 
\node[right of=1](2) {$2$}; 
\node[right of=2,xshift=.3cm](3) {$3$}; 
\node[right of=3,xshift=.3cm](4) {$4$}; 

\coordinate(a) at ($(0)!0.5!(1)+(0,.7cm)$);
\coordinate(b) at ($(0)!0.5!(1)-(0,.7cm)$);

\draw[->] (0) to node(c){} (1); 
\draw[->,bend left=25,rounded corners=1.2mm] ($(0.east)+(0,5pt)$) to (a) to ($(1.west)+(0,5pt)$);
\draw[->,bend right=25,rounded corners=1.2mm] ($(0.east)-(0,5pt)$) to (b) to ($(1.west)-(0,5pt)$);

\cell[la,right][n][0.5]{a}{c}{};
\cell[la,right][n][0.5]{c}{b}{};

\draw[->] (1) to (2);

\coordinate(a) at ($(2)!0.5!(3)+(0,1cm)$); 
\coordinate(b) at ($(2)!0.5!(3)+(0,.4cm)$); 
\coordinate(c) at ($(2)!0.5!(3)-(0,.4cm)$); 
\coordinate(d) at ($(2)!0.5!(3)-(0,1cm)$); 
 
\draw[->,bend left=30,rounded corners=1.2mm] ($(2.east)+(0,6pt)$) to (a) to ($(3.west)+(0,6pt)$); 
\draw[->,bend left=20,rounded corners=1mm] ($(2.east)+(0,2pt)$) to (b) to ($(3.west)+(0,2pt)$);
\draw[->,bend right=20,rounded corners=1mm] ($(2.east)-(0,2pt)$) to (c) to ($(3.west)-(0,2pt)$); 
\draw[->,bend right=30,rounded corners=1.2mm] ($(2.east)-(0,6pt)$) to (d) to  ($(3.west)-(0,6pt)$);

\cell[la,right][n][0.5]{a}{b}{};
\cell[la,right][n][0.5]{b}{c}{};
\cell[la,right][n][0.5]{c}{d}{};

\coordinate(b) at ($(3)!0.5!(4)+(0,.4cm)$); 
\coordinate(c) at ($(3)!0.5!(4)-(0,.4cm)$); 
 
\draw[->,bend left=20,rounded corners=1mm] ($(3.east)+(0,2pt)$) to (b) to ($(4.west)+(0,2pt)$);
\draw[->,bend right=20,rounded corners=1mm] ($(3.east)-(0,2pt)$) to (c) to ($(4.west)-(0,2pt)$); 

\cell[la,right][n][0.5]{b}{c}{};
\end{tz}

The canonical inclusion $\set\hookrightarrow\sset$ of sets as discrete simplicial sets induces a canonical inclusion $\psh{\Theta_2}\hookrightarrow\spsh{\Theta_2}$, which preserves limits and colimits. In particular, we often regard $\Theta_2$-sets as discrete $\Theta_2$-spaces without further specification.

For any object $\theta$ in $\Theta_2$, we denote by $\Theta_2[\theta]$ the $\Theta_2$-set represented by $\theta$ via the Yoneda embedding $\Theta_2\hookrightarrow\psh{\Theta_2}$.

\subsection{The models}
The following mathematical object was identified by Rezk \cite{rezkTheta}
as a model for $(\infty,2)$-categories. We recall the definition for completeness, but it will not be needed in this paper.

\begin{defn}
\label{thetaanodyne}
A \emph{complete Segal $\Theta_2$-space} is a $\Theta_2$-space $X\colon\Theta_2^{\op}\to\sset$ that is local\footnote{Given any small category $\cC$, there are well-defined derived mapping spaces $\Map^h_{\spsh{\cC}}(B,X)$ with respect to the homotopical structure on $\spsh{\cC}$ given by levelwise weak equivalences in $\sset_{(\infty,0)}$. For an explicit construction see e.g.~\cite[\S2.8]{rezkTheta}. We then say that a presheaf $X\colon \cC^{\op}\to\sset_{(\infty,0)}$ is \emph{local} with respect to a set of maps $S$ of $\spsh{\cC}$ if for every $f\colon A\to B$ in $S$ the induced map on derived mapping spaces $\Map^h_{\spsh{\cC}}(f,X)\colon\Map^h_{\spsh{\cC}}(B,X)\to\Map^h_{\spsh{\cC}}(A,X)$ with respect to levelwise weak homotopy equivalences is a weak equivalence of Kan complexes.}
with respect to the class of maps 
\begin{enumerate}[leftmargin=*]
    \item for all $i,j_1,\ldots,j_i\geq 0$, the \emph{horizontal Segality extension}
\[\Theta_2[1|j_1]\aamalg{\Theta_2[0]}\ldots\aamalg{\Theta_2[0]}\Theta_2[1|j_i]\to \Theta_2[i|j_1,\ldots,j_i]\]
    induced by the inclusions $[\langle s-1,s\rangle|j_s]\colon [1|j_s]\to [i|j_1,\ldots,j_i]$ where $\langle s-1,s\rangle\colon [1]\to [i]$ sends $0\mapsto s-1$ and $1\mapsto s$, for $1\leq s\leq i$;
    \item for all $j\geq 0$ the \emph{vertical Segality extension}
\[\Theta_2[1|1]\aamalg{\Theta_2[1|0]}\ldots\aamalg{\Theta_2[1|0]} \Theta_2[1|1]\to \Theta_2[1|j]\]
    induced by the inclusions $[1|\langle t-1,t\rangle]\colon [1|1]\to [1|j]$ where $\langle t-1,t\rangle\colon [1]\to [j]$ sends $0\mapsto t-1$ and $1\mapsto t$, for $1\leq t\leq j$;
    \item the \emph{horizontal completeness extension}
    \[\Theta_2[0]\to \Theta_2[0]\aamalg{\Theta_2[1|0]}\Theta[3|0,0,0]\aamalg{\Theta_2[1|0]}\Theta_2[0], \]
where the right-hand side is the colimit of the diagram
\begin{tz}
\node[](1) {$\Theta_2[0]$}; 
\node[right of=1,xshift=.5cm](2) {$\Theta_2[1|0]$}; 
\node[right of=2,xshift=1.2cm](3) {$\Theta_2[3|0,0,0]$}; 
\node[right of=3,xshift=1.2cm](4) {$\Theta_2[1|0]$}; 
\node[right of=4,xshift=.5cm](5) {$\Theta_2[0]$}; 

\draw[->] (2) to node[above,la]{$!$} (1);
\draw[->] (2) to node[above,la]{$\langle 0,2\rangle$} (3);
\draw[->] (4) to node[above,la]{$\langle 1,3 \rangle$} (3);
\draw[->] (4) to node[above,la]{$!$} (5);
\end{tz}
and the map is induced by the inclusion
$\langle 0\rangle\colon [0]\to [3]$;
    \item the \emph{vertical completeness extension}\footnote{While the completeness conditions are not the same as in \cite[\S11.4]{rezkTheta}, one can use \cite[\S4.4,\S10]{rezkTheta} to see that the two descriptions localizations are defining the same model structure.}
    \[\Theta_2[1|0]\to \Theta_2[1|0]\aamalg{\Theta_2[1|1]}\Theta[1|3]\aamalg{\Theta_2[1|1]}\Theta_2[1|0],\]
    where the right-hand side is the colimit of the diagram 
    \begin{tz}
\node[](1) {$\Theta_2[1|0]$}; 
\node[right of=1,xshift=.8cm](2) {$\Theta_2[1|1]$}; 
\node[right of=2,xshift=1.5cm](3) {$\Theta_2[1|3]$}; 
\node[right of=3,xshift=1.5cm](4) {$\Theta_2[1|0]$}; 
\node[right of=4,xshift=.8cm](5) {$\Theta_2[1|0]$}; 

\draw[->] (2) to node[above,la]{$[1|!]$} (1);
\draw[->] (2) to node[above,la]{$[1|\langle 0,2\rangle]$} (3);
\draw[->] (4) to node[above,la]{$[1|\langle 1,3 \rangle]$} (3);
\draw[->] (4) to node[above,la]{$[1|!]$} (5);
\end{tz}
and the map is induced by the inclusion
$[1|\langle 0\rangle]\colon [1|0]\to [1|3]$.
\end{enumerate}
\end{defn}

The following model structure is obtained as a left Bousfield localization of the injective model structure on $\spsh{\Theta_2}$.

\begin{thm}[{\cite[Thm~8.1]{rezkTheta}}]
\label{RezkModelStructure}
The category $\spsh{\Theta_2}$ of $\Theta_2$-spaces admits a model structure, denoted $\spsh{\Theta_2}_{(\infty,2)}$, in which
\begin{itemize}[leftmargin=*]
    \item the fibrant objects are the injectively fibrant complete Segal $\Theta_2$-spaces, and
    \item the cofibrations are the monomorphisms, and in particular every object is cofibrant.
\end{itemize}
\end{thm}

The following mathematical object was envisioned by Joyal \cite{JoyalDisks} and formalized by Ara \cite[\S5]{ara}
as a model for $(\infty,2)$-categories.

\begin{defn}
A \emph{$2$-quasi-category} is a $\Theta_2$-set $X\colon \Theta_2^{\op}\to\set$ that has the right lifting property with respect to the class of maps (1)-(4) from \cref{thetaanodyne}.
\end{defn}

\begin{thm}[{\cite[\S5.17]{ara}}]
\label{AraModelStructure}
The category $\psh{\Theta_2}$ of $\Theta_2$-sets admits a model structure, denoted $\psh{\Theta_2}_{(\infty,2)}$, in which
\begin{itemize}[leftmargin=*]
    \item the fibrant objects are the $2$-quasi-categories, and
    \item the cofibrations are the monomorphisms, and in particular every object is cofibrant.
\end{itemize}
\end{thm}

Ara showed as \cite[Thm~8.4]{ara} that the functor $(-)_0\colon\sset\to\set$, which extracts the $0$-th component,
induces a right Quillen equivalence
\begin{equation}
    \label{AraQE}
(-)_{\bullet,0}\colon\spsh{\Theta_2}_{(\infty,2)}\to \psh{\Theta_2}_{(\infty,2)}.
\end{equation}

%%%
%%%

\subsection{The nerve}

A nerve construction $\NLein\cD$ for any $2$-category $\cD$ was identified by Leinster \cite[Def.~J]{LeinsterSurvey} and further studied by Campbell\footnote{In the original source, the nerve is denoted $N\cD$, as opposed to $\NLein\cD$.} \cite{CampbellHoCoherent}.
Its construction is based on the notion of a normal pseudofunctor, which we recall later as \cref{normalpseudofunctor}.
Roughly speaking, those are maps between $2$-categories that preserve identities strictly and preserve compositions up to coherent
isomorphism.

\begin{const}[{\cite[Def.~J]{LeinsterSurvey}}]
\label{LeinsterNerve}
Let $\cD$ be a $2$-category. The nerve $\NLein\cD$ is the $\Theta_2$-set given for any $\theta\in\Theta_2$ by
the set of normal pseudofunctors from $\theta$ to $\cD$
\[ \NLein_\theta\cD\coloneqq(\NLein\cD)_\theta\coloneqq\twocat_{\nps}(\theta,\cD).\]
The assignment extends to a functor $\NLein\colon \twocat\to\psh{\Theta_2}$.
\end{const}

The homotopical properties of these nerve constructions follow from a combination of work by Campbell \cite{CampbellHoCoherent} and Lack \cite{lack2}, as explained in \cite[Rmk 5.16]{CampbellHoCoherent}.

\begin{thm}
\label{LeinsterEmbedding}
The functor $\NLein\colon \twocat\to\psh{\Theta_2}_{(\infty,2)}$ is a right Quillen embedding, and in particular a right Quillen and homotopical functor.
\end{thm}

\section{Nerves in bisimplicial models}

We devote this section to briefly recalling the main bisimplicial models of $(\infty,2)$-categories, and the relevant nerve constructions and model comparisons with the material from the previous section. We also prove the compatibility with the nerve construction from the previous section.

 The canonical inclusion $\set\hookrightarrow\sset$ of sets as discrete simplicial sets induces a canonical inclusion $\psh{(\Delta\times\Delta)}\hookrightarrow\spsh{(\Delta\times\Delta)}$, which preserves limits and colimits. In particular, we often regard bisimplicial sets as discrete bisimplicial spaces without further specification.

 For any $[i,j]$ in $\Delta\times\Delta$, we denote by $\Delta[i,j]$ the bisimplicial set represented by $[i,j]$ via the Yoneda embedding $\Delta\times\Delta\hookrightarrow\psh{(\Delta\times\Delta)}$.

\subsection{The first model}
The following mathematical object was identified by Barwick \cite[\S2.3]{BarwickThesis}
as a model of $(\infty,2)$-categories. It was also further studied with slightly different presentations, by Lurie \cite[Def.\ 1.3.6]{LurieGoodwillie},
Johnson-Freyd--Scheim\-bauer \cite[\S2]{JFSMorita} and Bergner--Rezk
\cite[Def.\ 5.3]{br2}.
See also \cite[\S2.2.2]{HaugsengThesis}
and \cite[\S4.2]{MoserNerve}. We recall the definition for completeness, but it will not be needed in this paper.

\begin{defn}
A \emph{$2$-fold complete Segal space} is a bisimplicial space $X\colon (\Delta\times\Delta)^{\op}\to\sset$ that is local with respect to the all maps of the following types:
\begin{enumerate}[leftmargin=*]
    \item\label{horSegCSS} for all $i,j\geq 0$, the \emph{horizontal Segality extension}
\[\Delta[1,j]\aamalg{\Delta[0,j]}\ldots \aamalg{\Delta[0,j]} \Delta[1,j]\to \Delta[i,j]\] 
    induced by the inclusions $\langle s-1,s\rangle\colon [1]\to [i]$ sending $0\mapsto s-1$, $1\mapsto s$, for $1\leq s\leq i$;
    \item for all $i,j\geq 0$, the \emph{vertical Segality extension}
\[\Delta[i,1]\aamalg{\Delta[i,0]}\ldots \aamalg{\Delta[i,0]} \Delta[i,1]\to \Delta[i,j]\]  
    induced by the inclusions $\langle t-1,t\rangle\colon [1]\to [j]$ sending $0\mapsto t-1$, $1\mapsto t$, for $1\leq t\leq j$; 
    \item for all $j\geq 0$ the \emph{horizontal completeness extension}
\[\Delta[0,j]\to \Delta[0,j]\aamalg{\Delta[1,j]}\Delta[3,j]\aamalg{\Delta[1,j]}\Delta[0,j], \]
    where the right-hand side is defined similarly to \cref{thetaanodyne}(3);
\item for all $i\geq 0$ the \emph{vertical completeness extension}
\[\Delta[i,0]\to \Delta[i,0]\aamalg{\Delta[i,1]}\Delta[i,3]\aamalg{\Delta[i,1]}\Delta[i,0],\]
where the right-hand side is defined similarly to \cref{thetaanodyne}(3);
    \item for all $j\geq 0$, the \emph{vertical homotopical constantness extension}
\[\Delta[0,0]\to \Delta[0,j]\]
    induced by the inclusion
    $\langle 0\rangle \colon [0]\to [j]$.
\end{enumerate}
\end{defn}

\begin{rmk}
\label{rmkSegal}
If $X$ is a $2$-fold complete Segal space, the fact that $X$ is local with respect to maps of the form \ref{horSegCSS} implies that for any $i,j\ge0$ the simplicial space $X_{\bullet,j}$ is local with respect to the map
\[\Delta[1]\aamalg{\Delta[0]}\ldots \aamalg{\Delta[0]} \Delta[1]\to \Delta[i] \]
induced by the inclusion $\langle s-1,s\rangle\colon [1]\to [i]$ for $1\leq s\leq i$. In particular, $X_{\bullet,j}$ is a Segal space. 
\end{rmk}

\begin{thm}[{\cite[Ch.~2]{BarwickThesis}}]
\label{BarwickModelStructure}
The category $\spsh{(\Delta\times\Delta)}$ of bisimplicial spaces admits a model structure, denoted $\spsh{(\Delta\times\Delta)}_{(\infty,2)}$, in which
\begin{itemize}[leftmargin=*]
    \item the fibrant objects are the injectively fibrant $2$-fold complete Segal spaces, and
    \item the cofibrations are the monomorphisms, and in particular every object is cofibrant.
\end{itemize}
\end{thm}

Precomposition with the functor $d\colon\Delta\times\Delta\to\Theta_2$, given on objects by $[i,j]\mapsto[i|j,\ldots,j]$, induces a functor $d^*\colon\spsh{\Theta_2}\to\spsh{(\Delta\times\Delta)}$. Explicitly,
\[(d^*X)_{i,j}\coloneqq X_{d[i,j]}=X_{[i|j,\ldots,j]}.\]
This functor admits a right adjoint $d_*$, which was proven by Bergner--Rezk as \cite[Cor.~7.1]{br2} to be a right Quillen equivalence
\begin{equation}
\label{d*QE}
d_*\colon\spsh{(\Delta\times\Delta)}_{(\infty,2)}\to\spsh{\Theta_2}_{(\infty,2)}.
\end{equation}

\subsection{The nerve}

A homotopically well-behaved nerve construction for bisimplicial models was studied by the first author in \cite{MoserNerve}, relying on the language of double categories. In this paper we aim at giving a presentation that is self-contained in the $2$-categorical world, so we take a slightly different viewpoint in recalling the necessary ingredients to describe the aforementioned nerve construction.

Three functors involving the category $\dblcat$ of double categories, namely the functors $L,L^{\simeq}\colon\dblcat\to\twocat$ and $\mathbb C\colon\spsh{\Delta\times \Delta}\to\dblcat$, are considered in \cite[\S\S2,5,6]{MoserNerve}. The composite functors $L\mathbb C,L^{\simeq}\mathbb C\colon\spsh{(\Delta\times\Delta)}\cong\psh{(\Delta\times\Delta\times\Delta)}\to \twocat$ are described more explicitly in \cite[Descr.~6.3.1, 6.3.2]{MoserNerve}. The following relation between the two functors is discussed in the proof of \cite[Thm\ 6.2.5]{MoserNerve}. 

\begin{prop} 
\label{LandLsimeq}
For $i,j,k\ge0$, there is a natural biequivalence of $2$-categories 
\[L^{\simeq}\bC \Delta[i,j,k] \to L\bC \Delta[i,j,k].\] 
\end{prop}

 The following nerve construction $\NMos\cD$ for any $2$-category $\cD$ was constructed by the first author\footnote{In the original source, the nerve is denoted $\mathbb N\mathbb H^{\simeq}\cD$, as opposed to $\NMos\cD$.} \cite[\S5.1]{MoserNerve} using the functor $L^{\simeq}\mathbb C$.

\begin{const}[{\cite[\S5.1]{MoserNerve}}]
\label{MoserNerve}
Let $\cD$ be a $2$-category. The nerve $\NMos\cD$ is the bisimplicial space given for any $i,j,k\ge0$ by
\[ \NMos_{i,j,k}\cD\coloneqq(\NMos\cD)_{i,j,k}\coloneqq \twocat(L^\simeq\bC\Delta[i,j,k], \cD).\]
The assignment extends to a functor $\NMos\colon \twocat\to\spsh{(\Delta\times\Delta)}$, which is the right adjoint to the functor $L^\simeq\bC\colon\spsh{(\Delta\times\Delta)}\to\twocat$.
\end{const}

\begin{thm}[{\cite[Thms~6.1.1,~6.1.3]{MoserNerve}}]
\label{MoserEmbedding}
The functor $\NMos\colon \twocat\to\spsh{(\Delta\times\Delta)}_{(\infty,2)}$ is a right Quillen embedding, and in particular a homotopical and right Quillen functor.
\end{thm}

Although the functor $L^{\simeq}\mathbb C$ is the one actually featuring in the definition of the nerve $\NMos$ from \cref{MoserNerve}, for the purpose of this paper it will be sufficient to have an explicit description of its easier version, the functor $L\mathbb C$. In order to give such a description, which is achieved in \cref{LCexplicit}, we first need to discuss preliminary material, including the $2$-categories $\Osim{i}$ and $\Otilde{k}$, and several flavors of tensor products between $2$-categories.

We denote by $\Sigma[1]=[1|1]$ the free-living $2$-cell, by $\Sigma\cI$ the free-living $2$-isomorphism, and by $\cO_2[i]$ the $i$-th $2$-truncated oriental; see e.g.~\cite[Def.\ 5.1.1]{MoserNerve} for an explicit description of this $2$-category.\footnote{As simplicial categories, there is an isomorphism $N_*\cO_2[i]\cong\mathfrak C\Delta[i]$, where $\mathfrak C\colon\sset\to\scat$ is the left adjoint to the homotopy coherent nerve functor.}

\begin{notn}[{\cite[Def.\ 5.1.1]{MoserNerve}}]
For $i\ge0$, let $\Osim{i}$ be the $2$-category obtained by gluing an invertible $2$-cell $\Sigma\cI$ on each generating $2$-cell $\Sigma[1]$ of the $2$-truncated $i$-oriental~$\cO_2[i]$; it can be expressed as the pushout of $2$-categories
  \begin{tz}
  \node[](1) {$\coprod\limits_{\genfrac{(}{)}{0pt}{2}{[i]}{[2]}}\Sigma[1]$}; 
 \node[below of=1](2) {$\coprod\limits_{\genfrac{(}{)}{0pt}{2}{[i]}{[2]}}\Sigma\cI$}; 
 \node[right of=1,xshift=.5cm,yshift=6pt](3) {$\cO_2[i]$}; 
 \node[below of=3](4) {$\Osim{i}$};
 
 \draw[->] ($(1.south)+(0,5pt)$) to (2);
 \draw[->] ($(1.east)+(0,6pt)$) to (3); 
 \draw[->] ($(2.east)+(0,6pt)$) to (4);
 \draw[->] (3) to (4);
 
\pushout{4};
 \end{tz}
 where the coproducts are indexed over the set of generating $2$-cells of $\cO_2[i]$.
\end{notn}

We denote by  
$[1]$ the free-living $1$-cell, and by $\cE$ the free adjoint equivalence; see e.g.~\cite[\S6]{lack2} or \cite[Not.\ 1.9]{Nerves2Cat} for an explicit description of this $2$-category.

\begin{notn}[{\cite[Def.\ 5.1.1]{MoserNerve}}]
For $k\ge0$, let $\Otilde{k}$ be the $2$-category obtained by gluing an adjoint equivalence $\cE$ on each generating $1$-cell $[1]$ of the $2$-category
$\Osim{k}$; it can be expressed as the pushout of $2$-categories
   \begin{tz}
 \node[](1) {$\coprod\limits_{\genfrac{(}{)}{0pt}{2}{[k]}{[1]}}[1]$}; 
 \node[below of=1](2) {$\coprod\limits_{\genfrac{(}{)}{0pt}{2}{[k]}{[1]}}\cE$}; 
 \node[right of=1,xshift=.5cm,yshift=6pt](3) {$\Osim{k}$}; 
 \node[below of=3,yshift=2pt](4) {$\Otilde{k}$};
 
 \draw[->] ($(1.south)+(0,2pt)$) to (2);
 \draw[->] ($(1.east)+(0,6pt)$) to (3); 
 \draw[->] ($(2.east)+(0,6pt)$) to ($(4.west)-(0,2pt)$);
 \draw[->] (3) to (4);
 
\pushout{4};
 \end{tz}
 where the coproducts are indexed over the set of generating $1$-cells of $\Osim{k}$.
\end{notn}

Next, we explore choices of tensor products for $2$-categories.

\begin{rmk}
\label{tensorproducts}
We consider several choices to form a $2$-category of $2$-functors between two $2$-categories $\cB$ and $\cD$, which all have $2$-functors from $\cB$ to $\cD$ as objects, and modifications as $2$-cells, but differ in the $1$-cells.
\begin{enumerate}[leftmargin=*]
\addtocounter{enumi}{-1}
\item The $2$-category $[\cB,\cD]$ consists of $2$-functors, (strict) natural transformations, and modifications.
    \item The $2$-category $[\cB,\cD]_{\textrm{ps}}$ consists of $2$-functors, pseudonatural transformations, and modifications.
    \item The $2$-category $[\cB,\cD]_{\textrm{lax}}$ consists of $2$-functors, lax\footnote{There are different conventions in the literature for the meaning of the word \emph{lax} (as opposed to \emph{oplax} or \emph{colax}), with equivalent resulting theories. The convention that we follow in this paper is consistent with the one used in e.g.\ in \cite{LackIcons, JFSMorita, HaugsengLax}, and it is opposite to the conventions of e.g.\ \cite{Gurski3D, AraLucas, AraMaltsiniotisJoin}.}  natural transformations, and modifications.
    \item The $2$-category $\iconI{\cB}{\cD}$ consists of $2$-functors, icons, and modifications. We recall that an \emph{icon} is a lax natural transformation for which each component is an identity.
\end{enumerate}
The first three notions were first discussed by Gray in \cite[\S{}I.2]{GrayFormal}, and the last one by Lack in \cite{LackIcons}. We refer the reader to the recent paper by Johnson--Yau \cite[Ch.~4]{JohnsonYau} for explicit definitions. The definition of the different kinds of natural transformations appears as \cite[Def.\ 4.2.1]{JohnsonYau}, where they call pseudonatural transformations ``strong'', the one of icons as \cite[Def.\ 4.6.2]{JohnsonYau}, and the one of modifications as \cite[Def.\ 4.4.1]{JohnsonYau}.

There are canonical and natural maps of $2$-categories
\[[\cB,\cD]\hookrightarrow[\cB,\cD]_{\textrm{ps}}\hookrightarrow[\cB,\cD]_{\textrm{lax}}\hookleftarrow[\cB,\cD]_{\textrm{ic}}.\]
Those constructions define functors
\[[-,-],[-,-]_{\textrm{ps}},[-,-]_{\textrm{lax}},[-,-]_{\textrm{ic}}\colon\twocat^{\op}\times\twocat\to\twocat\]

Each of those constructions is the internal hom functor for a corresponding tensor product which is
part of a two-variable adjunction $\twocat\times \twocat \to \twocat$, some of which are discussed in \cite[Thm~I.4.9, Thm~I.4.14, Cor.~I.4.17]{GrayFormal},
and \cite[Thm~3.16]{Gurski3D}.

The corresponding tensor products
\[\times,\pseudo,\otimes,\tensoricon\colon\twocat\times\twocat\to\twocat.\]
are, respectively:
\begin{enumerate}[leftmargin=*]
\setcounter{enumi}{-1}
\item\label{CartProd} the \emph{cartesian product} of $2$-categories $\cA\times\cB$;
    \item\label{PsProd} the \emph{pseudo Gray tensor product} of $2$-categories $\cA\pseudo\cB$;
    \item\label{LaxProd} the \emph{lax Gray tensor product} of $2$-categories $\cA\otimes\cB$; and
    \item\label{IconProd} a construction that we may call the \emph{icon tensor product} of $2$-categories $\cA\tensoricon\cB$.
\end{enumerate}
They are related via canonical and natural maps of $2$-categories
\[\cA\times\cB\leftarrow\cA\pseudo \cB\leftarrow\cA\otimes \cB\to\cA\tensoricon\cB.\]
Here, the left-pointing maps are classical (see e.g.~\cite[\S I.4.24]{GrayFormal} and \cite[Cor.\ 3.22]{Gurski3D}), 
and the right-pointing map is a consequence of \cref{tensoriconpushout}.

To highlight the difference between the four flavors, the four tensor products of the category $[1]$ with itself, or equivalently the corresponding naturality square, look as follows.
\begin{tz}[node distance=1cm]
\node[](0) {$\bullet$}; 
\node[below of=0](1) {$\bullet$}; 
\node[right of=0](2) {$\bullet$}; 
\node[below of=2](3) {$\bullet$}; 
\draw[->] (0) to (1);
\draw[->] (0) to (2); 
\draw[->] (1) to (3);
\draw[->] (2) to (3); 
\cell[la][e][0.5]{2}{1}{};
\node[la] at ($(1)!0.5!(3)-(0,.4cm)$) {$[1]\times [1]$};
\node[la] at ($(0)!0.5!(1)-(.6cm,0)$) {\eqref{CartProd}};

\node[right of=2,xshift=1cm](0) {$\bullet$}; 
\node[below of=0](1) {$\bullet$}; 
\node[right of=0](2) {$\bullet$}; 
\node[below of=2](3) {$\bullet$}; 
\draw[->] (0) to (1);
\draw[->] (0) to (2); 
\draw[->] (1) to (3);
\draw[->] (2) to (3); 
\cell[la,left,xshift=5pt,yshift=5pt][n][0.5]{2}{1}{$\cong$};
\node[la] at ($(1)!0.5!(3)-(0,.4cm)$) {$[1]\pseudo [1]$};
\node[la] at ($(0)!0.5!(1)-(.6cm,0)$) {\eqref{PsProd}};

\node[right of=2,xshift=1cm](0) {$\bullet$}; 
\node[below of=0](1) {$\bullet$}; 
\node[right of=0](2) {$\bullet$}; 
\node[below of=2](3) {$\bullet$}; 
\draw[->] (0) to (1);
\draw[->] (0) to (2); 
\draw[->] (1) to (3);
\draw[->] (2) to (3); 
\cell[la][n][0.5]{2}{1}{};
\node[la] at ($(1)!0.5!(3)-(0,.4cm)$) {$[1]\otimes [1]$};
\node[la] at ($(0)!0.5!(1)-(.6cm,0)$) {\eqref{LaxProd}};

\node[right of=2,xshift=1cm](0) {$\bullet$}; 
\node[below of=0](1) {$\bullet$}; 
\node[right of=0](2) {$\bullet$}; 
\node[below of=2](3) {$\bullet$}; 
\draw[d] (0) to (1);
\draw[->] (0) to (2); 
\draw[->] (1) to (3);
\draw[d] (2) to (3); 
\cell[la][n][0.5]{2}{1}{};
\node[la] at ($(1)!0.5!(3)-(0,.4cm)$) {$[1]\tensoricon [1]$};
\node[la] at ($(0)!0.5!(1)-(.6cm,0)$) {\eqref{IconProd}};
\end{tz}
\end{rmk}

Recall that the inclusion functor $\set\hookrightarrow\cat$ that regards every set as a discrete category admits left and right adjoint functors $\pi_0,\Ob\colon \cat\to\set$. They send a category $\cD$ to the set $\pi_0\cD$ of equivalence classes of its objects modulo the relation of being connected by a zigzag of $1$-morphisms, and to the set of objects~$\Ob\cD$ of $\cD$, respectively. The functor $\Ob\colon\cat\to\set$ also admits a right adjoint functor $\chaos\colon\set\to\cat$, which sends a set $S$ to the $1$-category $\chaos\,S$ whose set of objects is $S$ and which has exactly one morphism between any pair of objects.

These functors induce by base-change functors
$(\pi_0)_*,\Ob_*\colon \twocat\to \cat$ which are left and right adjoint to the inclusion functor $\cat\hookrightarrow\twocat$ that regards every category as a discrete $2$-category. They send a $2$-category~$\cD$ to the category $(\pi_0)_*\cD$ with the same objects as $\cD$ and hom-sets between two objects $c,d$ in $\cD$ given by $\pi_0(\cD(c,d))$, and to its underlying category~$\Ob_*\cD$ obtained by forgetting the $2$-morphisms, respectively.\footnote{The constructions $\Ob\cC$, $\pi_0\cC$, $\Ob_*\cD$ and $(\pi_0)_*\cD$ correspond to $\tau_{\leq0}^\mathrm{b}\cC$, $\tau_{\leq0}^\mathrm{i}\cC$, $\tau_{\leq1}^\mathrm{b}\cD$, and $\tau_{\leq1}^\mathrm{i}\cD$, respectively, following \cite[\S1.2]{AraMaltsiniotisJoin}, for a category $\cC$ and a $2$-category $\cD$.}
By base-change, we also get a right adjoint $\chaos_*\colon\cat\to\twocat$ for the functor $\Ob_*\colon\twocat\to\cat$, which sends a category $\cC$ to the $2$-category $\chaos_*\,\cC$ whose underlying category is $\cC$ and which has exactly one $2$-cell between any pair of parallel $1$-cells.

\begin{lem}
\label{tensoriconpushout}
For any $2$-categories $\cA$ and $\cB$ there is a pushout of $2$-categories
 \begin{tz}
 \node[](1) {$\Ob_*\cA \otimes \Ob \cB$}; 
 \node[below of=1](2) {$\pi_0(\Ob_*\cA) \otimes \Ob \cB$}; 
 \node[right of=1,xshift=1cm](3) {$\cA\otimes \cB$}; 
 \node[below of=3](4) {$\cA\tensoricon \cB$};
 \punctuation{4}{.};
 
 \draw[->] (1) to (2);
 \draw[->] (1) to (3); 
 \draw[->] (2) to (4);
 \draw[->] (3) to (4);
 
\pushout{4};
 \end{tz}
\end{lem}

\begin{proof}
The statement follows formally from the fact that for any $2$-categories $\cB$
 and $\cD$ the commutative square
 \begin{tz}
 \node[](1) {$\iconI{\cB}{\cD}$}; 
 \node[below of=1](2) {$\prod\limits_{\Ob\cB}\Ob\cD$}; 
 \node[right of=1,xshift=1.5cm](3) {$[\cB,\cD]_{\textrm{lax}}$}; 
 \node[below of=3](4) {$\prod\limits_{\Ob\cB}\chaos_* \Ob_*\cD$};
 
 \draw[->] (1) to (2);
 \draw[->] (1) to (3); 
 \draw[->] ($(2.east)+(0,4pt)$) to ($(4.west)+(0,4pt)$);
 \draw[->] (3) to (4);
 
\pullback{1};
 \end{tz}
is a pullback of $2$-categories.
\end{proof}

The following lemma can be understood as a special instance of \cite[Prop.~4.5]{AraLucas}.
\begin{lem}
\label{ReplacePseudo}
Given any $2$-category $\cA$ in which any $1$-morphism is an equivalence, and any $2$-category $\cB$, the canonical map is an isomorphism of $2$-categories
\[
\cA\otimes \cB \xrightarrow{\cong} \cA\pseudo \cB.
\]
\end{lem}

\begin{proof}
First, observe we have the following commutative diagram of $2$-categories.
\begin{tz}
 \node[](1) {$\coprod\limits_{\Mor\cA\times\Mor\cB}\Sigma[1]$}; 
 \node[below of=1](2) {$\coprod\limits_{\Mor\cA\times\Mor\cB}\Sigma \cI$}; 
 \node[right of=1,xshift=2cm](3) {$\coprod\limits_{\Mor\cA\times\Mor\cB}[1]\otimes[1]$}; 
 \node[below of=3](4) {$\coprod\limits_{\Mor\cA\times\Mor\cB}[1]\pseudo[1]$};
 \node[right of=3,xshift=1.8cm,yshift=5pt](5) {$\cA\otimes\cB$}; 
 \node[below of=5](6) {$\cA\pseudo\cB$};
 
 \draw[->] (1) to (2);
 \draw[->] ($(1.east)+(0,5pt)$) to ($(3.west)+(10pt,5pt)$); 
 \draw[->] ($(2.east)+(0,5pt)$) to ($(4.west)+(12pt,5pt)$); 
 \draw[->] (3) to (4);
 \draw[->] ($(3.east)+(0,5pt)$) to (5);
 \draw[->] ($(4.east)+(0,5pt)$) to (6);
 \draw[->] (5) to (6);
 
 \coordinate(a) at ($(4)+(0,5pt)$);
\pushout{a};
\pushout{6};
\end{tz}
Here, the external and left-hand commutative squares are pushouts, so the right-hand one is too. If $\cB$ meets the assumptions of the lemma, the map $\coprod [1]\otimes [1]\to \cA\otimes \cB$ factors through the canonical inclusion $\coprod [1]\otimes [1]\to \coprod \cE\otimes [1]$ and we obtain a lift in the above right-hand square, constructed as follows.
\begin{tz}
 \node[](1) {$\coprod\limits_{\Mor\cA\times\Mor\cB}[1]\otimes[1]$}; 
 \node[below of=1](2) {$\coprod\limits_{\Mor\cA\times\Mor\cB}[1]\pseudo [1]$}; 
 \node[right of=1,xshift=2.5cm](3) {$\coprod\limits_{\Mor\cA\times\Mor\cB}\cE\otimes [1]$}; 
 \node[below of=3](4) {$\coprod\limits_{\Mor\cA\times\Mor\cB}\cE\pseudo [1]$};
 \node[right of=3,xshift=1.8cm,yshift=5pt](5) {$\cA\otimes\cB$}; 
 \node[below of=5](6) {$\cA\pseudo\cB$};

 \draw[->] (1) to (2);
 \draw[->] ($(1.east)+(0,5pt)$) to ($(3.west)+(10pt,5pt)$); 
 \draw[->] ($(2.east)+(0,5pt)$) to ($(4.west)+(13pt,5pt)$); 
 \draw[->] (3) to (4);
 \draw[->,dashed] ($(3.east)+(0,5pt)$) to (5);
 \draw[->] ($(4.east)+(0,5pt)$) to (6);
 \draw[->] (5) to (6);
 \draw[->,bend left=18] (1) to (5);
 
 \draw[->,dashed] (2) to (3);
\end{tz}
It follows from the universal property of pushouts that in any pushout square that admits a diagonal lift the right vertical map is an isomorphism of $2$-categories, which concludes the proof.
\end{proof}

With the following proposition, we can now give an explicit description of the functor~$L$ on representable presheaves $\Delta[i,j,k]$. It could be taken as a definition by the reader who encounters it for the first time, or as a statement for the reader who is familiar with the double categorical framework, whose necessary ingredients we recall in the proof.

\begin{prop}
\label{LCexplicit}
For $i,j,k\ge0$ there is a natural isomorphism of $2$-categories
\[L\mathbb C\Delta[i,j,k]\cong
\Osim{j} \tensoricon (\Osim{i} \pseudo \Otilde{k}).
\]
\end{prop}

\begin{proof}
First, we recall the relevant constructions and definitions from \cite{MoserNerve} needed to prove the desired claim.
\begin{enumerate}[leftmargin=*, label=(Recall\ \arabic*)]
    \item The \emph{horizontal and vertical embeddings} $\mathbb H,\mathbb V\colon\twocat\to\dblcat$, which regard any $2$-category $\cD$ as a \emph{horizontal} and \emph{vertical double category}, recalled as \cite[Def.\ 2.1.7, Rmk\ 2.1.10]{MoserNerve}.
    \item\label{ItembfV} Their respective right adjoint functors $\mathbf{H},\mathbf{V}\colon \dblcat\to\twocat$, namely the \emph{underlying horizontal} and \emph{vertical $2$-category} functors, are discussed in \cite[Def.\ 2.1.8, Rmk\ 2.1.10]{MoserNerve}.
        \item The functor $\mathbb C\colon\psh{(\Delta\times\Delta\times\Delta)}\to\dblcat$ from \cite[Prop.\ 5.1.4]{MoserNerve}.
    \item\label{ItemL} The left adjoint functor $L\colon\dblcat\to\twocat$ of $\mathbb H$, discussed in \cite[\S6]{MSV1}.
     \item\label{ItemDblHom} The pseudo hom double category $\llbracket-,-\rrbracket\colon\dblcat^{\op}\times\dblcat\to\dblcat$ from \cite[\S2.2]{BoehmGray}.
    \item The corresponding pseudo Gray tensor product of double categories of \cite{BoehmGray}
    $\dblgray\colon\dblcat\times\dblcat\to\dblcat$.
\end{enumerate}
Next, we collect a few important facts that we will use.
\begin{enumerate}[leftmargin=*, label=(Obs.\ \arabic*)]
    \item\label{ObsIcon} For any $2$-categories $\cB$ and $\cD$, there is a natural isomorphism of $2$-categories
    \[\mathbf{V}\llbracket\bH\cB, \bH\cD\rrbracket\cong\iconI{\cB}{\cD}.\]
    This can be deduced from a careful analysis of the involved $2$-categories.
    \item\label{ObsX} For $i,j,k\ge0$, by \cite[Def.\ 2.2.4, Def.\ 5.1.3]{MoserNerve} the value of $\mathbb C$ at $\Delta[i,j,k]$ is given by the $2$-category
\[
\bX_{i,j,k}=\mathbb C\Delta[i,j,k]=\bV\Osim{j} \dblgray \bH\Osim{i} \dblgray\bH\Otilde{k}.
\]
\item\label{Lem79} By \cite[Lem.~7.8]{MSV1}, for any $2$-categories $\cA$ and $\cB$, there is an isomorphism of double categories
\[\bH\cA\dblgray\bH\cB\cong\bH(\cA \pseudo\cB).\]
\end{enumerate}
Now, for any $2$-category $\cD$ and $i,j,k\ge0$, we obtain natural bijections
\begin{align*}
\twocat(L\mathbb C\Delta[i,j,k], \cD) & \cong \dblcat(\mathbb C\Delta[i,j,k], \bH\cD) & \text{\ref{ItemL}}\\
 & \cong \dblcat(\bV\Osim{j} \dblgray \bH\Osim{i} \dblgray\bH\Otilde{k}, \bH\cD) & \text{\ref{ObsX}}\\
&\cong \dblcat(\bV\Osim{j} \dblgray \bH(\Osim{i} \pseudo\Otilde{k}), \bH\cD) &\text{\ref{Lem79}}\\
&\cong \dblcat(\bV\Osim{j}, \llbracket\bH(\Osim{i} \pseudo\Otilde{k}), \bH\cD\rrbracket)& \text{\ref{ItemDblHom}}\\
&\cong \twocat(\Osim{j}, \mathbf{V}\llbracket\bH(\Osim{i} \pseudo\Otilde{k}), \bH\cD\rrbracket)&\text{\ref{ItembfV}}\\
&\cong \twocat(\Osim{j}, \iconI{\Osim{i} \pseudo\Otilde{k}}{\cD})&\text{\ref{ObsIcon}}\\
&\cong \twocat(\Osim{j} \tensoricon (\Osim{i} \pseudo \Otilde{k}), \cD).&\text{Rmk~\ref{tensorproducts}}
\end{align*}
The claim then follows from the Yoneda lemma.
\end{proof}

%%%
%%%

\subsection{Nerve comparison}

We study the compatibility between $\NMos$ and $\NLein$.

\begin{lem}
\label{optimistic lemma}
For any $2$-category $\theta$ in $\Theta_2$ there is an isomorphism of bisimplicial sets
\[d^*\Theta_2[\theta]=d^*\NLein\theta\cong \NMos_{i,j,0}\theta.\]
\end{lem}

\begin{proof}
As a preliminary observation, there is a canonical map
\[ \Osim{j} \tensoricon \Osim{i}\to [j]\tensoricon [i]\cong [i|j,\ldots,j],\]
that can be seen by inspection to be a biequivalence.

Given that the $2$-category $\theta$ is gaunt, namely it does not have any non-identity $2$-isomorphisms and any non-identity $1$-equivalences, the functor $\twocat(-,\theta)$ sends biequivalences to bijections, and the functors $\twocat(-,\theta)$ and $\twocat_\mathrm{nps}(-,\theta)$ are isomorphic.
Hence, for any $i,j\geq 0$ we find natural bijections
\begin{align*}
    (d^*\NLein\theta)_{i,j} &=\NLein_{d[i,j]}\theta=\NLein_{[i|j,\ldots,j]}\theta \cong\twocat_\mathrm{nps}([i|j,\ldots,j],\theta) \\
    & \cong\twocat([i|j,\ldots,j],\theta)\cong\twocat(\Osim{j} \tensoricon \Osim{i} ,\theta)  \\
    &\cong\twocat(L\mathbb C\Delta[i,j,0],\theta) \cong \NMos_{i,j,0}\theta, 
\end{align*}
as desired. 
\end{proof}

\begin{thm}
\label{LeinsterVSMoser}
For any $2$-category $\cD$ there is a natural isomorphism of $\Theta_2$-sets
\[d_*\NMos_{\bullet,0}\cD\cong \NLein\cD.\]
\end{thm}

The crucial technical computation occurring in the proof is proven later as \cref{bisimplicialnerveFF}.

\begin{proof}[Proof of \cref{LeinsterVSMoser}]
For any $2$-category $\cD$ and any object $\theta$ in $\Theta_2$ there is a natural bijection
\begin{align*}
    d_*\NMos_{\theta,0} \cD &\cong \spsh{\Theta_2}(\Theta_2[\theta],d_*\NMos\cD) & \\
     &\cong \spsh{(\Delta\times\Delta)}(d^*\Theta_2[\theta],\NMos\cD) & d^*\dashv d_*\\
     &\cong \spsh{(\Delta\times\Delta)}(d^*\NLein\theta,\NMos\cD) & \text{\cref{optimistic lemma}} \\
     &\cong \spsh{(\Delta\times\Delta)}(\NMos\theta,\NMos\cD) & \text{\cref{optimistic lemma}}\\
&\cong \twocat_{\nps}(\theta,\cD) & \text{\cref{bisimplicialnerveFF}}\\
&\cong \NLein_\theta\cD, &
\end{align*}
as desired.
\end{proof}

Recall the right Quillen equivalences from \eqref{AraQE} and \eqref{d*QE} and the nerve constructions from \cref{LeinsterNerve,MoserNerve}.

\begin{cor}
\label{corLeinsterVSMoser}
The diagram of $\infty$-categories
\begin{tz}
\node[](1) {$[\twocat]_\infty$}; 
\node[below of=1](2) {$[\spsh{\Theta_2}_{(\infty,2)}]_\infty$}; 
\node[left of=2,xshift=-2cm](3) {$[\spsh{(\Delta\times \Delta)}_{(\infty,2)}]_\infty$}; 
\node[right of=2,xshift=2cm](4) {$[\psh{\Theta_2}_{(\infty,2)}]_\infty$}; 

\draw[->] (1) to node[above,la,xshift=-15pt]{$[\NMos]_\infty$} (3);
\draw[->] (1) to node[above,la,xshift=15pt]{$[\NLein]_\infty$} (4); 
\draw[->] (3) to node[below,la]{$[d_*]_\infty$} (2); 
\draw[->] (2) to node[below,la]{$[(-)_{\bullet,0}]_\infty$} (4);
\end{tz}
commutes up to equivalence.
\end{cor}

\begin{proof}
The corollary is an application of the ``right Quillen'' version of \cref{NerveRecognition} to the diagram
\begin{tz}
\node[](1) {$\twocat$}; 
\node[below of=1](2) {$\spsh{\Theta_2}_{(\infty,2)}$}; 
\node[left of=2,xshift=-1.5cm](3) {$\spsh{(\Delta\times \Delta)}_{(\infty,2)}$}; 
\node[right of=2,xshift=1.5cm](4) {$\psh{\Theta_2}_{(\infty,2)}$}; 
\punctuation{4}{.};

\draw[->] (1) to node[above,la,xshift=-15pt]{$\NMos$} (3);
\draw[->] (1) to node[above,la,xshift=15pt]{$\NLein$} (4); 
\draw[->] (3) to node[below,la]{$d_*$} (2); 
\draw[->] (2) to node[below,la]{$(-)_{\bullet,0}$} (4);
\end{tz}
The fact that all the assumptions of the lemma are met is from \cref{LeinsterVSMoser,MoserEmbedding,LeinsterEmbedding}.
\end{proof}

%%%
%%%
%%%

\subsection{A further model}
\label{PrecatSubsection}

Alternative models of $(\infty,2)$-categories, due to Bergner--Rezk, arises as the class of fibrant objects of two model structures on the category $P\cat(\spsh{\Delta})$, see \cite[\S6.7, \S6.11]{br1}. Here, $P\cat(\spsh{\Delta})$ denotes the full subcategory of $\spsh{(\Delta\times\Delta)}$ spanned by the bisimpicial spaces $X$ for which $X_0$ is a set. We refer to those as \emph{precategory objects in simplicial spaces}.\footnote{In the original source those are referred to as \emph{Segal precategories in simplicial spaces}.}
One is referred to as the \emph{injective-like} model structure on $P\cat(\spsh{\Delta})$, and one as the \emph{projective-like}. In this paper, we make use of the projective-like, which we denote $P\cat(\spsh{\Delta})_{(\infty,2)}$. We will never need an explicit description of this model structure, and we only use the fact that it comes with two Quillen equivalences, which will be recalled as \eqref{precatQE} and \eqref{PrecatQEenriched}.

The canonical inclusion functor $I\colon P\cat(\spsh{\Delta})\hookrightarrow\spsh{(\Delta\times\Delta)}$ admits a right adjoint~$R$, which was proven by Bergner--Rezk as \cite[Prop.~9.5, Thm~9.6]{br2} to be a right Quillen equivalence
\begin{equation}
    \label{precatQE}
R\colon \spsh{(\Delta\times\Delta)}_{(\infty,2)}\to P\cat(\spsh{\Delta})_{(\infty,2)}.
\end{equation}

The functor $I\colon P\cat(\spsh{\Delta})_{(\infty,2)}\to \spsh{(\Delta\times\Delta)}_{(\infty,2)}$ reflects weak equivalences between precategories,
in the sense of the following lemma.

\begin{lem}\label{Ireflects}
The functor $I\colon P\cat(\spsh{\Delta})_{(\infty,2)}\to \spsh{(\Delta\times\Delta)}_{(\infty,2)}$ reflects weak equivalences. That is,
if $f\colon X\to Y$ is a map in $P\cat(\spsh{\Delta})_{(\infty,2)}$ such that its image $If\colon IX\to IY$ is a weak equivalence in $\spsh{(\Delta\times\Delta)}_{(\infty,2)}$, then $f\colon X\to Y$ is a weak equivalence in $P\cat(\spsh{\Delta})_{(\infty,2)}$.
\end{lem}

This statement already occurs in the proof of \cite[Thm~9.6]{br2}, recalling from \cite[\S6.7,\ \S6.11]{br1} that the weak equivalences of the two model structures from \cite[\S6]{br1} considered on $P\cat(\spsh{\Delta})$ coincide. We recollect an outline of the argument here for the reader's convenience.

\begin{proof}
Bergner--Rezk introduce a functor $L\colon P\cat(\spsh{\Delta}) \to P\cat(\spsh{\Delta})$ and natural weak equivalence $X\xrightarrow{\simeq} LX$ in $\spsh{(\Delta\times\Delta)}_{(\infty,2)}$ for every $X$ in $P\cat(\spsh{\Delta})$, in \cite[\S6.7]{br1}.
By construction, $LX$ is an injectively fibrant Segal space.
It is discussed in \cite[\S6.7,\ \S6.11]{br1} that the functor $L$ detects weak equivalences of $P\cat(\spsh{\Delta})_{(\infty,2)}$ in the following sense: 
a map $f\colon X\to Y$ in $P\cat(\spsh{\Delta})$ is a weak equivalence in $P\cat(\spsh{\Delta})_{(\infty,2)}$ if and only if the induced map $Lf \colon LX\to LY$ is a Dwyer--Kan equivalence in the sense of \cite[Def.\ 8.2]{br2}.

Now assume a map $f\colon X\to Y$ in $P\cat(\spsh{\Delta})$ is a weak equivalence viewed as $If\colon IX\to IY$ in $\spsh{(\Delta\times\Delta)}_{(\infty,2)}$. We consider the commutative diagram 
\begin{tz}
 \node[](1) {$X$}; 
 \node[below of=1](2) {$LX$}; 
 \node[right of=1,xshift=.5cm](3) {$Y$}; 
 \node[below of=3](4) {$LY$};
 \punctuation{4}{.};
 
 \draw[->] (1) to (2);
 \draw[->] (1) to (3); 
 \draw[->] (2) to (4);
 \draw[->] (3) to (4);
\end{tz}
By assumption and by construction of $L$, all but possibly the lower horizontal map are weak equivalences in $\spsh{(\Delta\times\Delta)}_{(\infty,2)}$. By $2$-out-of-$3$, the lower horizontal map must also be. Now once again by construction, its source and target are (injectively fibrant) Segal spaces. By \cite[Thm~8.18]{br2}, this being a weak equivalence in  $\spsh{(\Delta\times\Delta)}_{(\infty,2)}$ is equivalent to being a Dwyer--Kan equivalence, thus showing that $f$ is a weak equivalence in $P\cat(\spsh{\Delta})_{(\infty,2)}$, as desired.
\end{proof}

From \cite[\S9]{br2}, for any bisimplicial space $X$, the value $RX$ of the right adjoint $R$ to the inclusion of $P\cat(\spsh{\Delta})$ into trisimplicial sets can be understood as the following pullback in trisimplicial sets.
\begin{diagram}
\label{descriptionR}
 \node[](1) {$RX$}; 
 \node[below of=1](2) {$X$}; 
 \node[right of=1,xshift=1cm](3) {$\cosk_0(X_{0,0,0})$}; 
 \node[below of=3](4) {$\cosk_0(X_{0,\bullet})$};
 
 \draw[->] (1) to (2);
 \draw[->] (1) to (3); 
 \draw[->] (2) to (4);
 \draw[->] (3) to (4);
 
\pullback{1};
\end{diagram}
Here, $\cosk_0\colon\spsh{\Delta}\to\spsh{\Delta\times \Delta}$ denotes the $0$-th coskeleton functor used in \cite[\S9]{br2}.
We give an alternative description of $RX$.

\begin{rmk}
\label{Descriptionq*}
For $i,j,k\ge0$, let $\DelQuot[i,j,k]$ denote the following pushout of bisimplicial spaces.
\begin{tz}
 \node[](1) {$\coprod\limits_{i+1}\Delta[0,j,k]$}; 
 \node[below of=1](2) {$\coprod\limits_{i+1} \Delta[0,0,0]$}; 
 \node[right of=1,xshift=1.2cm,yshift=4pt](3) {$\Delta[i,j,k]$}; 
 \node[below of=3,yshift=.5pt](4) {$\DelQuot[i,j,k]$};
 
 \draw[->] ($(1.south)+(0,6pt)$) to (2);
 \draw[->] ($(1.east)+(0,4pt)$) to (3); 
 \draw[->] ($(2.east)+(0,4.5pt)$) to (4);
 \draw[->] (3) to (4);
 
\pushout{4};
\end{tz}
Notice that, although not all objects occurring in the span belong to
$P\cat(\spsh{\Delta})$, the pushout $\DelQuot[i,j,k]$
does in fact belong to
$P\cat(\spsh{\Delta})$.

Given the pullback \eqref{descriptionR}, we deduce that for any bisimplicial space $X$ there is a natural bijection
\[(RX)_{i,j,k}\cong\spsh{(\Delta\times\Delta)}(\DelQuot[i,j,k],X).\]
 \end{rmk}

\begin{rmk}
\label{DelQuotZero}
For any $j,k\ge0$ there is an isomorphism of bisimplicial spaces
\[\DelQuot[0,j,k]\cong\Delta[0,0,0].\]
\end{rmk}

\subsection{The nerve}

In the remainder of this subsection, we study the bisimplicial space $R\NMos\cD$, which will be relevant in addressing the compatibility of the nerve constructions for bisimplicial and enriched models.

\begin{rmk}
Observe that for any $i,j,k\ge0$ and $\cD$ a $2$-category there is a natural bijection
\[(R\NMos\cD)_{i,j,k}\cong \spsh{(\Delta\times\Delta)}(\DelQuot[i,j,k], \NMos\cD)\cong
\twocat(L^{\simeq}\mathbb C\DelQuot[i,j,k],\cD).\]
\end{rmk}

\begin{prop} \label{lem:q_*NHvsq_*NHsim}
For any $i,j,k\ge0$ there is a natural isomorphism of $2$-categories
\[L^{\simeq}\mathbb C\DelQuot[i,j,k]\cong L\mathbb C\DelQuot[i,j,k].\]
In particular, for any $2$-category $\cD$ and $i,j,k\ge0$ there is a natural bijection
\[
(R\NMos\cD)_{i,j,k}\cong
\twocat(L^{\simeq}\mathbb C\DelQuot[i,j,k],\cD)\cong\twocat(L\mathbb C\DelQuot[i,j,k],\cD).\]
\end{prop}

The proof relies on the following lemma.

\begin{lem}
\label{lemLyne1}
For any $i,j,k\geq0$ there is a pushout of $2$-categories
\begin{tz}
 \node[](1) {$\coprod\limits_{i+1}\coprod\limits_{k+1} L^\simeq\bC\Delta[0,j,0]$}; 
 \node[below of=1](2) {$\coprod\limits_{i+1}\coprod\limits_{k+1}[0]$}; 
 \node[right of=1,xshift=2cm,yshift=4pt](3) {$L^{\simeq}\bC\Delta[i,j,k]$}; 
 \node[below of=3](4) {$L\bC\Delta[i,j,k]$};
 \punctuation{4}{.};
 
 \draw[->] ($(1.south)+(0,8pt)$) to (2);
 \draw[->] ($(1.east)+(0,4pt)$) to (3); 
 \draw[->] ($(2.east)+(0,4pt)$) to (4);
 \draw[->] (3) to (4);
 
\pushout{4};
\end{tz}
\end{lem}

\begin{proof}
Let $\cP$ be the pushout of the span of $2$-categories
\begin{tz}
\node[](1) {$\coprod\limits_{i+1}\coprod\limits_{k+1}[0]$}; 
\node[right of=1,xshift=1.8cm](2) {$\coprod\limits_{i+1}\coprod\limits_{k+1} L^\simeq\bC\Delta[0,j,0]$};
\node[right of=2,yshift=4pt,xshift=2cm](3) {$L^{\simeq}\bC\Delta[i,j,k]$}; 
\punctuation{3}{.};

\draw[->] ($(2.west)+(0,4pt)$) to ($(1.east)+(0,4pt)$);
\draw[->] ($(2.east)+(0,4pt)$) to (3);
\end{tz}
Using the naturality of the map in \cref{LandLsimeq}, we get an induced commutative diagram of $2$-categories.
\begin{tz}
 \node[](1) {$\coprod\limits_{i+1}\coprod\limits_{k+1} L^\simeq\bC\Delta[0,j,0]$}; 
 \node[below of=1](2) {$\coprod\limits_{i+1}\coprod\limits_{k+1}[0]$}; 
 \node[right of=1,xshift=2cm,yshift=4pt](3) {$L^{\simeq}\bC\Delta[i,j,k]$}; 
 \node[below of=3](4) {$\cP$};
 \node[below of=4,xshift=1.5cm,yshift=.3cm](5) {$L\bC\Delta[i,j,k]$};
 
 \draw[->] ($(1.south)+(0,8pt)$) to (2);
 \draw[->] ($(1.east)+(0,4pt)$) to (3); 
 \draw[->] ($(2.east)+(0,4pt)$) to (4);
 \draw[->] (3) to (4);
 
 \draw[->,dashed] (4) to (5);
 \draw[->,bend right=15] ($(2.south)+(18pt,8pt)$) to (5.west);
 \draw[->,bend left=25] (3) to (5);
 
\pushout{4};
\end{tz}
In the diagram, the left vertical map is a coproduct of the biequivalence
\[L^\simeq \bC \Delta[0,j,0]\to L \bC \Delta[0,j,0] \cong
\Osim{j} \tensoricon (\Osim{0} \pseudo \Otilde{0})\cong[0],\]
built using
\cref{LandLsimeq,LCexplicit}, so it is a biequivalence itself.

Also, the top horizontal map is obtained by applying the composite left Quillen functor $L^{\simeq}\bC\colon \spsh{(\Delta\times\Delta)}\to\twocat$ from \cref{MoserEmbedding} to the cofibration
\[\coprod\limits_{i+1}\coprod\limits_{k+1}\Delta[0,j,0]\to\Delta[i,j,k],\]
so it is a cofibration itself.

Since the model structure on $\twocat$ is left proper by \cite[Thm\ 6.3]{lack1} (see also \cite[\S2]{lack2}), it follows that the bottom horizontal map
\[L^{\simeq}\bC\Delta[i,j,k]\to\cP\]
is also a biequivalence.
Since the map $L^\simeq\bC\Delta[i,j,k]\to L\bC\Delta[i,j,k]$ is a biequivalence by \cref{LandLsimeq}, then by $2$-out-of-$3$ the comparison map
\[\cP\to L\bC\Delta[i,j,k]\]
is also a biequivalence. Using the explicit description from \cite[Desc.\ 6.3.1]{MoserNerve},
we now see that this comparison map is an isomorphism on underlying $1$-categories, which is sufficient to conclude that it must in fact be an isomorphism of $2$-categories, as biequivalences are in particular isomorphisms on $2$-morphisms.
\end{proof}

We can now prove the proposition.

\begin{proof}[Proof of \cref{lem:q_*NHvsq_*NHsim}]
We argue that for any $i,j,k\ge0$ there is an isomorphism of $2$-categories 
\begin{equation}
\label{equation}
L^{\simeq}\bC \DelQuot[i,j,k]\cong L\bC \DelQuot[i,j,k]
\end{equation}
that is natural in $i,j,k$. 

To this end, we consider the following commutative diagram in of $2$-categories.
\begin{tz}
\node[](1) {$\coprod\limits_{i+1}\coprod\limits_{k+1} [0]$}; 
\node[right of=1,xshift=2.5cm](1') {$\coprod\limits_{i+1}\coprod\limits_{k+1} [0]$}; 
\node[right of=1',xshift=2.5cm](1'') {$\coprod\limits_{i+1}\coprod\limits_{k+1}[0]$};

\node[below of=1](2) {$\coprod\limits_{i+1}\coprod\limits_{k+1} [0]$};
\node[below of=1'](2') {$\coprod\limits_{i+1}\coprod\limits_{k+1} L^\simeq\bC\Delta[0,j,0]$}; 
\node[below of=1''](2'') {$\coprod\limits_{i+1}\coprod\limits_{k+1} L^\simeq\bC\Delta[0,j,0]$};

\node[below of=2](3) {$\coprod\limits_{i+1} [0]$};
\node[below of=2'](3') {$\coprod\limits_{i+1}L^\simeq\bC\Delta[0,j,k]$};
\node[below of=2'',yshift=4pt](3'') {$L^\simeq\bC\Delta[i,j,k]$};

\draw[d] (1) to (2);
\draw[d] ($(1.east)+(0,4pt)$) to ($(1'.west)+(0,4pt)$);
\draw[->] (2) to (3);
\draw[d] ($(1'.east)+(0,4pt)$) to ($(1''.west)+(0,4pt)$);
\draw[->] ($(2'.west)+(0,4pt)$) to ($(2.east)+(0,4pt)$);
\draw[d] ($(2'.east)+(0,4pt)$) to ($(2''.west)+(0,4pt)$);
\draw[->] (2'') to (1'');
\draw[->] ($(2''.south)+(0,8pt)$) to (3''); 
\draw[->] (2') to (1'); 
\draw[->] ($(2'.south)+(0,8pt)$) to (3');
\draw[->] ($(3'.west)+(0,4pt)$) to ($(3.east)+(0,4pt)$);
\draw[->] ($(3'.east)+(0,4pt)$) to (3'');
\end{tz}
The colimit of this diagram can be equivalently computed by either taking the colimit of the colimit of each row, or by taking the colimit of the colimit of each column.

    On the one hand, by \cref{Descriptionq*} and using the fact that $L^{\simeq}$ is a left adjoint functor, the colimit of each row produces the span of $2$-categories
\begin{tz}
\node[](1) {$\coprod\limits_{i+1}\coprod\limits_{k+1} [0]$}; 
\node[below of=1](2) {$\coprod\limits_{i+1}\coprod\limits_{k+1} [0]$};
\node[below of=2](3) {$L^{\simeq}\bC \DelQuot[i,j,k]$};

\draw[d] (1) to (2);
\draw[->] (2) to (3);
\end{tz}
% \[\begin{tikzcd}
%     \coprod\limits_{i+1}\coprod\limits_{k+1}[0]\\
%     \coprod\limits_{i+1}\coprod\limits_{k+1}[0]\arrow[u,equal]\arrow[d]\\
%     L^{\simeq}\bC \DelQuot[i,j,k],
% \end{tikzcd}\]
whose pushout is $L^{\simeq}\bC \DelQuot[i,j,k]$.

    On the other hand, by \cref{lemLyne1} we see that the pushout of each column produces the following span of $2$-categories
\begin{tz}
\node[](1) {$\coprod\limits_{i+1} [0]$}; 
\node[right of=1,xshift=2.5cm](1') {$\coprod\limits_{i+1}L\bC\Delta[0,j,k]$}; 
\node[right of=1',xshift=2.5cm,yshift=4.5pt](1'') {$L\bC\Delta[i,j,k]$}; 
\draw[->] ($(1'.west)+(0,4.5pt)$) to ($(1.east)+(0,4.5pt)$);
\draw[->] ($(1'.east)+(0,4.5pt)$) to (1'');
\end{tz}
and by \cref{Descriptionq*} and using the fact that $L$ and $\mathbb C$ are left adjoint functors, its colimit is $L\bC \DelQuot[i,j,k]$.

Hence, the isomorphism \eqref{equation} follows. 
\end{proof}

\begin{rmk}
\label{SegalObjectRmk1}
Given any $2$-category $\cD$, combining \cref{MoserEmbedding} and \eqref{precatQE} we know that $R\NMos\cD$ is fibrant in $P\cat(\spsh{\Delta})_{(\infty,2)}$.
By \cref{rmkSegal}, for any $j\ge0$ we know that $(R\NMos\cD)_{\bullet,j}$ is a Segal space. It follows that for any $i,j\ge0$ we have a weak equivalence of spaces
\[(R\NMos\cD)_{i,j}\simeq(R\NMos\cD)_{1,j}\ttimes{(R\NMos\cD)_{0,j}}\ldots\ttimes{(R\NMos\cD)_{0,j}}(R\NMos\cD)_{1,j}.\]
\end{rmk}

This motivates us to understand better the sets $(R\NMos\cD)_{0,j,k}$ and $(R\NMos\cD)_{1,j,k}$, which we achieve in \cref{RofNLyne0,RofNLyne1}.

\begin{prop}
\label{RofNLyne0}
For any $2$-category $\cD$ and $j,k\ge0$ there is a natural bijection
\[(R\NMos\cD)_{0,j,k}\cong\Ob\cD.\]
\end{prop}

\begin{proof}
By \cref{Descriptionq*,DelQuotZero}, for any $j,k\ge0$, we have a natural bijection
\begin{align*}
    (R\NMos\cD)_{0,j,k} &\cong \spsh{(\Delta\times\Delta)}(\DelQuot[0,j,k],\NMos\cD)\\
     & \cong \spsh{(\Delta\times \Delta)} (\Delta[0,0,0], \NMos\cD)\\
    & \cong \NMos_{0,0,0}\cD \cong \twocat ([0], \cD)\cong\Ob\cD,
\end{align*}
as desired.
\end{proof}

We now proceed to describing $(R\NMos\cD)_{1,j,k}$, which requires some extra work.

Given any category $\cA$, we denote by $\Sigma\cA$ the $2$-point suspension of $\cA$, which consists of two distinct objects and a single interesting hom-category given by $\cA$. The construction extends to a left adjoint functor $\Sigma\colon\cat\to\twocat_{*,*}$.

\begin{lem}
\label{lemmasuspensionhigher}
For any $2$-category $\cA$ there is a pushout of $2$-categories
 \begin{tz}
 \node[](1) {$\cA\amalg\cA$}; 
 \node[below of=1](2) {$\cA\otimes [1]$}; 
 \node[right of=1,xshift=1cm](3) {$[0]\amalg[0]$}; 
 \node[below of=3](4) {$\Sigma(\pi_0)_*\cA$};
 \punctuation{4}{.};
 
 \draw[->] (1) to (2);
 \draw[->] (1) to (3); 
 \draw[->] (2) to (4);
 \draw[->] (3) to (4);
 
\pushout{4};
 \end{tz} 
\end{lem}

\begin{proof}
If we denote by $\cP\cA$ the following pushout of $2$-categories,
\begin{tz}
 \node[](1) {$\cA\amalg\cA$}; 
 \node[below of=1](2) {$\cA\otimes [1]$}; 
 \node[right of=1,xshift=1cm](3) {$[0]\amalg[0]$}; 
 \node[below of=3](4) {$\cP\cA$};
 
 \draw[->] (1) to (2);
 \draw[->] (1) to (3); 
 \draw[->] (2) to (4);
 \draw[->] (3) to (4);
 
\pushout{4};
 \end{tz}
this construction can also be regarded as a left adjoint functor $\cP\colon\twocat\to\twocat_{*,*}$. At the same time, also $\Sigma(\pi_0)_*$ defines a left adjoint functor $\Sigma(\pi_0)_*\colon\twocat\to\twocat_{*,*}$.
One can now prove by direct inspection that for any $i$-cell $\Sigma^i[1]$ for $i=0,1,2$, there is a natural isomorphism of bipointed $2$-categories
\[\cP\Sigma^i[1]\cong \Sigma(\pi_0)_*\Sigma^i[1].\]
It follows by cocontinuity that for every $2$-category $\cA$ there is an isomorphism of (bipointed) $2$-categories
\[\cP\cA\cong \Sigma(\pi_0)_*\cA,\]
concluding the proof. 
\end{proof}

For $j,k\ge0$, we let $\widetilde{[k]}$ denote the unique contractible groupoid with $k+1$ objects, namely the category with $k+1$ objects and a unique morphism between any two objects,
and $\Sigma([j]\times \widetilde{[k]})$ the $2$-point suspension of the $1$-category $[j]\times \widetilde{[k]}$. This is the $2$-category with two objects and a single interesting hom-category given by $[j]\times \widetilde{[k]}$.

\begin{lem}
\label{lemLyne2}
For any $j,k\ge0$ there is a pushout of $2$-categories
\begin{tz}
 \node[](1) {$L\bC\Delta[0,j,k]\amalg L\bC\Delta[0,j,k]$}; 
 \node[below of=1](2) {$L\bC\Delta[1,j,k]$}; 
 \node[right of=1,xshift=2cm](3) {$[0]\amalg[0]$}; 
 \node[below of=3,yshift=1pt](4) {$\Sigma([j]\times \widetilde{[k]})$};
 \punctuation{4}{.};
 
 \draw[->] (1) to (2);
 \draw[->] (1) to (3); 
 \draw[->] (2) to ($(4.west)-(0,1pt)$);
 \draw[->] (3) to ($(4.north)-(0,3pt)$);
 
\pushout{4};
 \end{tz}
\end{lem}

\begin{proof} 
Denote by $\cP$ the following pushout of $2$-categories.
\begin{tz}
 \node[](1) {$L\bC\Delta[0,j,k]\amalg L\bC\Delta[0,j,k]$}; 
 \node[below of=1](2) {$L\bC\Delta[1,j,k]$}; 
 \node[right of=1,xshift=2cm](3) {$[0]\amalg[0]$}; 
 \node[below of=3](4) {$\cP$};
 
 \draw[->] (1) to (2);
 \draw[->] (1) to (3); 
 \draw[->] (2) to (4);
 \draw[->] (3) to (4);
 
\pushout{4};
 \end{tz} 
 Consider the following commutative diagram of $2$-categories.
 \begin{tz}
\node[](1) {$\pi_0 (\Ob_*\Osim{j})\otimes \Ob([1]\pseudo \Otilde{k})$}; 
\node[right of=1,xshift=4.5cm,yshift=-4pt](1') {$\coprod\limits_{2}\pi_0 (\Ob_* \Osim{j})\otimes \Ob\Otilde{k}$}; 
\node[right of=1',xshift=2.2cm,yshift=-2pt](1'') {$\coprod\limits_{2} [0]$};

\node[below of=1](2) {$\Ob_* \Osim{j}\otimes \Ob([1]\pseudo \Otilde{k})$};
\node[below of=1'](2') {$\coprod\limits_{2}\Ob_* \Osim{j}\otimes \Ob\Otilde{k}$}; 
\node[below of=1''](2'') {$\coprod\limits_{2} [0]$};

\node[below of=2](3) {$\Osim{j}\otimes ([1]\pseudo \Otilde{k})$};
\node[below of=2'](3') {$\coprod\limits_{2}\Osim{j}\otimes \Otilde{k}$};
\node[below of=2''](3'') {$\coprod\limits_{2}[0]$};

\draw[->] ($(2.north)-(0,3pt)$) to (1);
\draw[->] ($(1'.west)+(0,2pt)$) to node[above,la]{$\cong$} ($(1.east)-(0,2pt)$);
\draw[->] ($(1'.east)+(0,2pt)$) to ($(1''.west)+(0,4pt)$);
\draw[->] (2) to ($(3.north)-(0,3pt)$);
\draw[->] ($(2'.west)+(0,2pt)$) to node[above,la]{$\cong$} ($(2.east)-(0,2pt)$);
\draw[->] ($(2'.east)+(0,2pt)$) to ($(2''.west)+(0,4pt)$);
\draw[d] (2'') to ($(1''.south)+(0,2pt)$);
\draw[d] ($(2''.south)+(0,2pt)$) to (3''); 
\draw[->] ($(2'.north)-(0,3pt)$) to ($(1'.south)+(0,8pt)$); 
\draw[->] ($(2'.south)+(0,8pt)$) to ($(3'.north)-(0,3pt)$);
\draw[->] ($(3'.west)+(0,2pt)$) to ($(3.east)-(0,2pt)$);
\draw[->] ($(3'.east)+(0,2pt)$) to ($(3''.west)+(0,4pt)$);
\end{tz}
The colimit of this diagram can be equivalently computed by either taking the colimit of the colimit of each row, or by taking the colimit of the colimit of each column.

By doing pushouts of each column first, we get using \cref{tensoriconpushout} the pushout of the span
    \begin{tz}
\node[](1) {$\Osim{j}\tensoricon ([1]\pseudo \Otilde{k})$}; 
\node[right of=1,xshift=3.5cm,yshift=-4pt](1') {$\coprod\limits_{2}\Osim{j}\tensoricon \Otilde{k}$}; 
\node[right of=1',xshift=2cm,yshift=-2pt](1'') {$\coprod\limits_{2}[0]$};
\punctuation{1''}{.};

\draw[->] ($(1'.west)+(0,2pt)$) to ($(1.east)-(0,2pt)$);
\draw[->] ($(1'.east)+(0,2pt)$) to ($(1''.west)+(0,4pt)$);
\end{tz}
By \cref{LCexplicit}, we identify the pushout to be computed as the pushout of the span
\begin{tz}
\node[](1) {$L\bC\Delta[1,j,k]$}; 
\node[right of=1,xshift=2cm,yshift=-4pt](1') {$\coprod\limits_{2}L\bC\Delta[0,j,k]$}; 
\node[right of=1',xshift=1.5cm](1'') {$\coprod\limits_{2}[0]$};

\draw[->] ($(1'.west)+(0,4pt)$) to (1);
\draw[->] ($(1'.east)+(0,4pt)$) to ($(1''.west)+(0,4pt)$);
\end{tz}
which gives precisely $\cP$. 

Now note that there are natural isomorphisms of $2$-categories
   \begin{align*}
        \Osim{j}\otimes ([1]\pseudo \Otilde{k}) &\cong \Osim{j}\otimes (\Otilde{k}\pseudo [1]) & \text{symmetry of } \pseudo \\
        &\cong \Osim{j}\otimes (\Otilde{k}\otimes [1]) & \text{\cref{ReplacePseudo}} \\
        &\cong (\Osim{j}\otimes \Otilde{k})\otimes [1]. & \text{associativity of } \otimes 
    \end{align*}
    By doing pushouts of each row, we get using the above isomorphism and \cref{lemmasuspensionhigher} applied to $\cA=\Osim{j}\otimes \Otilde{k}$ the pushout of the span 
    \begin{tz}
\node[right of=1',xshift=2.7cm,yshift=-2pt](1'') {$\coprod\limits_{i+1}[0]$};
\node[below of=1''](2'') {$\coprod\limits_{i+1}[0]$};
\node[below of=2''](3'') {$\Sigma (\pi_0)_* (\Osim{j}\otimes \Otilde{k})$};

\draw[d] (2'') to ($(1''.south)+(0,2pt)$);
\draw[->] ($(2''.south)+(0,2pt)$) to (3'');
\end{tz}
    which gives precisely $\Sigma (\pi_0)_* (\Osim{j}\otimes \Otilde{k})$. Combining \cite[Prop.\ A.27,~\S A.31]{AraMaltsiniotisJoin}, remembering that $(\pi_0)_*\cD\cong\tau_{\leq 1}^i\cD$, we obtain that for any $j,k\ge0$, there are natural isomorphisms of $2$-categories
\begin{align*}
    \Sigma(\pi_0)_*(\Osim{j} \otimes \Otilde{k})  &  \cong \Sigma(((\pi_0)_*\Osim{j})\times ((\pi_0)_*\Otilde{k}))\\
          &  \cong \Sigma([j]\times\widetilde{[k]}).
\end{align*}
    
So the desired isomorphism follows.
\end{proof}

%%%

%%%

We can now describe $(R\NMos\cD)_{1,j,k}$.

\begin{prop}
\label{RofNLyne1}
For any $j,k\ge0$ there is a natural
isomorphism of $2$-categories
\[L\bC\DelQuot[1,j,k]\cong\Sigma([j]\times \widetilde{[k]}).\]
In particular, for any $2$-category $\cD$ and $j,k\ge0$ there is a natural bijection
\[
(R\NMos\cD)_{1,j,k}\cong
\twocat(L^{\simeq}\mathbb C\DelQuot[1,j,k],\cD)\cong\twocat(\Sigma([j]\times \widetilde{[k]}),\cD).\]
\end{prop}

\begin{proof}
We show that for $j,k\ge0$ there is a natural isomorphism of $2$-categories
\begin{equation}
\label{equationsuspension}
L\bC\DelQuot[1,j,k]\cong\Sigma([j]\times \widetilde{[k]}).
\end{equation}
By \cref{Descriptionq*}, we know that $\DelQuot[1,j,k]$ is the pushout of the span
 \begin{tz}
\node[](1) {$\Delta[0,0,0]\amalg\Delta[0,0,0]$}; 
\node[right of=1,xshift=2.5cm](1') {$\Delta[0,j,k]\amalg\Delta[0,j,k]$}; 
\node[right of=1',xshift=1.6cm](1'') {$\Delta[1,j,k]$}; 
\punctuation{1''}{.};
\draw[->] (1') to (1);
\draw[->] (1') to (1'');
\end{tz}
% \[
% \begin{tikzcd}
%      {\Delta[0,0,0]\amalg\Delta[0,0,0]}& \arrow[l]\Delta[0,j,k]\amalg\Delta[0,j,k]\arrow[r]&\Delta[1,j,k]
% \end{tikzcd}
% \]
Since $L$ and $\bC$ are left adjoint functors, we obtain that $L\bC \DelQuot[1,j,k]$ is the pushout of the span
 \begin{tz}
\node[](1) {$[0]\amalg [0]$}; 
\node[right of=1,xshift=2.1cm](1') {$L\bC\Delta[0,j,k]\amalg L\bC\Delta[0,j,k]$}; 
\node[right of=1',xshift=2.4cm](1'') {$L\bC\Delta[1,j,k]$}; 
\punctuation{1''}{.};
\draw[->] (1') to (1);
\draw[->] (1') to (1'');
\end{tz}
% \[
% \begin{tikzcd}
%      {[0]\amalg [0]}& \arrow[l]L\bC\Delta[0,j,k]\amalg L\bC\Delta[0,j,k]\arrow[r]&L\bC\Delta[1,j,k]
% \end{tikzcd}
% \]
By \cref{lemLyne2}, its pushout is $\Sigma([j]\times \widetilde{[k]})$.
Hence, the isomorphism \eqref{equationsuspension} follows.

The second part of the statement is a consequence of the above isomorphism and  \cref{lem:q_*NHvsq_*NHsim}.
\end{proof}

\section{Nerves in categories enriched over \texorpdfstring{$(\infty,1)$}{(infinity,1)}-categories}

We refer the reader to \cite[Def.\ A.3.2.16]{htt} for the definition of an excellent monoidal model category. The following cases are relevant in this paper.
\begin{enumerate}[leftmargin=*]
\setcounter{enumi}{-1}
    \item Let $\cV=\cat$ be the canonical model structure on the category $\cat$ of small categories (see e.g.~\cite{RezkCat}), which is seen to be excellent using the fact that the ordinary nerve functor ${\bf N}\colon\cat\to\sset_{(\infty,1)}$ creates weak equivalences and commutes with filtered colimits.
    \item Let $\cV=\sset_{(\infty,1)}$ be the Joyal model structure on the category $\sset$ of simplicial sets from \cite[Thm~6.12]{JoyalVolumeII}, which is excellent by \cite[Ex.\ A.3.2.23]{htt}.
    \item Let $\cV=\spsh{\Delta}_{(\infty,1)}$ being the Rezk model structure from
\cite[Thm~7.2]{rezk}
on the category $\spsh{\Delta}$ of simplicial spaces, which is discussed to be excellent in \cite[Thm 3.11]{br1}.
    \item Let $\cV=\sset^+_{(\infty,1)}$ be the Lurie model structure on the category $\sset^+$ of marked simplicial sets from \cite[Prop.~3.1.3.7]{htt}, which is excellent by \cite[Ex.\ A.3.2.22]{htt}.
\end{enumerate}

\subsection{The models}

All enriched models of $(\infty,2)$-categories will be a special case of the following.

\begin{defn}
Let $\cV$ be an excellent monoidal model category. A \emph{locally fibrant $\cV$-category} is a $\cV$-category $\cD$ for which for any pair of objects $c,d$ in $\cD$ the hom-object $\cD(c,d)$ is fibrant in $\cV$.
\end{defn}

\begin{thm}[{\cite[Thm A.3.2.24]{htt}}] \label{EnrichedCatMS}
Let $\cV$ be an excellent monoidal model category. The category of small
categories enriched over $\cV$ admits a model structure in which
\begin{itemize}[leftmargin=*]
    \item the fibrant objects are the \emph{locally fibrant $\cV$-categories}, and 
\item the trivial fibrations are precisely the $\cV$-functors that are surjective on objects, and locally a trivial fibration in $\cV$.
\end{itemize}
We denote this model structure by $\vcat{\cV}$.
\end{thm}

We specialize this construction to the following cartesian model categories.

\begin{enumerate}[leftmargin=*]
\setcounter{enumi}{-1}
\item Let $\cV=\cat$ be the canonical model structure.
We then obtain precisely the model category $\vcat{\cat}=\twocat$ from \cref{lackMS}, as discussed in \cite[Ex.\ 1.8]{BergerMoerdijkEnriched}, in which every object is fibrant.

\item Let $\cV=\sset_{(\infty,1)}$ be the Joyal model structure.
We then obtain the model category $\vcat{\sset_{(\infty,1)}}$, in which the fibrant objects are the categories enriched over quasi-categories.

\item Let $\cV=\spsh{\Delta}_{(\infty,1)}$ being the Rezk model structure.
We then obtain the model category $\vcat{\spsh{\Delta}_{(\infty,1)}}$, in which the fibrant objects are the categories enriched over complete Segal spaces.
\item Let $\cV=\sset^+_{(\infty,1)}$ be the Lurie model structure on the category $\sset^+$.
We then obtain the model category $\vcat{\sset^+_{(\infty,1)}}$, in which the fibrant objects are the categories enriched over naturally marked quasi-categories.
\end{enumerate}

We recall from \cite[Thm~4.2.4]{CruttwellThesis} or \cite{EilenbergKelly}
that any lax monoidal functor $F\colon\cV\to\cV'$ induces a base-change functor $F_*\colon\vcat{\cV}\to\vcat{\cV'}$. This is in particular the case when $F$ is (strong) monoidal.
For any $\cV$-category $\cD$, the $\cV'$-category $F_*\cD$ has the same set of objects as $\cD$, and for any two objects $c,d$ in $\cD$ the hom-categories are defined by $(F_*\cD)(c,d)\coloneqq F(\cD(c,d))$. If $F\colon\cV\to\cV'$ is a right adjoint functor with a monoidal left adjoint functor $L\colon\cV'\to\cV$, then $L_*$ is the left adjoint of $F_*$.

\begin{prop}
\label{F*properties}
Let $\cV$, $\cV'$ be excellent monoidal model categories, and $F\colon\cV\to\cV'$ a right adjoint functor whose left adjoint functor is monoidal.
Denote by $F_*\colon\vcat{\cV}\to\vcat{\cV'}$ the induced base-change functor.
\begin{enumerate}[leftmargin=*]
    \item If $F$ is a right Quillen functor, then $F_*$ is a right Quillen functor.
    \item If $F$ is a right Quillen embedding, then $F_*$ is a right Quillen embedding.
    \item If $F$ is a Quillen equivalence, then $F_*$ is a Quillen equivalence.
\end{enumerate}
\end{prop}

\begin{proof}
Parts (1) and (3) are treated in \cite[Rmk.~A.3.2.6]{htt}, while Part (2) can easily be verified as a variant of (3).
\end{proof}

As special cases, we obtain the following model comparison functors.
\begin{enumerate}[label=(\alph*), leftmargin=*]
    \item The functor $(-)_{\bullet,0}\colon\spsh{\Delta}\to\sset$ is shown to be a right Quillen equivalence in \cite[\S4]{JT} and its left adjoint is product-preserving because it is a right adjoint itself, as discussed e.g.~in \cite[\S2]{JT}. We then obtain a right Quillen equivalence
\[
((-)_{\bullet,0})_*\colon\vcat{\spsh{\Delta}_{(\infty,1)}}\to \vcat{\sset_{(\infty,1)}}.
\]
    \item The underlying simplicial set functor $U\colon\sset^+\to\sset$ is a right Quillen equivalence by \cite[Thm~3.1.5.1]{htt} and 
 its left adjoint, given by the functor $(-)^{\flat}\colon\sset\to\sset^+$ which marks a simplicial set minimally, preserves finite products. We then obtain a right Quillen equivalence
\[
    U_*\colon\vcat{\sset^+_{(\infty,1)}}\to \vcat{\sset_{(\infty,1)}}.
    \]
\end{enumerate}

\subsection{The nerves}

The proposition can also be used to produce valuable nerve constructions.

\begin{const}
\label{EnrichedNerves} All the following base-change functors are special instances of \cref{F*properties}.
\begin{enumerate}[leftmargin=*]
    \item The ordinary nerve functor $\bf N\colon\cat\to\sset$ is a right Quillen embedding and its left adjoint functor preserves finite products by \cite[Prop.\ B.0.15]{JoyalVolumeII},
    there attributed to Gabriel--Zisman. We then obtain a right Quillen embedding
\[{\bf N}_*\colon\twocat\to \vcat{\sset_{(\infty,1)}}.\]
    \item The natural nerve functor \footnote{In the original source, ${\mathbf N}^{\natural}\cD$ is obtained as the value of a composite functor $\mathrm{N}^+\iota\cD$.} $\bf N^{\natural}\colon\cat\to\sset^+$ from \cite[Formula~(1.1)]{GHL} is a right Quillen embedding by \cite[Lem.~1.9]{GHL} and its left adjoint preserves finite products by \cite[\S1.1]{GHL}. We then obtain a right Quillen embedding
\[{\bf N}^{\natural}_*\colon\twocat\to \vcat{\sset^+_{(\infty,1)}}.\]
   \item The Rezk nerve functor\footnote{In the original source, ${\bf N}^R\cC$ is the \emph{classifying diagram} of $\cC$, denoted $N\cD$.} $\NRezk\colon\cat\to\spsh{\Delta}$ from \cite[\S3.5]{rezk}and recalled in \cref{AppendixRezk} is a right Quillen embedding by \cref{RezkRightQuillen} and we verify that its left adjoint preserves finite products in \cref{RezkInducesMonoidal}. We then obtain a right Quillen embedding
\[\NRezk_*\colon\twocat\to \vcat{\spsh{\Delta}_{(\infty,1)}}.\]
\end{enumerate}
\end{const}

The three nerve constructions are compatible with each other, as the next corollary shows.

\begin{cor}
\label{corcomparisonenriched}
The diagram of $\infty$-categories
\begin{tz}
\node[](1) {$[\twocat]_\infty$}; 
\node[below of=1](2) {$[\vcat{\sset_{(\infty,1)}}]_{\infty}$}; 
\node[left of=2,xshift=-2.3cm,yshift=-2pt](3) {$[\vcat{\spsh{\Delta}_{(\infty,1)}}]_{\infty}$}; 
\node[right of=2,xshift=2.3cm,yshift=-2pt](4) {$[\vcat{\sset^+_{(\infty,1)}}]_{\infty}$}; 

\draw[->] (1) to node[above,la,xshift=-15pt]{$[\NRezk_*]_\infty$} (3);
\draw[->] (1) to node[above,la,xshift=15pt]{$[\bf N^{\natural}_*]_\infty$} (4); 
\draw[->] ($(3.east)+(0,2.5pt)$) to node[below,la,xshift=2pt]{$[((-)_{\bullet,0})_*]_\infty$} ($(2.west)+(0,.5pt)$); 
\draw[->] ($(4.west)+(0,2.5pt)$) to node[below,la]{$[U_*]_\infty$} ($(2.east)+(0,.5pt)$);
\draw[->] (1) to node[right,la,xshift=-2pt]{$[\bfN_*]_\infty$} (2);
\end{tz}
commutes up to equivalence.
\end{cor}

\begin{proof}
The corollary is an application of the ``right Quillen'' version of \cref{NerveRecognition} to the following diagram.
\begin{tz}
\node[](1) {$\twocat$}; 
\node[below of=1](2) {$\vcat{\sset_{(\infty,1)}}$}; 
\node[left of=2,xshift=-1.5cm,yshift=-2pt](3) {$\vcat{\spsh{\Delta}_{(\infty,1)}}$}; 
\node[right of=2,xshift=1.5cm,yshift=-2pt](4) {$\vcat{\sset^+_{(\infty,1)}}$}; 

\draw[->] (1) to node[above,la,xshift=-10pt]{$\NRezk_*$} (3);
\draw[->] (1) to node[above,la,xshift=10pt]{$\bf N^{\natural}_*$} (4); 
\draw[->] ($(3.east)+(0,2.5pt)$) to node[below,la]{$((-)_{\bullet,0})_*$} ($(2.west)+(0,.5pt)$); 
\draw[->] ($(4.west)+(0,2.5pt)$) to node[below,la]{$U_*$} ($(2.east)+(0,.5pt)$);
\draw[->] (1) to node[right,la]{$\bfN_*$} (2);
\end{tz}
The fact that the diagram commutes up to isomorphism is a consequence of the fact that the diagram
\begin{tz}
\node[](1) {$\cat$}; 
\node[below of=1](2) {$\sset_{(\infty,1)}$}; 
\node[left of=2,xshift=-1.5cm](3) {$\spsh{\Delta}_{(\infty,1)}$}; 
\node[right of=2,xshift=1.5cm](4) {$\sset^+_{(\infty,1)}$}; 

\draw[->] (1) to node[above,la,xshift=-10pt]{$\NRezk$} (3);
\draw[->] (1) to node[above,la,xshift=10pt]{$\bf N^{\natural}$} (4); 
\draw[->] ($(3.east)+(0,.5pt)$) to node[below,la]{$((-)_{\bullet,0})$} ($(2.west)+(0,.5pt)$); 
\draw[->] ($(4.west)+(0,.5pt)$) to node[below,la]{$U$} ($(2.east)+(0,.5pt)$);
\draw[->] (1) to node[right,la]{$\bfN$} (2);
\end{tz}
commutes up to isomorphism.
\end{proof}

%%%
%%%
%%%

\subsection{Nerve comparison}

Bergner--Rezk consider an enriched nerve functor in \cite[Def.~7.3]{br1}, obtained by regarding a bisimplicial category as a simplicial object in simplicial spaces, and show that it defines
a right Quillen equivalence
\begin{equation}
    \label{PrecatQEenriched}
\BRequi\colon\vcat{\spsh{\Delta}_{(\infty,1)}}\to P\cat(\spsh{\Delta})_{(\infty,2)}.
\end{equation}
If $\cQ$ is a category enriched over simplicial spaces with object set $\cQ_0$, and $\cQ_1$ denotes the simplicial space \[\cQ_1=\coprod_{a,b \in \cQ_0}{\cQ}(a,b),\]
by definition of $\BRequi$ (as given in \cite[Def.~7.3]{br1})
there are isomorphisms of bisimplicial sets
\[(\BRequi\cQ)_{0}\cong\cQ_0\text{ and }(\BRequi\cQ)_{1}\cong\cQ_1,\]
and for any $i\ge0$
\begin{equation}
    \label{SegalObjectRmk2}
\left(\BRequi\cQ\right)_{i}\cong \underbrace{\cQ_1 \underset{\cQ_0}{\times}\cQ_1 \underset{\cQ_0}{\times} \ldots \underset{\cQ_0}{\times} \cQ_1}_{i}.
\end{equation}

First, we aim at giving an explicit description for $(\BRequi\NRezk_*\cD)_{i,j,k}$, which we achieve in \cref{RofNRezk}.

Given any category $\cA$ and $i\ge0$, we define inductively a $2$-category $\Sigma_i\cA$, called the \emph{$(i+1)$-point suspension} of~$\cA$. We set $\Sigma_0\cA\coloneqq [0]$, and for $i\ge1$ the $2$-category $\Sigma_i\cA$ can be understood as the pushout of $2$-categories
 \begin{tz}
 \node[](1) {$[0]$}; 
 \node[below of=1](2) {$\Sigma\cA$}; 
 \node[right of=1,xshift=.5cm](3) {$\Sigma_{i-1}\cA$}; 
 \node[below of=3,yshift=-1pt](4) {$\Sigma_{i}\cA$};
 \punctuation{4}{.}[3pt];
 
 \draw[->] (1) to (2);
 \draw[->] (1) to (3); 
 \draw[->] (2) to ($(4.west)+(0,1pt)$);
 \draw[->] (3) to (4);
 
\pushout{4};
 \end{tz}
The construction extends to a functor $\Sigma_i\colon\cat\to\twocat_{*,*}$.

\begin{lem}
\label{MultipleSuspensionPushout}
Given a category $\cA$ and $i\ge1$ there is a pushout of $2$-categories
 \begin{tz}
 \node[](1) {$\coprod\limits_{i+1} \cA$}; 
 \node[below of=1](2) {$\cA\otimes [i]$}; 
 \node[right of=1,xshift=1cm](3) {$\coprod\limits_{i+1}[0]$}; 
 \node[below of=3](4) {$\Sigma_i(\pi_0)_*\cA$};
 \punctuation{4}{.};
 
 \draw[->] ($(1.south)+(0,2pt)$) to (2);
 \draw[->] ($(1.east)+(0,4pt)$) to ($(3.west)+(0,4pt)$); 
 \draw[->] (2) to (4);
 \draw[->] ($(3.south)+(0,2pt)$) to (4);
 
\pushout{4};
 \end{tz} 
\end{lem}

\begin{proof}
The statement can be proven by induction on $i\ge1$. The basis of the induction, namely the case $i=1$, is precisely \cref{lemmasuspensionhigher}, and we now show the inductive step. 

For $i>1$, denote by $\cP$ the following pushout.
 \begin{tz}
 \node[](1) {$\coprod\limits_{i+1} \cA$}; 
 \node[below of=1](2) {$\cA\otimes [i]$}; 
 \node[right of=1,xshift=1cm](3) {$\coprod\limits_{i+1}[0]$}; 
 \node[below of=3](4) {$\cP$};
 
 \draw[->] ($(1.south)+(0,2pt)$) to (2);
 \draw[->] ($(1.east)+(0,4pt)$) to ($(3.west)+(0,4pt)$); 
 \draw[->] (2) to (4);
 \draw[->] ($(3.south)+(0,2pt)$) to (4);
 
\pushout{4};
 \end{tz}
Consider the following commutative diagram of $2$-categories.
\begin{tz}
\node[](1) {$\cA\otimes [1]$}; 
\node[right of=1,xshift=1cm](1') {$\cA\amalg\cA$}; 
\node[right of=1',xshift=1cm](1'') {$[0]\amalg [0]$};

\node[below of=1](2) {$\cA$};
\node[below of=1'](2') {$\cA$}; 
\node[below of=1''](2'') {$[0]$};

\node[below of=2,yshift=4pt](3) {$\cA\otimes [i-1]$};
\node[below of=2'](3') {$\coprod\limits_{i}\cA$};
\node[below of=2''](3'') {$\coprod\limits_{i}[0]$};

\draw[->] (2) to (1);
\draw[->] (1') to (1);
\draw[->] (2) to (3);
\draw[->] (1') to (1'');
\draw[d] (2') to (2);
\draw[->] (2') to (2'');
\draw[->] (2'') to (1'');
\draw[->] (2'') to (3''); 
\draw[->] (2') to (1'); 
\draw[->] (2') to (3');
\draw[->] ($(3'.west)+(0,4pt)$) to (3);
\draw[->] ($(3'.east)+(0,4pt)$) to ($(3''.west)+(0,4pt)$);
\end{tz}
The colimit of this diagram can be equivalently computed by either taking the colimit of the colimits of each row, or by taking the colimit of the colimits of each column. Following the first procedure, the resulting $2$-category is the pushout of the span
\begin{tz}
\node[](3) {$\cA\otimes [i]$};
\node[right of=3,xshift=1cm,yshift=-4pt](3') {$\coprod\limits_{i+1}\cA$};
\node[right of=3',xshift=1cm](3'') {$\coprod\limits_{i+1}[0]$};

\draw[->] ($(3'.west)+(0,4pt)$) to (3);
\draw[->] ($(3'.east)+(0,4pt)$) to ($(3''.west)+(0,4pt)$);
\end{tz}
    which gives precisely $\cP$.
    
    Instead, following the second procedure,
    the resulting $2$-category is by induction hypothesis the pushout of the span
    \begin{tz}
\node[](1) {$\Sigma(\pi_0)_*\cA$}; 
\node[below of=1](2) {$[0]$};
\node[below of=2](3) {$\Sigma_{i-1} (\pi_0)_*\cA$};

\draw[->] (2) to (1);
\draw[->] (2) to (3);
\end{tz}
  which is $\Sigma_{i}(\pi_0)_*\cA$. 
So the desired isomorphism follows.
\end{proof}

%%%%

For $i,j,k\ge0$, let $\Sigma_i([j]\times \widetilde{[k]})$ denote the $(i+1)$-point 
suspension of $[j]\times \widetilde{[k]}$, which is obtained by gluing $i$ consecutive copies of $\Sigma([j]\times \widetilde{[k]})$.

\begin{prop}
\label{RofNRezk}
For any $2$-category $\cD$ and $i,j,k\ge0$ we have
a natural bijection
\[(\BRequi\NRezk_*\cD)_{i,j,k}\cong\twocat(\Sigma_i([j]\times \widetilde{[k]}),\cD).\]
\end{prop}

\begin{proof}
For any $i\ge0$ we have a natural isomorphism of bisimplicial spaces
\begin{align*}
  \left(\BRequi\NRezk_*\cD\right)_{i}   &  \cong (\NRezk_*\cD)_1 \underset{(\NRezk_*\cD)_0}{\times}(\NRezk_*\cD)_1 \underset{(\NRezk_*\cD)_0}{\times} \ldots \underset{(\NRezk_*\cD)_0}{\times} (\NRezk_*\cD)_1\\
&  \cong (\NRezk_*\cD)_1 \underset{\Ob\cD}{\times}(\NRezk_*\cD)_1 \underset{\Ob\cD}{\times} \ldots \underset{\Ob\cD}{\times} (\NRezk_*\cD)_1\\
&  \cong \coprod\limits_{d_0,\ldots,d_i\in \Ob\cD}\NRezk\cD(d_0,d_1) \times\NRezk\cD(d_1,d_2)\times\ldots\times\NRezk\cD(d_{i-1},d_i)
\end{align*}
induced by the Segal maps. So for any $j,k\ge0$ we get a natural bijection 
\begin{align*}
  \left(\BRequi\NRezk_*\cD\right)_{i,j,k}&  \cong \coprod\limits_{d_0,\ldots,d_i\in \Ob\cD}\NRezk\cD(d_0,d_1)_{j,k} \times\ldots\times\NRezk\cD(d_{i-1},d_i)_{j,k}\\
&\cong \coprod\limits_{d_0,\ldots,d_i\in \Ob\cD}\cat([j]\times \widetilde{[k]},\cD(d_0,d_1))\times\ldots\times\cat([j]\times \widetilde{[k]},\cD(d_{i-1},d_i))\\
&\cong \twocat(\Sigma_i([j]\times \widetilde{[k]}),\cD),
\end{align*}
as desired.
\end{proof}

Next, we show the comparison between $\BRequi\NRezk_*\cD$ and $R\NMos\cD$.

\begin{thm}
\label{bisimplicialVSenriched}
For any $2$-category $\cD$ there is a natural map of bisimplicial spaces
\[\BRequi\NRezk_*\cD\to R\NMos\cD\]
 that is a weak equivalence in $\spsh{(\Delta\times\Delta)}_{(\infty,2)}$ and in $P\cat(\spsh{\Delta})_{(\infty,2)}$.
\end{thm}

First, we give a more general version of \cref{lemLyne2}.

\begin{lem}
\label{lemLyne3}
For any $i,j,k\ge0$ there is a pushout of $2$-categories
 \begin{tz}
 \node[](1) {$\coprod\limits_{i+1} \Osim{j} \tensoricon \Otilde{k}$}; 
 \node[below of=1](2) {$\Osim{j}\tensoricon ([i] \pseudo \Otilde{k})$}; 
 \node[right of=1,xshift=2.5cm,yshift=-2pt](3) {$\coprod\limits_{i+1}[0]$}; 
 \node[below of=3](4) {$\Sigma_i([j]\times \widetilde{[k]})$};
 \punctuation{4}{.};
 
 \draw[->] ($(1.south)+(0,8pt)$) to ($(2.north)-(0,3pt)$);
 \draw[->] ($(1.east)+(0,2pt)$) to ($(3.west)+(0,4pt)$); 
 \draw[->] ($(2.east)-(0,2pt)$) to (4);
 \draw[->] ($(3.south)+(0,2pt)$) to (4);
 
\pushout{4};
 \end{tz}
\end{lem}

\begin{proof}
The proof is similar to \cref{lemLyne2} replacing $[1]$ with $[i]$ and using \cref{MultipleSuspensionPushout}.
\end{proof}

We can now prove the theorem.

\begin{proof}[Proof of \cref{bisimplicialVSenriched}]
We first build the desired map. To this end, consider the following map of spans.
\begin{tz}
\node[](2) {$\Osim{j}\tensoricon (\Osim{i}\pseudo \Otilde{k})$};
\node[right of=2,xshift=3.5cm,yshift=-4pt](2') {$\coprod\limits_{i+1} \Osim{j} \tensoricon \Otilde{k}$}; 
\node[right of=2',xshift=2cm,yshift=-2pt](2'') {$\coprod\limits_{i+1} [0]$};

\node[below of=2](3) {$\Osim{j}\tensoricon ([i] \pseudo \Otilde{k})$};
\node[below of=2'](3') {$\coprod\limits_{i+1}\Osim{j}\tensoricon \Otilde{k}$};
\node[below of=2''](3'') {$\coprod\limits_{i+1} [0]$};

\draw[->] (2) to ($(3.north)-(0,3pt)$);
\draw[->] ($(2'.west)+(0,2pt)$) to ($(2.east)-(0,2pt)$);
\draw[->] ($(2'.east)+(0,2pt)$) to ($(2''.west)+(0,4pt)$);
\draw[d] ($(2''.south)+(0,2pt)$) to (3''); 
\draw[d] ($(2'.south)+(0,8pt)$) to ($(3'.north)-(0,3pt)$);
\draw[->] ($(3'.west)+(0,2pt)$) to ($(3.east)-(0,2pt)$);
\draw[->] ($(3'.east)+(0,2pt)$) to ($(3''.west)+(0,4pt)$);
\end{tz}
By \cref{LCexplicit}, the top row is given by the span 
\begin{tz}
\node[](2) {$L\bC\Delta[i,j,k]$};
\node[right of=2,xshift=2cm,yshift=-4pt](2') {$\coprod\limits_{i+1} L\bC\Delta [0,j,k]$}; 
\node[right of=2',xshift=2cm](2'') {$\coprod\limits_{i+1} L\bC\Delta[0,0,0]$};

\draw[->] ($(2'.west)+(0,4pt)$) to (2);
\draw[->] ($(2'.east)+(0,4pt)$) to ($(2''.west)+(0,4pt)$);
\end{tz}
and using \cref{Descriptionq*} and the fact that $L\bC$ commutes with colimits, its pushout is precisely $L\bC\DelQuot[i,j,k]$. By \cref{lemLyne3} the pushout of the bottom row is $\Sigma_i([j]\times \widetilde{[k]})$. Hence the map of spans yields the unique induced map of pushouts
\[
 L\bC\DelQuot[i,j,k] \to \Sigma_i([j]\times \widetilde{[k]}).
\]
Composing with the map in \cref{LandLsimeq}, we get a map
\[ L^{\simeq}\mathbb C\DelQuot[i,j,k]\to  L\mathbb C\DelQuot[i,j,k]\to  \Sigma_i([j]\times \widetilde{[k]}) \]
which induces by \cref{lem:q_*NHvsq_*NHsim,RofNRezk} a map of sets
\[(I\BRequi\NRezk_*\cD)_{i,j,k}=(\BRequi\NRezk_*\cD)_{i,j,k}\to\NMos_{i,j,k}\cD\]
which induces a map of bisimplicial spaces
\[I\BRequi\NRezk_*\cD\to\NMos\cD\]
which induces a map in $P\cat(\spsh{\Delta})$
\[\BRequi\NRezk_*\cD\to R\NMos\cD\]
as desired.

We now argue this map is a levelwise weak equivalence. By \cref{RofNLyne0,RofNLyne1,RofNRezk},
it induces isomorphisms in $\sset$ for $i=0,1$ and $j\ge0$
\[(\NRezk_*\cD)_{0,j}=(\BRequi\NRezk_*\cD)_{0,j}\xrightarrow{\cong} (R\NMos\cD)_{0,j},\;\; (\NRezk_*\cD)_{0,1}=(\BRequi\NRezk_*\cD)_{1,j}\xrightarrow{\cong} (R\NMos\cD)_{1,j}.\]
Using the fact that $R\NMos\cD$ and $\BRequi\NRezk_*\cD$ are Segal objects by \cref{SegalObjectRmk1} and \eqref{SegalObjectRmk2}, it follows that for any $i,j\ge0$ it induces a weak equivalence in $\sset_{(\infty,0)}$
\[(\BRequi\NRezk_*\cD)_{i,j}\to (R\NMos\cD)_{i,j},\]
showing that the desired map is a weak equivalence in $\spsh{(\Delta\times\Delta)}_{(\infty,2)}$.

Finally, the fact that the desired map is a weak equivalence in $P\cat(\spsh{\Delta})_{(\infty,2)}$ is a consequence of \cref{Ireflects}.
\end{proof}

Finally, we compare the nerves from \cref{EnrichedNerves}.

\begin{cor}
\label{corcomparisonPrecat}
The diagram of $\infty$-categories
\begin{tz}
\node[](1) {$[\twocat]_\infty$}; 
\node[below of=1](2) {$[P\cat(\spsh{\Delta})_{(\infty,2)}]_\infty$}; 
\node[left of=2,xshift=-2.3cm](3) {$[\spsh{(\Delta\times \Delta)}_{(\infty,2)}]_\infty$}; 
\node[right of=2,xshift=2.3cm,yshift=-2pt](4) {$[\vcat{\spsh{\Delta}_{(\infty,1)}}]_\infty$}; 

\draw[->] (1) to node[above,la,xshift=-15pt]{$[\NMos]_\infty$} (3);
\draw[->] (1) to node[above,la,xshift=15pt]{$[\NRezk_*]_\infty$} (4); 
\draw[->] (3) to node[below,la]{$[R]_\infty$} (2); 
\draw[->] ($(4.west)+(0,2pt)$) to node[below,la]{$[\BRequi]_\infty$} (2);

\draw[->] (1) to node[right,la,xshift=-2pt]{$[R\NMos]_\infty$} (2);
\end{tz}
commutes up to equivalence.
\end{cor}

\begin{proof}
The corollary follows from applying twice the ``right Quillen'' version of \cref{NerveRecognition} to the following diagram.
\begin{tz}
\node[](1) {$\twocat$}; 
\node[below of=1](2) {$P\cat(\spsh{\Delta})_{(\infty,2)}$}; 
\node[left of=2,xshift=-1.8cm](3) {$\spsh{(\Delta\times \Delta)}_{(\infty,2)}$}; 
\node[right of=2,xshift=1.8cm,yshift=-2pt](4) {$\vcat{\spsh{\Delta}_{(\infty,1)}}$}; 

\draw[->] (1) to node[above,la,xshift=-15pt]{$\NMos$} (3);
\draw[->] (1) to node[above,la,xshift=15pt]{$\NRezk_*$} (4); 
\draw[->] (3) to node[below,la]{$R$} (2); 
\draw[->] ($(4.west)+(0,2pt)$) to node[below,la]{$\BRequi$} (2);

\draw[->] (1) to node[right,la,xshift=-2pt]{$R\NMos$} (2);
\end{tz}
The fact that all the assumptions of the lemma are met are from \cref{MoserEmbedding,EnrichedNerves,bisimplicialVSenriched}.
\end{proof}

\section{Nerves in simplicial models}

\subsection{The models}

Verity envisioned a model of $(\infty,2)$-categories (part of a family of $(\infty,n)$-categories for general $n$) based on simplicial sets endowed with a subset of distinguished simplices.

\begin{defn}
A \emph{simplicial set with marking}\footnote{Originally referred to as \emph{stratified simplicial set} e.g.~in \cite[Def.\ 96]{VerityComplicialAMS}, \emph{simplicial sets with normality} \cite{StreetHandwritten} and \emph{hollow simplicial sets} \cite{StreetOrientedSimplexes}.}
is a simplicial set with a set of distinguished simplices -- called \emph{marked} -- in positive dimension and containing degenerate simplices.
\end{defn}

Amongst all simplicial sets with marking, the following identify those that are $(\infty,2)$-categories.
The following mathematical object was identified by Verity \cite{VeritySlides} as a model for $(\infty,2)$-categories, and was further studied in \cite[\S3.3]{EmilyNotes}, \cite[\S1.3]{or} and \cite[App.~D]{RiehlVerityBook}.

\begin{defn}
\label{defcomplicial}
A \emph{saturated $2$-complicial set}\footnote{Sometimes for brevity referred to as \emph{$2$-complicial set}.}
is a simplicial set that has the right lifting property with respect to all maps of the following kinds:
\begin{enumerate}[leftmargin=*]
 \item for $m> 1$ and $0< k< m$, the \emph{complicial inner horn extension}
\[\Lambda^k[m]\to \Delta^k[m];\]
here, $\Delta^k[m]$ is the standard $m$-simplex in which a non-degenerate simplex is marked if and only if it contains the vertices $\{k-1,k,k+1\}\cap [m]$, and $\Lambda^k[m]$ is the regular sub-simplicial set with marking of $\Delta^k[m]$ whose simplicial set is the $k$-horn $\Lambda^k[m]$;
 \item for $m\geq 2$ and $0< k < m$, the \emph{complicial thinness extension}
\[\Delta^k[m]' \to \Delta^k[m]'';\]
here, $\Delta^k[m]'$ is the standard $m$-simplex with marking obtained from $\Delta^k[m]$ by additionally marking the $(k-1)$-st and $(k+1)$-st face of $\Delta[m]$, and $\Delta^k[m]''$ is the standard $m$-simplex with marking obtained from $\Delta^k[m]'$ by additionally marking the $k$-th face of $\Delta[m]$;
\item for $m>2$, the \emph{triviality extension}
\[\Delta[m]\to\Delta[m]_t;\]
here, $\Delta[m]$ is the minimally marked $m$-simplex, and $\Delta[m]_t$ is the thin $m$-simplex in which the only non-degenerate simplex marked is the unique $m$-simplex;
\item for $m\ge-1$, the \emph{complicial saturation extension}
\[\Delta[3]_{\mathrm{eq}}\star\Delta[m]\to\Delta[3]_{\sharp}\star\Delta[m];\]
here, $\Delta[3]_{\mathrm{eq}}$ is the standard $3$-simplex with marking given by all simplices in dimension at least $2$, as well as the $1$-simplices $[0,2]$ and $[1,3]$, and $\Delta[3]_{\sharp}$ is the standard $3$-simplex with the maximal marking.
\end{enumerate}
See e.g.~\cite[Def.~1.19]{or} for more details. We refer the reader to \cite{VerityComplicialAMS} for the join $\star\colon\msset\times\msset\to\msset$ of marked simplicial sets.
\end{defn}

The following model structure is obtained as an application of Verity's machinery from \cite[Thm~100]{VerityComplicialI}, and was further studied in \cite[\S4.3]{EmilyNotes}, and \cite[Thm~1.25]{or}.

\begin{thm}
\label{VerityModelStructure}
The category $\msset$ of simplicial sets with marking admits a model structure, denoted $\msset_{(\infty,2)}$, in which
\begin{itemize}[leftmargin=*]
    \item the fibrant objects are the saturated $2$-complicial sets, and
    \item the cofibrations are the monomorphisms \textnormal{(}of underlying simplicial sets\textnormal{)}, and in particular every object is cofibrant.
\end{itemize}
\end{thm}

Lurie proposed a simplified variant of this idea that focuses on the study of $(\infty,2)$-categories (as opposed to $(\infty,n)$-categories for general $n$), based on simplicial sets with marking only in dimension $2$.

\begin{defn}[{\cite[Def.\ 3.1.1]{LurieGoodwillie}}]
A \emph{scaled simplicial set} is a simplicial set with a \emph{scaling}, namely a set of distinguished $2$-simplices -- called \emph{marked} or \emph{thin} -- containing degenerate $2$-simplices.
\end{defn}

Amongst all scaled simplicial sets, the following identify those that are $(\infty,2)$-categories.
To recall this definition, we use the author's original convention that we denote a simplicial set with marking by listing a pair $(X,T)$ where $X$ is the underlying simplicial set, and $T$ is the set of non-degenerate scaled simplices.

\begin{defn}[{\cite[Def.~4.1.1]{LurieGoodwillie}}]
An \emph{$\infty$-bicategory}\footnote{This was originally referred to as a \emph{weak $\infty$-bicategory}, but was shown by Gagna--Harpaz--Lanari in \cite[Thm~5.1]{GHL} to agree with the original definition of \emph{$\infty$-bicategory} from \cite[Def.~4.2.8]{LurieGoodwillie}} is a simplicial set that has the right lifting property with respect to all maps indicated in \cite[Def.~3.1.3]{LurieGoodwillie}, namely
\begin{enumerate}[leftmargin=*]
 \item for $m\ge2$ and $0<k<m$ the \emph{scaled inner horn extension}
 \[(\Lambda^k[m],\{[k-1,k,k+1]\})\to (\Delta[m],\{[k-1,k,k+1]\});\]
  \item for $n\ge3$ the \emph{scaled outer horn extension}
 \[(\Lambda^0[m]\aamalg{\Delta[1]}\Delta[0],\{[0,1,n]\})\to (\Delta[m]\aamalg{\Delta[1]}\Delta[0],\{[0,1,n]\}), \]
where the pushouts are induced by the map $\langle 0,1\rangle \colon\Delta[1]\to\Delta[m]$;
 \item the \emph{scaled saturation extension}
\[(\Delta[4],T) \to (\Delta[4],\{T\cup\{[0,3,4],[0,1,4]\}), \]
where $T=\{[0,2,4],[1,2,3],[0,1,3],[1,3,4],[0,1,2]\}$.
\end{enumerate}
\end{defn}

The following model structure is obtained as an application of Smith Theorem.

\begin{thm}[{\cite[Thm~4.2.7]{LurieGoodwillie}}]
\label{LurieModelStructure}
The category $\sset^{sc}$ of scaled simplicial sets admits a model structure, denoted $\sset^{sc}_{(\infty,2)}$, in which
\begin{itemize}[leftmargin=*]
    \item the fibrant objects are the $\infty$-bicategories, and
    \item the cofibrations are the monomorphisms \textnormal{(}of underlying simplicial sets\textnormal{)}, and in particular every object is cofibrant.
\end{itemize}
\end{thm}

Gagna--Harpaz--Lanari prove in \cite[Thm 7.9]{GHL} that the canonical forgetful functor defines a right Quillen equivalence
\[
    U\colon\msset_{(\infty,2)}\to{\sset}^{sc}_{(\infty,2)},
\]

A further variant of Verity's original framework is given by working with \emph{$t\Delta$-sets}, where $t\Delta$ is an enlargement of the ordinary simplex category $\Delta$. More precisely, the category $t\Delta$ contains~$\Delta$ as a non-full subcategory, and in addition to the objects $[n]$ for $n\ge0$ it also contains objects of the form $[n]_t$ together with a map $[n]\to[n]_t$ for each $n\geq 1$.  We refer the reader to \cite[Not.\ 1.1]{or} or \cite[Not.\ D.1.4]{RiehlVerityBook} for more details on the category $t\Delta$.

Any $t\Delta$-set $X\colon t\Delta^{\op}\to\set$ can be seen as a \emph{simplicial set with multiple marking}. The underlying simplicial set of $X$ is the restriction of $X$ along the inclusion $\Delta^{\op}\to t\Delta^{\op}$,
so $X([n])=X_n$ is the set of $n$-simplices, while $X([n]_t)$ is the set of marked $n$-simplices; by definition, there is a structure map $X([n]_t)\to X([n])=X_n$ for every $n\ge1$, that remembers which simplex each marking belongs to. Notice that an $n$-simplex can be marked multiple times, namely, multiple elements of $X([n]_t)$ can map to the same element in $X_n$. According to this interpretation, simplicial sets with marking are precisely the $t\Delta$-sets for which all structure maps $X([n]_t)\to X([n])=X_n$ are monomorphisms\footnote{This approach looks more complicated at first glance, but offers certain technical advantages because, unlike the category of simplicial sets with marking, the category of $t\Delta$-sets is a category of presheaves.}, and there is an inclusion $\msset\hookrightarrow\psh{t\Delta}$.

\begin{defn}[{\cite[Def.\ 1.23]{or}}]
A \emph{$2$-precomplicial set}\footnote{We warn the reader that the same terminology is also used in \cite[\S6]{VerityComplicialAMS} to mean something unrelated.}
is a $t\Delta$-set that has the right lifting property with respect to the maps of the kinds (1)-(4) from \cref{defcomplicial}.
\end{defn}

The following model structure is an application of Cisinski's machinery from \cite[\S1.3]{cisinski}.

\begin{thm}[{\cite[Thm 1.28]{or}}]
\label{ORModelStructure}
The category $\psh{t\Delta}$ of simplicial sets with multiple marking admits a model structure, denoted $\psh{t\Delta}_{(\infty,2)}$, in which
\begin{itemize}[leftmargin=*]
    \item the fibrant objects are the saturated $2$-precomplicial sets, and
    \item the cofibrations are the monomorphisms \textnormal{(}of underlying simplicial sets\textnormal{)}, and in particular every object is cofibrant.
\end{itemize}
\end{thm}

The inclusion $\msset\hookrightarrow\psh{t\Delta}$ admits a left adjoint $\mathrm{Refl}$, which was proven by the second and third author as \cite[Prop.~1.31]{or} to be a left Quillen equivalence
\begin{equation}
\label{ORQE}
    \mathrm{Refl}\colon\psh{t\Delta}_{(\infty,2)}\to\msset.
\end{equation}
Given a $t\Delta$-set $X$, the functor $\mathrm{Refl}$ preserves the underlying simplicial set, so that we have $(\mathrm{Refl}X)_n=X_n=X([n])$, and the set of marked $n$-simplices $(\mathrm{Refl}X)([n]_t)$ is determined by the epi-mono factorization of the structure map
\[X([n]_t)\twoheadrightarrow(\mathrm{Refl}X)([n]_t)\hookrightarrow X([n]).\]
This means that an $n$-simplex is marked in $\mathrm{Refl}X$ if and only if it has at least one marking in $X$.

\subsection{The nerves}

Nerve constructions have been identified for the three discussed simplicial models of $(\infty,2)$-categories, and they are all based on the same underlying simplicial set: the Duskin nerve\footnote{In the original source, $\Nduskin\cD$ is denoted $\mathbf{Ner}\cD$.} $\Nduskin\cD$ of a $2$-category $\cD$ from \cite[\S6]{duskin}.

The \emph{Duskin nerve} $\Nduskin\cD$ of a $2$-category $\cD$ is the ($3$-coskeletal) simplicial set in which the set of $n$-simplices is given by
\[(\Nduskin\cD)_n\coloneqq\twocat(\cO_2[n],\cD).\]
 The assignment extends to a functor $\Nduskin\colon \twocat\to\sset$. In particular,
\begin{enumerate}[leftmargin=*]
\setcounter{enumi}{-1}
\item a $0$-simplex consists of an object $x$ of $\cD$;
    \item a $1$-simplex consists of
    a $1$-morphism $a\colon x\to y$ of $\cD$;
    \item a $2$-simplex consists of
    a $2$-cell  $\varphi\colon c\Rightarrow b\circ a$ of $\cD$ of the form
\begin{tz}
\node[](1) {$x$};
\node[above right of=1,xshift=.3cm](2) {$y$};
\node[below right of=2,xshift=.3cm](3) {$z$};
\punctuation{3}{;}[1pt];
\draw[->] (1) to node[above,la,xshift=-2pt]{$a$} (2);
\draw[->] (2) to node[above,la,xshift=2pt]{$b$} (3);
\draw[->] (1) to node[below,la]{$c$} (3);
\coordinate(a) at ($(1)!0.5!(3)$);
\cell[la,right][n][0.37]{a}{2}{$\varphi$};
\end{tz}
%     $$\begin{tikzcd}[baseline=(current  bounding  box.center)]
%  & y \arrow[rd, "b"]  & \\%{\alpha_s}
%     x \arrow[ru, "{a}"]
%   \arrow[rr, "c"{below}, ""{name=D,inner sep=1pt}]
%   && z;
%   \arrow[Rightarrow, from=D, 
%  to=1-2, shorten >= 0.1cm, shorten <= 0.1cm, "\varphi"]
% \end{tikzcd}$$
\item a $3$-simplex consists of four $2$-cells of $\cD$ that satisfy the following pasting equality.
\begin{tz}
\node[](1) {$x$};
\node[above of=1](2) {$y$};
\node[right of=2](3) {$z$};
\node[below of=3](4) {$w$};
\draw[->] (1) to node[left,la]{$a$} (2);
\draw[->] (2) to node[above,la]{$b$} (3);
\draw[->] (3) to node[right,la]{$d$} (4);
\draw[->] (1) to node[below,la]{$f$} (4);
\draw[->] (1) to node[over,la]{$c$} (3);
\coordinate(a) at ($(1)!0.5!(3)$);
\cell[la,right,yshift=7pt][n][0.45]{a}{2}{};
\coordinate(a) at ($(1)!0.7!(4)$);
\coordinate(b) at ($(a)+(0,1cm)$);
\cell[la,left][n][0.37]{a}{b}{};

\node[right of=4](1) {$x$};
\node at ($(3)!0.5!(1)$) {$=$};
\node[above of=1](2) {$y$};
\node[right of=2](3) {$z$};
\node[below of=3](4) {$w$};
\draw[->] (1) to node[left,la]{$a$} (2);
\draw[->] (2) to node[above,la]{$b$} (3);
\draw[->] (3) to node[right,la]{$d$} (4);
\draw[->] (1) to node[below,la]{$f$} (4);
\draw[->] (2) to node[over,la]{$e$} (4);
\coordinate(a) at ($(2)!0.5!(4)$);
\cell[la,right,yshift=-7pt][n][0.45]{a}{3}{};
\coordinate(a) at ($(1)!0.3!(4)$);
\coordinate(b) at ($(a)+(0,1cm)$);
\cell[la,right][n][0.37]{a}{b}{};
\end{tz}
\end{enumerate}
The face maps can be read off from the pictures.

\begin{const}[{\cite[Const.~4.8]{Nerves2Cat}}]
\label{ORNerve}
Let $\cD$ be a $2$-category. The nerve $\NOR\cD$ is the simplicial set $\Nduskin\cD$ with marking given by the following:
\begin{enumerate}[leftmargin=*]
\item all $1$-simplices inhabited by equivalences, each marked as many times as ways of completing the equivalence to an adjoint equivalence;
\item all $2$-simplices inhabited by isomorphisms, each marked uniquely;
    \item all simplices in dimension higher than $2$, each marked uniquely.
\end{enumerate}
This assignment extends to a functor $\NOR\colon \twocat\to\psh{t\Delta}$.
\end{const}

\begin{rmk}
Given $\cD$ a $2$-category, $\mathrm{Refl}\,\NOR\cD$ is the simplicial set $\Nduskin\cD$ endowed with the marking 
described in \cite[Prop.~3.1.10]{EmilyNotes}.
Essentially, the difference between $\NOR\cD$ and $\mathrm{Refl}\,\NOR\cD$ is that in the former each $1$-equivalence is marked many times, while in the latter it is marked only once (without remembering the data of any specific adjoint equivalence).\footnote{Another marking on $\NOR\cD$ considered in the literature is the Roberts--Street nerve from e.g.~\cite{VerityComplicialAMS}, for which the marked simplices are those inhabited by an identity cell. This nerve has important properties, but is not homotopically well-behaved, and does not play a role in this paper.}
\end{rmk}

\begin{const}[{\cite[Def.\ 8.1]{GHL}}]
\label{GHLNerve}
Let $\cD$ be a $2$-category. The nerve\footnote{In the original source, $\NGHL\cD$ is denoted $\mathcal N_2\cD$.} $\NGHL\cD$ is the simplicial set $\Nduskin\cD$ with scaling given by the set of all $2$-simplices inhabited by isomorphisms.
 The assignment extends to a functor $\NGHL\colon \twocat\to\sset^{sc}$.
\end{const}

These nerve constructions are well behaved homotopically.

\begin{thm}[{\cite[Thms\ 4.10,4.12]{Nerves2Cat}}]
\label{OREmbedding}
The functor $\NOR\colon \twocat\to\spsh{t\Delta}_{(\infty,2)}$ is a right Quillen embedding, and in particular a homotopical and right Quillen functor.
\end{thm}

\begin{thm}[{\cite[Prop.~8.2, 8.3]{GHL}}] \label{GHLEmbedding}
The functor $\NGHL\colon \twocat\to\sset^{sc}_{(\infty,2)}$ is a right Quillen embedding, and in particular a homotopical and right Quillen functor.
\end{thm}

\begin{rmk}
The functor $\mathrm{Refl}\,\NOR\colon \twocat\to\msset$ is \emph{not} a right adjoint functor. Indeed, if it admitted a left adjoint $L\colon \msset\to\twocat$, then we would have a natural bijection for any $2$-category $\cD$
\[\twocat(L\Delta[1]_t,\cD)\cong\msset(\Delta[1]_t,\mathrm{Refl}\,\NOR\cD)\cong (\mathrm{Refl}\,\NOR\cD)([1]_t)\cong\mathrm{eq}\cD,\]
where $\mathrm{eq}\cD$ denotes the set of equivalences in $\cD$.
However, one can use e.g.~\cite[Prop.\ 2.4.8]{RiehlCTC} to see that the functor $\mathrm{eq}\colon\twocat\to\set$ given by $\cD\mapsto\mathrm{eq}\cD$ is not corepresentable, obtaining a contradiction.
\end{rmk}

\begin{prop}
\label{RVEmbedding}
The functor $\mathrm{Refl}\,\NOR\colon \twocat\to\msset_{(\infty,2)}$ is homotopical and induces a fully faithful functor at the level of $\infty$-categories.
\end{prop}

\begin{proof}
The functor $\mathrm{Refl}\,\NOR\colon \twocat\to\msset_{(\infty,2)}$ is the composite of the right Quillen functor $\NOR\colon \twocat\to\psh{t\Delta}_{(\infty,2)}$ from \cref{OREmbedding}, followed by the left Quillen functor $\mathrm{Refl}\colon\psh{t\Delta}_{(\infty,2)}\to\msset_{(\infty,2)}$ from \eqref{ORQE}, which are both in particular homotopical and homotopically fully faithful.
Hence, $\mathrm{Refl}\,\NOR$ is homotopical and homotopically fully faithful.
\end{proof}

\subsection{Nerve comparisons}

The nerve constructions are compatible with each other as follows.

\begin{prop}
\label{WithinSimplicial}
For any $2$-category $\cD$ there is an isomorphism of scaled simplicial sets
\[\NGHL\cD\cong U\mathrm{Refl}\,\NOR\cD.\]
\end{prop}

\begin{proof}
The two scaled simplicial sets $\NGHL\cD$ and $U\mathrm{Refl}\,\NOR\cD$ have the same underlying simplicial set, given by the Duskin nerve $\Nduskin\cD$, and by reading through the relevant definitions and the explicit description of the reflector one can see that the marked $2$-simplices are precisely those inhabited by a $2$-isomorphism of $\cD$.
\end{proof}

\begin{cor}
\label{corsimplicial}
The diagram of $\infty$-categories 
\begin{tz}
\node[](1) {$[\twocat]_\infty$}; 
\node[below of=1](2) {$[\msset_{(\infty,2)}]_{\infty}$}; 
\node[left of=2,xshift=-2cm](3) {$[\psh{t\Delta}_{(\infty,2)}]_{\infty}$}; 
\node[right of=2,xshift=2cm](4) {$[\sset^{sc}_{(\infty,2)}]_{\infty}$}; 

\draw[->] (1) to node[above,la,xshift=-15pt]{$[\NOR]_\infty$} (3);
\draw[->] (1) to node[above,la,xshift=15pt]{$[\NGHL]_\infty$} (4); 
\draw[->] ($(3.east)+(0,1pt)$) to node[below,la]{$[\mathrm{Refl}]_\infty$} ($(2.west)+(0,1pt)$); 
\draw[->] ($(2.east)+(0,1pt)$) to node[below,la]{$[U]_\infty$} ($(4.west)+(0,1pt)$);

\draw[->] (1) to node[right,la,xshift=-2pt]{$[\mathrm{Refl}\,\NOR]_\infty$} (2);
\end{tz}
commutes up to equivalence.
\end{cor}

\begin{proof}
The commutativity of the left triangle is an application of the ``left Quillen'' version of \cref{NerveRecognition} to the diagram 
\begin{tz}
\node[](1) {$\twocat$}; 
\node[below of=1](2) {$\msset_{(\infty,2)}$}; 
\node[left of=2,xshift=-1.3cm](3) {$\psh{t\Delta}_{(\infty,2)}$}; 

\draw[->] (1) to node[above,la,xshift=-15pt]{$\NOR$} (3);
\draw[->] ($(3.east)+(0,1pt)$) to node[below,la]{$\mathrm{Refl}$} ($(2.west)+(0,1pt)$); 

\draw[->] (1) to node[right,la]{$\mathrm{Refl}\,\NOR$} (2);
\end{tz}
where the assumptions of the lemma are met by \cref{OREmbedding,RVEmbedding}. Then, 
the commutativity of the right triangle is an application of the ``right Quillen'' version of \cref{NerveRecognition} to the diagram
\begin{tz}
\node[](1) {$\twocat$}; 
\node[below of=1](2) {$\msset_{(\infty,2)}$}; 
\node[right of=2,xshift=1.5cm](4) {$\sset^{sc}_{(\infty,2)}$}; 

\draw[->] (1) to node[above,la,xshift=15pt]{$\NGHL$} (4); 
\draw[->] ($(2.east)+(0,1pt)$) to node[below,la]{$U$} ($(4.west)+(0,1pt)$);

\draw[->] (1) to node[left,la]{$\mathrm{Refl}\,\NOR$} (2);
\end{tz}
where the assumptions of the lemma are met by \cref{GHLEmbedding,RVEmbedding,WithinSimplicial}.
\end{proof}

We now discuss how the nerve constructions of simplicial models compare with those from the enriched models. Lurie showed as \cite[Thm 0.0.3]{LurieGoodwillie} that
the \emph{scaled homotopy coherent nerve}\footnote{In the original source, $\mathfrak N^{sc}\cD$ is denoted $N^{sc}\cD$.} functor introduced as \cite[Def.\ 3.1.10]{LurieGoodwillie} defines a right Quillen equivalence
\[\schoco\colon\vcat{{\sset^+}_{(\infty,1)}}\to\sset^{sc}_{(\infty,2)}.\]

\begin{prop}[{\cite[Prop.~8.2]{GHL}}]
\label{enrichedVSsimplicial}
For any $2$-category $\cD$ there is an isomorphism of scaled simplicial sets
\[\schoco{\bf N}^\natural_*\cD\cong \NGHL\cD. \]
\end{prop}

\begin{cor}
\label{corComplicialEnriched}
The diagram of $\infty$-categories 
\begin{tz}
\node[](1) {$[\twocat]_\infty$}; 
\node[below of=1,xshift=-2.5cm,yshift=-1pt,yshift=.2cm](3) {$[\vcat{\sset^+_{(\infty,1)}}]_{\infty}$}; 
\node[below of=1,xshift=2.5cm,yshift=.2cm](4) {$[\sset^{sc}_{(\infty,2)}]_{\infty}$}; 

\draw[->] (1) to node[above,la,xshift=-15pt]{$[\bf N^\natural_*]_\infty$} (3);
\draw[->] (1) to node[above,la,xshift=15pt]{$[\NGHL]_\infty$} (4); 
\draw[->] ($(3.east)+(0,2pt)$) to node[below,la]{$[\schoco]_\infty$} ($(4.west)+(0,1pt)$); 
\end{tz}
commutes up to equivalence.
\end{cor}

\begin{proof}
The corollary is an application of the ``right Quillen'' version of \cref{NerveRecognition} to the following diagram.
\begin{tz}
\node[](1) {$\twocat$}; 
\node[below of=1,xshift=-2.2cm,yshift=-1pt,yshift=.2cm](3) {$\vcat{\sset^+_{(\infty,1)}}$}; 
\node[below of=1,xshift=2.2cm,yshift=.2cm](4) {$\sset^{sc}_{(\infty,2)}$}; 

\draw[->] (1) to node[above,la,xshift=-15pt]{$\bf N^\natural_*$} (3);
\draw[->] (1) to node[above,la,xshift=15pt]{$\NGHL$} (4); 
\draw[->] ($(3.east)+(0,2pt)$) to node[below,la]{$\schoco$} ($(4.west)+(0,1pt)$); 
\end{tz}
The fact that all the assumptions of the lemma are met are from \cref{EnrichedNerves,enrichedVSsimplicial,GHLEmbedding}.
\end{proof}

%%%

\section{Nerves of \texorpdfstring{$2$}{2}-categories as local \texorpdfstring{$(\infty,2)$}{(infinity,2)}-categories}
\label{GHcharacterization}
\setcounter{shorten}{1}

The goal of this subsection is to prove \cref{MainTheorem2}, which will be completed in \cref{2dim}.
The ingredients for the proof are
\cref{suspensioncomparison,QE2,equivalencewithGH2}. We also use some of their $0$-dimensional analogs -- \cref{QE0,equivalencewithGH0} -- and $1$-dimensional analogs -- \cref{QE1,equivalencewithGH1} -- which are treated in \cref{0dim,1dim}, respectively.

\subsection{The \texorpdfstring{$0$}{0}-dimensional case}

\label{0dim}

The goal of this subsection is to show that the Quillen pair
\[ \pi_0\colon\sset_{(\infty,0)}\rightleftarrows\set\colon\mathrm{disc}\]
is equivalent to the left Bousfield 
localization of the Kan--Quillen model structure $\sset_{(\infty,0)}$ with respect to a set $\Lambda$ of maps. We also discuss in \cref{equivalencewithGH0} that this entails that the discrete embedding realizes sets as local $(\infty,0)$-categories with respect to the set of maps~$\Lambda$.

Recall from e.g.~\cite{CamarenaSets} that there is a canonical model structure on $\set$ in which the weak equivalences are the bijections, and every object is fibrant and cofibrant. Recall from \cite{QuillenHA} that the category $\sset$ of simplicial sets admits the Kan--Quillen model structure $\sset_{(\infty,0)}$, in which the weak equivalences are the weak homotopy equivalences, everything is cofibrant and the fibrant objects are precisely the Kan complexes.

The functor $\mathrm{disc}\colon\set\to\sset$ that regards each set as a discrete simplicial set admits a left adjoint given by the functor $\pi_0\colon\sset\to\set$ that takes a simplicial set to its set of connected components. The following is a straightforward verification.

\begin{prop}
The functor $\mathrm{disc}
\colon\set\to\sset_{(\infty,0)}$ is a right Quillen embedding.
\end{prop}

In particular, we have a Quillen reflection pair
\[
    \pi_0\colon\sset_{(\infty,0)}\rightleftarrows\set\colon\mathrm{disc}.
\]

\begin{rmk}
The essential image of the functor $[\mathrm{disc}]_\infty\colon[\set]_\infty\to[\sset_{(\infty,0)}]_\infty$ is the full sub-$\infty$-category of $[\sset_{(\infty,0)}]_\infty$ generated by the homotopically discrete $(\infty,0)$-categories.
\end{rmk}

For $k>0$, let $S^k\coloneqq\partial\Delta[k]$ denote the simplicial $k$-sphere. Since the model structure $\sset_{(\infty,0)}$ is combinatorial and left proper, the following model structure exists.

\begin{prop}
The category $\sset$ admits the left Bousfield localization $\cL_{\Lambda}\sset_{(\infty,0)}$ of the model structure $\sset_{(\infty,0)}$ with respect to the set $\Lambda$ of maps of the form
\[\Delta[0]\hookrightarrow S^k, \quad\text{for $k>0$}.\]
\end{prop}

In particular, there is a Quillen reflection pair
\[
\mathrm{Id}\colon\sset_{(\infty,0)}\rightleftarrows \cL_{\Lambda}\sset_{(\infty,0)}\colon \mathrm{Id}.
\]

The following is a straightforward verification.

\begin{prop}
\label{QE0}
The functor $\mathrm{disc}\colon\set\to \cL_{\Lambda}\sset_{(\infty,0)}$ defines a right Quillen equivalence.
\end{prop}

The following relates two approaches to localizations of $
\infty$-categories and is classical, but it is described e.g.~in the proof of \cite[Prop.\ A.3.7.8]{htt}.

We refer the reader to \cite[Def.\ 5.2.7.2,~Prop.\ 5.5.4.15]{htt} for a discussion on the localization $\mathscr L_S\mathscr Q$ of a quasi-category $\mathscr Q$ with respect to a set of edges $S$, and to \cite[Ch.\ 3]{hirschhorn}
for the left Bousfield localization $\cL_S\cM$ of a model category $\cM$ with respect to a set of morphisms $S$, namely the localization in the context of model categories.

\begin{prop}
\label{localization compatible}
Given a combinatorial left proper model category $\cM$ and a set of maps~$S$, denote by $\cL_{S}\cM$ the left Bousfield localization and by $\mathscr L_{S}[\cM]_\infty$ the localization in the sense of $\infty$-categories. Then there is a diagram of $\infty$-categories
\begin{tz}
\node[](1) {$[\cL_S\cM]_\infty$}; 
\node[right of=1,xshift=1.8cm](2) {$\mathscr L_S[\cM]_\infty$};
\node[below of=1,xshift=1.8cm,yshift=.3cm](3) {$[\cM]_\infty$};
\draw[->] (1) to node[above,la]{$\simeq$} (2);
\draw[right hook->] (1) to node[below,la,xshift=-10pt]{$[\mathrm{Id}_\cM]_\infty$}(3);
\draw[left hook->] (2) to (3);
\end{tz}
that commutes up to equivalence.
\end{prop}

With the following remark we verify that the map of $\infty$-categories induced by \cref{QE0} does implement the inclusion of the $\infty$-category of sets into the $\infty$-category of spaces considered by Gepner--Haugseng in \cite[\S6]{GH}.

\begin{rmk}
\label{equivalencewithGH0}
We know -- and it is also mentioned in \cite[\S6]{GH} -- that the underlying $\infty$-category of the Kan--Quillen model structure $\sset_{(\infty,0)}$ models the established $\infty$-category~$\mathscr S$ of spaces, meaning there exists an equivalence of $\infty$-categories
\begin{equation}
    \label{correctspaces}
  [\sset_{(\infty,0)}]_\infty\simeq \mathscr S.
\end{equation}
Any such equivalence can be used to construct a specific
equivalence of $\infty$-categories
\begin{align}\label{correctspaces2}
    [\set]_\infty&\simeq [\cL_\Lambda \sset_{(\infty,0)}]_\infty&\text{\cref{QE0}} \nonumber\\
    &\simeq\mathscr L_\Lambda[\sset_{(\infty,0)}]_\infty&\text{\cref{localization compatible}}\\
    &\simeq\mathscr L_\Lambda\mathscr S&\text{\eqref{correctspaces}}\nonumber\\
     &\simeq\mathscr Set&\text{\cite[Lem.~6.1.6(1)]{GH}}\nonumber
\end{align}
     between the ($\infty$-)category of sets $\mathscr Set$ and the underlying ($\infty$-)category of the model structure $\set$ on sets. Via the chosen identifications \eqref{correctspaces} and \eqref{correctspaces2},
we see that the functor $[\mathrm{disc}]_\infty\colon[\set]_\infty\to [\sset_{(\infty,0)}]_\infty$ and the canonical inclusion $\mathscr Set\hookrightarrow\mathscr S$ from \cite[Def.~6.1.6(i)]{GH} are equivalent. Indeed, this is witnessed by the following diagram of $\infty$-categories
\begin{tz}
\node[](1) {$[\set]_\infty$}; 
\node[below of=1](1') {$[\sset_{(\infty,0)}]_\infty$}; 
\node[right of=1,xshift=1.5cm](2) {$[\cL_\Lambda\sset_{(\infty,0)}]_\infty$}; 
\node[below of=2](2') {$[\sset_{(\infty,0)}]_\infty$}; 
\node[right of=2,xshift=2cm](3) {$\mathscr L_\Lambda[\sset_{(\infty,0)}]_\infty$}; 
\node[below of=3](3') {$[\sset_{(\infty,0)}]_\infty$}; 
\node[right of=3,xshift=1.5cm,yshift=1pt](4) {$\mathscr L_\Lambda \mathscr S$}; 
\node[below of=4](4') {$\mathscr S$}; 
\node[right of=4,xshift=.8cm](5) {$\mathscr Set$}; 
\node[below of=5](5') {$\mathscr S$}; 

\draw[->] ($(1.east)+(0,1pt)$) to node[above,la]{$\simeq$} ($(2.west)+(0,1pt)$);
\draw[->] ($(2.east)+(0,1pt)$) to node[above,la]{$\simeq$} ($(3.west)+(0,1pt)$);
\draw[->] ($(3.east)+(0,1pt)$) to node[above,la]{$\simeq$} (4);
\draw[->] (4) to node[above,la]{$\simeq$} (5);

\draw[d] ($(1'.east)+(0,1pt)$) to ($(2'.west)+(0,1pt)$);
\draw[d] ($(2'.east)+(0,1pt)$) to ($(3'.west)+(0,1pt)$);
\draw[->] ($(3'.east)+(0,1pt)$) to node[below,la]{$\simeq$} (4');
\draw[d] (4') to (5');

\draw[right hook->] (1) to node[left,la]{$[\mathrm{disc}]_\infty$} (1');
\draw[right hook->] (2) to (2');
\draw[right hook->] (3) to (3');
\draw[right hook->] (4) to (4');
\draw[right hook->] (5) to (5');
\end{tz}
% \[\begin{tikzcd}
% \set_\infty\arrow[r,"\simeq"]\arrow[d,"\mathrm{disc}_\infty" swap]&(\cL_?(\sset_{(\infty,0)}))_\infty\arrow[d,""]\arrow[r,"\simeq"]&\mathscr L_?((\sset_{(\infty,0)})_\infty)\arrow[d]\arrow[r,"\simeq"]&\mathscr L_?\mathscr S\arrow[d,""]\arrow[r,"\simeq"]&\mathscr Set\arrow[d,""]\\
% (\sset_{(\infty,0)})_\infty\arrow[r,equal]&(\sset_{(\infty,0)})_\infty\arrow[r,equal]&(\sset_{(\infty,0)})_\infty\arrow[r,"\simeq" swap]&\mathscr S\arrow[r,equal]&\mathscr S
% \end{tikzcd}\]
% \[\begin{tikzcd}
% \set_\infty\arrow[r,"\mathrm{disc_\infty}"]\arrow[d,"\simeq"]&\sset_{(\infty,0)}\arrow[d,equal]\\
% (\cL(\sset_{(\infty,0)}))_\infty\arrow[r]\arrow[d,"\simeq"]&\sset_{(\infty,0)}\arrow[d,equal]\\
% \mathscr L((\sset_{(\infty,0)}))_\infty\arrow[r]\arrow[d,"\simeq"]&\sset_{(\infty,0)}\arrow[d,"\simeq"]\\
% \mathscr L\mathscr S\arrow[r]\arrow[d,"\simeq"]&\mathscr S\arrow[d,equal]\\
% \mathscr Cat_{(0,0)}\arrow[r]\arrow[d,"\simeq"]&\mathscr S\arrow[d,equal]\\
% \mathscr Set\arrow[r]&\mathscr S,\\
% \end{tikzcd}\]
which commutes up to equivalence, using \cref{localization compatible} and \cite[Lem.~6.1.6(i)]{GH}.
\end{rmk}

%%%

\subsection{The \texorpdfstring{$1$}{1}-dimensional case}

\label{1dim}

The goal of this subsection is to show that the Quillen pair given by the ordinary nerve--categorification adjunction
\[c\colon \sset_{(\infty,1)}\rightleftarrows\cat\colon {\bf N}\]
is equivalent to the left Bousfield 
localization of the Joyal model structure $\sset_{(\infty,1)}$ with respect to a set $\Sigma\Lambda$ of maps. We also discuss in \cref{equivalencewithGH1} that this entails that the nerve embedding realizes $1$-categories as local $(\infty,1)$-categories with respect to the set of maps $\Sigma\Lambda$.

The following is a well-known fact, and of straightforward verification.

\begin{prop}
\label{nerverightquillen}
The functor ${\bf N}\colon\cat\to\sset_{(\infty,1)}$ is a right Quillen embedding.
\end{prop}

In particular, we have a Quillen reflection pair
\[
c\colon \sset_{(\infty,1)}\rightleftarrows\cat\colon {\bf N}.
\]

\begin{rmk}
The essential image of the functor $[\mathbf N]_\infty\colon[\cat]_\infty\to[\sset_{(\infty,1)}]_\infty$ is the full sub-$\infty$-category of $[\sset_{(\infty,1)}]_\infty$ generated by the locally homotopically discrete $(\infty,1)$-categories.
\end{rmk}

Recall that the (right-sided) suspension of simplicial sets defines a left adjoint functor $\Sigma\colon\sset\to\sset_{*,*}$. Given a simplicial set $X$, the suspension can be understood as the following pushout of simplicial sets.
 \begin{tz}
 \node[](1) {$X$}; 
 \node[below of=1](2) {$\Delta[0]$}; 
 \node[right of=1,xshift=.5cm](3) {$X\star\Delta[0]$}; 
 \node[below of=3](4) {$\Sigma X$};
 
 \draw[->] (1) to (2);
 \draw[->] (1) to (3); 
 \draw[->] (2) to (4);
 \draw[->] (3) to (4);
 
\pushout{4};
 \end{tz}
% \[\begin{tikzcd}
% X\arrow[r]\arrow[d]&X\star\Delta[0]\arrow[d]\\\Delta[0]\arrow[r]&\Sigma X.
% \end{tikzcd}\]

Recall from \cite{HirschhornOvercategories} that given any model category $\cM$, there is a model category $\cM_{*,*}$ of bipointed objects in $\cM$, in which fibrations, cofibrations, and weak equivalences are created by the forgetful functor $\cM_{*,*}\to\cM$.

The proof of the following could be adapted from \cite[Lemma 2.7]{ORfundamentalpushouts}, using ideas from \cite[Prop.\ 6.29]{JoyalVolumeII}.

\begin{prop}
\label{suspension1}
The suspension functor $\Sigma\colon\sset_{(\infty,0)}\to{(\sset_{(\infty,1)})}_{*,*}$ is a left Quillen functor.
\end{prop}

Since the model structure $\sset_{(\infty,1)}$ is combinatorial and left proper, the following model structure exists.

\begin{prop}
The category $\sset$ admits the left Bousfield localization $\cL_{\Sigma\Lambda}\sset_{(\infty,1)}$ of the Joyal model structure $\sset_{(\infty,1)}$ with respect to the set $\Sigma\Lambda$ of maps of the form
\begin{equation}
\label{suspensionofsphere}
   \Sigma\Delta[0]\hookrightarrow\Sigma S^k,\quad\text{for $k>0$}. 
\end{equation}
\end{prop}

So there is a Quillen reflection pair
\[
\mathrm{Id}\colon\sset_{(\infty,1)}\rightleftarrows \cL_{\Sigma\Lambda}\sset_{(\infty,1)}\colon \mathrm{Id}.
\]

To prove the desired result, we will show that the nerve functor induces a right Quillen equivalence ${\bf N}\colon\cat\to \cL_{\Sigma\Lambda}\sset_{(\infty,1)}$.

\begin{prop}
\label{Loc1Nerve}
The nerve functor ${\bf N}\colon\cat\to \cL_{\Sigma\Lambda}\sset_{(\infty,1)}$ defines a right Quillen embedding.
\end{prop}

\begin{rmk}
\label{catSigma}
For every simplicial set $X$ there is a natural isomorphism of categories
\[c\Sigma X \cong \Sigma \pi_0 cX \cong \Sigma\pi_0X.\]
\end{rmk}

\begin{proof}
[{Proof of \cref{Loc1Nerve}}]
By \cite[Prop.\ 3.3.18]{hirschhorn} and \cref{nerverightquillen}, it is sufficient to show that $c$ sends all maps from \eqref{suspensionofsphere} to (weak) equivalences in $\cat$.

Let $k>0$.
The functor $c$ sends the map
\[\Sigma\Delta[0]\hookrightarrow\Sigma S^k\]
to the map
\[c\Sigma\Delta[0]\hookrightarrow c\Sigma S^k,\]
which is
by \cref{catSigma}
\[\Sigma\pi_0\Delta[0]\hookrightarrow\Sigma\pi_0S^k,\]
which is the identity isomorphism at $\Sigma[0]$.
This concludes the proof that the desired functor is right Quillen.

The fact that it is a right Quillen embedding follows directly from \cref{nerverightquillen} as the derived counits of $\bfN\colon \cat\to \sset_{(\infty,1)}$ and $\bfN\colon \cat\to \cL_{\Sigma\Lambda}\sset_{(\infty,1)}$ coincide at a fibrant object in $\cL_{\Sigma \Lambda} \sset_{(\infty,1)}$.
\end{proof}

\begin{prop}
\label{suspension1localization}
The suspension functor $\Sigma\colon\cL_\Lambda\sset_{(\infty,0)}\to(\cL_{\Sigma\Lambda}{(\sset_{(\infty,1)})})_{*,*}$ is a left Quillen functor.
\end{prop}

\begin{proof}
As an instance of \cite[Theorem 3.3.20]{hirschhorn} combined with the fact that every object is cofibration in $\sset_{(\infty,0)}$, we know that \[\Sigma\colon\cL_\Lambda\sset_{(\infty,0)}\to\cL_{\Sigma\Lambda}({(\sset_{(\infty,1)})}_{*,*})\]is a left Quillen functor. Further, since left Bousfield localizations commute with taking bipointed model structures, the model structures
\[\cL_{\Sigma\Lambda}({(\sset_{(\infty,1)})}_{*,*})=(\cL_{\Sigma\Lambda}{(\sset_{(\infty,1)})})_{*,*}\]
are equal. This concludes the proof.
\end{proof}

The functor $\Sigma\colon\sset\to\sset_{*,*}$ admits a right adjoint $\mathrm{Hom}^R\colon\sset_{*,*}\to\sset$, used e.g.~in \cite[\S1.2.2]{htt}. For any simplicial set $X$ with given vertices $x$ and $y$
we write $X(x,y)\coloneqq \mathrm{Hom}^R_X(x,y)$.

\begin{rmk} \label{rem:lem12}
The following facts are of straightforward verifications. The first one uses the explicit description from e.g.~\cite[Prop.\ 4.12]{BoardmanVogt} of the category $cX$ in the case of $X$ being a quasi-category; see also \cite[Prop.\ 1.11]{JoyalVolumeII}.
\begin{enumerate}[leftmargin=*]
    \item For any quasi-category $X$ with vertices $x$ and $y$ there is a bijection
\[\pi_0(X(x,y))\cong(cX)(x,y).\]
    \item For any category $\cC$ there is an isomorphism of simplicial sets
\[\mathrm{disc}(\cC(x,y))
\cong({\bf N}\cC)(x,y).\]
\end{enumerate}
\end{rmk}

\begin{thm}
\label{QE1}
The nerve functor ${\bf N}\colon\cat\to \cL_{\Sigma\Lambda}\sset_{(\infty,1)}$ defines a right Quillen equivalence.
\end{thm}

\begin{proof}
By \cref{Loc1Nerve}, it remains to prove that the component of the derived unit at every object $X$ in $\cL_{\Sigma\Lambda}\sset_{(\infty,1)}$ is a weak equivalence. We do this by first proving it in the case of $X$ being fibrant in $\cL_{\Sigma\Lambda}\sset_{(\infty,1)}$, and then treating the general case.

Assume that $X$ is fibrant in $\cL_{\Sigma\Lambda}\sset_{(\infty,1)}$.
Then for any vertices $x$ and~$y$ in $X$, the tuple $(X,x,y)$ is fibrant in $(\cL_{\Sigma\Lambda}(\sset_{(\infty,1)}))_{*,*}$ so $X(x,y)$ is fibrant in $\cL_\Lambda\sset_{(\infty,0)}$ by \cref{suspension1localization}.

By \cref{QE0}, the (derived) unit at $X(x,y)$ is a weak equivalence in $\cL_\Lambda \sset_{(\infty,0)}$
\begin{align*}
X(x,y)  & \simeq \mathrm{disc}(\pi_0(X(x,y))&\\
   &\cong  \mathrm{disc}((cX)(x,y))&\text{\cref{rem:lem12}(1)}\\
  & \cong ({\bf N}cX)(x,y).&\text{\cref{rem:lem12}(2)}
\end{align*}
between fibrant objects. Hence, it is already a weak equivalence in $\sset_{(\infty,0)}$.

This weak equivalence
\[X(x,y)\to({\bf N}cX)(x,y)\]
is precisely the map  obtained by taking $\Hom^R$ of the (derived) unit of $(X,x,y)$.
This means that the (derived) unit of $X$
\[X\to {\bf N}cX\]
is locally a weak equivalence of simplicial sets, as well as a bijection on objects.
By the fundamental theorem of $(\infty,1)$-categories, originally due to Joyal \cite{JoyalVolumeII} and recalled e.g.~in \cite[Thm~3.9.7]{CisinskiBook},
we deduce that the (derived) unit is then a weak equivalence in $\sset_{(\infty,1)}$, so in particular in the localization $\cL_{\Sigma\Lambda}\sset_{(\infty,1)}$
as desired.

Now if $X$ is more generally any (cofibrant) simplicial set, we consider a fibrant replacement $X^{\mathrm{fib}}$ in  $\cL_{\Sigma\Lambda}\sset_{(\infty,1)}$ and the following naturality diagram.
\begin{tz}
\node[](1) {$X$}; 
\node[right of=1,xshift=.7cm](2) {${\bf N}cX$}; 
\node[below of=1](3) {$X^\mathrm{fib}$}; 
\node[below of=2](4) {${\bf N}c(X^\mathrm{fib})$};

\draw[->] (1) to (2);
\draw[->] (1) to (3);
\draw[->] (2) to (4);
\draw[->] (3) to (4); 
\end{tz}
% \[\begin{tikzcd}
%     X\arrow[d]\arrow[r]&{\bf N}c X\arrow[d]\\
%   X^{\mathrm{fib}}\arrow[r]&{\bf N}((c(X^{\mathrm{fib}}))).\\ 
% \end{tikzcd}\]
% \[\begin{tikzcd}
%     X\arrow[d]\arrow[r]&N((c X)^{\mathrm{fib}})\arrow[d]\\
%   X^{\mathrm{fib}}\arrow[r]&N((c(X^{\mathrm{fib}}))^{\mathrm{fib}})\\ 
% \end{tikzcd}\]
Here, the left vertical map is a weak equivalence in $\cL_{\Sigma\Lambda}\sset_{(\infty,1)}$ by construction, the right vertical map is a weak equivalence because both $\bf N$ and $c$ are homotopical, and the bottom horizontal arrow is a weak equivalence by the case that we already treated. It follows by $2$-out-of-$3$ that the top horizontal map, which is the (derived) unit of $X$, is a weak equivalence, as desired.
\end{proof}

\begin{rmk}
\label{equivalencewithGH1}
We know -- and it is also mentioned in \cite[Ch.~3]{htt} -- that the underlying $\infty$-category of the Joyal model structure $\sset_{(\infty,1)}$ models the established $\infty$-category $\mathscr Cat_{(\infty,1)}$ of $\infty$-categories, so there exists an equivalence of $\infty$-categories
\begin{equation}
    \label{correctinfinitycategories}
  [\sset_{(\infty,1)})]_\infty\simeq \mathscr Cat_{(\infty,1)}.
\end{equation}
Any such equivalence can be used to construct a specific
equivalence of $\infty$-categories
\begin{align}\label{correctinfinitycategories2}
    [\cat]_\infty&\simeq[\cL_{\Sigma\Lambda}\sset_{(\infty,1)}]_\infty&\text{\cref{QE1}} \nonumber\\
    &\simeq\mathscr L_{\Sigma\Lambda}[\sset_{(\infty,1)}]_\infty&\text{\cref{localization compatible}}\\
    &\simeq\mathscr L_{\Sigma\Lambda}\mathscr Cat_{(\infty,1)}&\text{\eqref{correctinfinitycategories}}\nonumber\\
     &\simeq\mathscr Cat_1&\text{\cite[Lem.~6.1.7(v)]{GH}}\nonumber
\end{align}
     between the established $\infty$-category $\mathscr Cat_1$ of categories and the underlying $\infty$-category of the model structure $\cat$ on categories. Via the chosen identifications \eqref{correctinfinitycategories} and \eqref{correctinfinitycategories2},
we see that the functor $[{\mathbf N}]_\infty\colon[\cat]_\infty\to[\sset_{(\infty,1)}]_\infty$ and the canonical inclusion functor $\mathscr Cat_1\hookrightarrow\mathscr Cat_{(\infty,1)}$ from \cite[Lem.~6.1.7(v)]{GH} -- used with $n=1$ -- are equivalent. Indeed, this is witnessed by the following diagram of $\infty$-categories
\begin{tz}
\node[](1) {$[\cat]_\infty$}; 
\node[below of=1](1') {$[\sset_{(\infty,1)}]_\infty$}; 
\node[right of=1,xshift=1.2cm](2) {$[\cL_{\Sigma\Lambda}\sset_{(\infty,1)}]_\infty$}; 
\node[below of=2](2') {$[\sset_{(\infty,1)}]_\infty$}; 
\node[right of=2,xshift=2cm](3) {$\mathscr L_{\Sigma\Lambda}[\sset_{(\infty,1)}]_\infty$}; 
\node[below of=3](3') {$[\sset_{(\infty,1)}]_\infty$}; 
\node[right of=3,xshift=1.7cm](4) {$\mathscr L_{\Sigma\Lambda} \mathscr Cat_{(\infty,1)}$}; 
\node[below of=4](4') {$\mathscr Cat_{(\infty,1)}$}; 
\node[right of=4,xshift=.8cm](5) {$\mathscr Cat_1$}; 
\node[below of=5](5') {$\mathscr Cat_{(\infty,1)}$}; 

\draw[->] (1) to node[above,la]{$\simeq$} (2);
\draw[->] (2) to node[above,la]{$\simeq$} (3);
\draw[->] (3) to node[above,la]{$\simeq$} (4);
\draw[->] (4) to node[above,la]{$\simeq$} (5);

\draw[d] (1') to (2');
\draw[d] (2') to (3');
\draw[->] (3') to node[below,la]{$\simeq$} (4');
\draw[d] (4') to (5');

\draw[right hook->] (1) to node[left,la]{$[\bfN]_\infty$} (1');
\draw[right hook->] (2) to (2');
\draw[right hook->] (3) to (3');
\draw[right hook->] (4) to (4');
\draw[right hook->] (5) to (5');
\end{tz}
% \[\begin{tikzcd}
% \cat_\infty\arrow[r,"\simeq"]\arrow[d,"\mathbf{N}_\infty" swap]&(\cL_{\Sigma\Lambda}(\sset_{(\infty,1)}))_\infty\arrow[d,""]\arrow[r,"\simeq"]&\mathscr L_{\Sigma\Lambda}((\sset_{(\infty,1)})_\infty)\arrow[d]\arrow[r,"\simeq"]&\mathscr L_{\Sigma\Lambda}(\mathscr Cat_{(\infty,1)})\arrow[d,""]\arrow[r,"\simeq"]&\mathscr Cat_1\arrow[d,""]\\
% (\sset_{(\infty,1)})_\infty\arrow[r,equal]&(\sset_{(\infty,1)})_\infty\arrow[r,equal]&(\sset_{(\infty,1)})_\infty\arrow[r,"\simeq" swap]&\mathscr Cat_{(\infty,1)}\arrow[r,equal]&\mathscr Cat_{(\infty,1)}
% \end{tikzcd}\]
% \[\begin{tikzcd}
% \cat_\infty\arrow[r,"\mathrm{\mathbf N_\infty}"]\arrow[d,"\simeq"]&\sset_{(\infty,1)}\arrow[d,equal]\\
% (\cL(\sset_{(\infty,1)}))_\infty\arrow[r]\arrow[d,"\simeq"]&\sset_{(\infty,1)}\arrow[d,equal]\\
% \mathscr L((\sset_{(\infty,1)}))_\infty\arrow[r]\arrow[d,"\simeq"]&\sset_{(\infty,1)}\arrow[d,"\simeq"]\\
% \mathscr L\mathscr Cat_\infty\arrow[r]\arrow[d,"\simeq"]&\mathscr Cat_\infty\arrow[d,equal]\\
% \mathscr Cat_{(1,1)}\arrow[r]\arrow[d,"\simeq"]&\mathscr Cat_\infty\arrow[d,equal]\\
% \mathscr Cat_\infty\arrow[r]&\mathscr Cat_\infty\\
% \end{tikzcd}\]
which commutes up to equivalence, using \cref{localization compatible}, \cite[Lem.~6.1.9]{GH}, and \cite[Lem.~6.1.7(i)]{GH}.
\end{rmk}

\subsection{The \texorpdfstring{$2$}{2}-dimensional case}
\label{2dim}

The goal of this subsection is to show that the Quillen reflection pair from \cref{EnrichedNerves}(1)
\[c_*\colon \vcat{\sset_{(\infty,1)}}\rightleftarrows\twocat\colon {\bf N}_*.\]
is equivalent to the left Bousfield localization of $\vcat{\sset_{(\infty,1)}}$ with respect to a set $\Sigma^2\Lambda$ of maps. We also discuss in \cref{equivalencewithGH2} that this entails that the nerve embedding realizes $2$-categories as local $(\infty,2)$-categories with respect to the set of maps $\Sigma^2\Lambda$.

\begin{rmk}
The essential image of the functor $[\mathbf N_*]_\infty\colon[\twocat]_\infty\to[\vcat{\sset_{(\infty,1)}}]_\infty$ is the full sub-$\infty$-category of $[\vcat{\sset_{(\infty,1)}}]_\infty$ generated by the $(\infty,2)$-categories that are locally equivalent to $1$-categories.
\end{rmk}

Recall that there is a suspension functor $\Sigma\colon\sset\to(\vcat{\sset})_{*,*}$ which is a left adjoint. Given a simplicial set $X$, the simplicial category $\Sigma X$ has two objects and a single non-trivial hom-simplicial set given by $X$. The following is briefly discussed e.g.~as \cite[Lem.~4.1.5]{HORR}.

\begin{prop}
\label{suspension12}
The suspension functor $\Sigma\colon\sset_{(\infty,1)}\to(\vcat{\sset_{(\infty,1)}})_{*,*}$ is a left Quillen functor.
\end{prop}

We consider the composite functor
\[\Sigma^2\colon\sset\xrightarrow{\Sigma}\sset_{*,*}\xrightarrow{U}\sset\xrightarrow{\Sigma}(\vcat{\sset})_{*,*}. \]

\begin{prop}
\label{Sigma2}
The $2$-fold suspension functor $\Sigma^2\colon\sset_{(\infty,0)}\to(\vcat{\sset_{(\infty,1)}})_{*,*}$ is a left Quillen functor.
\end{prop}

\begin{proof}
It is a composite of the left Quillen (hence homotopical) functor
\[\Sigma\colon\sset_{(\infty,0)}\to(\sset_{(\infty,1)})_{*,*}\]
from \cref{suspension1} with the homotopical functor
\[U\colon (\sset_{(\infty,1)})_{*,*}\to\sset_{(\infty,1)},\]
which just forgets the two base points, and with the left Quillen (hence homotopical) functor
\[\Sigma\colon\sset_{(\infty,1)}\to(\vcat{\sset_{(\infty,1)}})_{*,*}\]
from \cref{suspension12}.
\end{proof}

\begin{rmk}
\label{suspensioncomparison}
Let $\cV=\sset_{(\infty,1)}$, so that in particular $\mathscr V=[\cV]_\infty=[\sset_{(\infty,1)}]_\infty\simeq\mathscr Cat_{\infty}$.
The suspension functor from \cref{suspension12} is a left Quillen functor, and induces a functor of $\infty$-categories
\begin{equation}
    \label{oursusp}
[\Sigma]_\infty\colon[\cV]_\infty\to[(\vcat{\cV})_{*,*}]_\infty.
\end{equation}
In \cite[Def.\ 4.3.21]{GH} Gepner--Haugseng consider a functor
\begin{equation}
    \label{GHsusp}
\mathscr V
\to\mathscr{C}at_{\mathscr{V}}^{\{0,1\}}.\end{equation}
Here, $\mathscr{C}at_{\mathscr{V}}^{\{0,1\}}$ denotes the $\infty$-category of $\infty$-categories enriched over $\mathscr V$ with fixed set of objects $\{0,1\}$, as defined in \cite[Def.~5.4.3]{GH}. As shown in \cite[\textsection5]{HaugsengRectification}, this $\infty$-category can be realized as the underlying $\infty$-category $[\vcat{\cV}^{\{0,1\}}]_\infty\simeq\mathscr{C}at_{\mathscr{V}}^{\{0,1\}}$ of the model category $\vcat{\cV}^{\{0,1\}}$ of $\cV$-categories with set of objects $\{0,1\}$, considered in \cite[Lemma~3.20]{HaugsengRectification}.

Via the canonical map
\begin{equation}
    \label{comparisonbipointed}
[(\vcat{\cV})_{*,*}]_\infty\to [\vcat{\cV}^{\{0,1\}}]_\infty\simeq\mathscr{C}at_{\mathscr{V}}^{\{0,1\}}
\end{equation}
we will see in \cref{comparisonsuspensionlemma} that the two functors \eqref{oursusp} and \eqref{GHsusp} are compatible, as they fit in a diagram of $\infty$-categories that commutes up to equivalence.

\end{rmk}

\begin{prop}
\label{comparisonsuspensionlemma}
There is a diagram of $\infty$-categories
\begin{tz}
\node[](1) {$[(\vcat{\cV})_{*,*}]_\infty$}; 
\node[above of=1](2) {$\mathscr{C}at_{\mathscr{V}}^{\{0,1\}}$}; 
\node[left of=1,xshift=-1.3cm](3) {$[\cV]_\infty$}; 
\node[above of=3](4) {$\mathscr V$}; 

\draw[->] (3) to node[below,la]{$[\Sigma]_\infty$} (1);
\draw[->] (2) to (1);
\draw[d] (3) to (4);
\draw[->] (4) to (2);
\end{tz}
that commutes up to equivalence, where the functors involved are those of \eqref{oursusp}, \eqref{GHsusp}, and \eqref{comparisonbipointed}.
\end{prop}

\begin{proof}
Each of the functors of $\infty$-categories involved in the diagram admits a right adjoint. 
We prove that the diagram of right adjoints commutes up to equivalence:
\begin{tz}
\node[](1) {$[(\vcat{\cV})_{*,*}]_\infty$}; 
\node[below of=1](2) {$\mathscr{C}at_{\mathscr{V}}^{\{0,1\}}$}; 
\node[right of=1,xshift=1.3cm](3) {$[\cV]_\infty$}; 
\node[below of=3](4) {$\mathscr V$}; 

\draw[->] (1) to node[above,la]{$[\Hom]_\infty$} (3);
\draw[->] (1) to (2);
\draw[d] (3) to (4);
\draw[->] (2) to node[below,la]{$\mathscr{H}om$} (4);
\end{tz}
Building the desired commutative diagram out of smaller ones requires several ingredients, for which we provide references for the interested reader.
The diagram is:
\begin{tz}
\node[](1) {$[(\vcat{\cV})_{*,*}]_\infty$}; 
\node[below of=1](2) {$[\vcat{\cV}^{\{0,1\}}]_\infty$}; 
\node[below of=2](3) {$\mathscr{C}at_{\mathscr{V}}^{\{0,1\}}$};
\node[below of=3](4) {$\Alg_{\Delta_{\{0,1\}}^{\op}}(\mathscr V)$};
\node[below of=4,yshift=2pt](5) {$\Fun_{\Delta^{\op}}(\Delta_{\{0,1\}}^{\op}, \mathscr V^{\otimes})$}; 

\node[right of=1,xshift=3cm,yshift=1pt](1') {$[\cV^{\{0,1\}\times\{0,1\}}]_\infty$}; 
\node[below of=1'](2') {$[\cV]_\infty^{\{0,1\}\times\{0,1\}}$}; 
\node[below of=2'](3') {$\mathscr{V}^{\{0,1\}\times\{0,1\}}$};
\node[below of=3',yshift=-1pt](4') {$\mathrm{Alg}_{(\Delta_{\{0,1\}}^{\op})_{\mathrm{triv}}}(\mathscr V)$};
\node[below of=4',yshift=2pt](5') {$\Fun_{\Delta^{\op}}((\Delta_{\{0,1\}}^{\op})_{\mathrm{triv}}, \mathscr V^{\otimes})$};

\node[right of=1',xshift=3cm,yshift=-1pt](1'') {$[\cV]_\infty$};
\node[below of=1''](2'') {$[\cV]_\infty$}; 
\node[below of=2''](3'') {$\mathscr{V}$};
\node[below of=3'',yshift=2pt](4'') {$\Fun(\{0,1\}\times \{0,1\}, \mathscr{V})$};
\node[below of=4''](5'') {$\Fun((\Delta_{\{0,1\}}^{\op})_{\mathrm{triv}})_{[1]}, \mathscr V)$};

\setcounter{region}{0};
\refstepcounter{region}
\node[scale=0.9] at ($(3)+(1cm,-.6cm)$) {(\theregion)\label{regioneta}};
\refstepcounter{region}
\node[scale=0.9] at ($(4')+(-1.3cm,.6cm)$) {(\theregion)\label{regionvtau}};
\refstepcounter{region}
\node[scale=0.9] at ($(5)!0.5!(5')+(0,.6cm)$) {(\theregion)\label{regiontau}};

\draw[->,rounded corners] (3.west) to ($(3)-(2.2cm,0)$) to ($(5)-(2.2cm,.7
5cm)$) to node[over]{$\mathscr{H}om$} ($(5'')+(2.2cm,-.75cm)$) to ($(3'')+(2.2cm,0)$) to (3''.east);

\draw[->,rounded corners] (1.west) to ($(1)-(2.2cm,0)$) to ($(1)+(-2.2cm,.7
5cm)$) to node[over]{$[\Hom]_\infty$} ($(1'')+(2.2cm,.75cm)$) to ($(1'')+(2.2cm,0)$) to (1''.east);

\draw[->] (1) to (2); 
\draw[->] (2) to node[left,la]{$\simeq$} (3); 
\draw[->] (3) to node[left,la]{$\eta$} (4); 
\draw[right hook->] (4) to (5);

\draw[->] (1') to node[left,la]{$\simeq$} (2'); 
\draw[d] (2') to (3'); 
\draw[->] (3') to node[left,la]{$\simeq$} (4'); 
\draw[right hook->] (4') to (5');

\draw[d] (1'') to (2''); 
\draw[d] (2'') to (3'');
\draw[->] (4'') to node[right,la]{$\ev_{(0,1)}$} (3'');
\draw[->] (5'') to node[right,la]{$\cong$} (4'');

\draw[->] (1) to ($(1'.west)+(0,-1pt)$); 
\draw[->] ($(1'.east)+(0,-1pt)$) to node[above,la]{$\ev_{(0,1)}$} (1'');

\draw[->] (2) to ($(2'.west)+(0,-1pt)$); 
\draw[->] ($(2'.east)+(0,-1pt)$) to node[above,la]{$\ev_{(0,1)}$} (2'');

\draw[->] (3) to ($(3'.west)+(0,-1pt)$); 
\draw[->] ($(3'.east)+(0,-1pt)$) to node[above,la]{$\ev_{(0,1)}$} (3'');

\draw[->] ($(4.east)+(0,2pt)$) to node[below,la]{$\tau_{\Delta_{\{0,1\}}^{\op}}^*$} ($(4'.west)+(0,2pt)$); 
\draw[->] (4) to node[above,la]{$V$} (3');
\draw[d] (3') to (4''); 
\draw[->] (4') to node[above,la]{$\simeq$} (5'');

\draw[->] ($(5.east)+(0,1pt)$) to ($(5'.west)+(0,1pt)$); 
\draw[->] ($(5'.east)+(0,1pt)$) to ($(5''.west)+(0,1pt)$);
\end{tz}
% \[\begin{tikzcd}
% {[(\vcat{\cV})_{*,*}]_\infty}\arrow[r]&{[\cV^{\{0,1\}\times\{0,1\}}]_\infty}\arrow[d,"\simeq"]\arrow[rd,"\ev_{(0,1)}"]\\
% {[\vcat{\cV}^{\{0,1\}}]_\infty}\arrow[r]\arrow[d,"\simeq"]\arrow[u]&{{[\cV]_\infty^{\{0,1\}\times\{0,1\}}}}\arrow[d,"\simeq"]\arrow[r,"\ev_{(0,1)}"]&{[\cV]_\infty}\arrow[d, "\simeq"]\\
% {\mathscr{C}at_{\mathscr{V}}}^{\{0,1\}}\arrow[r]\arrow[d,"\eta"]&\mathscr{V}^{\{0,1\}\times\{0,1\}}\arrow[rd,"\simeq"]\arrow[d,"\simeq"]\arrow[r,"\ev_{(0,1)}"]&\mathscr{V}\\
% \Alg_{\Delta_{\{0,1\}}^{\op}}(\mathscr V)\arrow[ru,"V"]\arrow[d,hook]\arrow[r,"\tau_{\Delta_{\{0,1\}}^{\op}}^*" swap]&\mathrm{Alg}_{(\Delta_{\{0,1\}}^{\op})_{\mathrm{triv}}}(\mathscr V)\arrow[d,hook]\arrow[dr, "\simeq"]&\Fun(\{0,1\}\times\{0,1\},\mathscr{V})\arrow[u]\arrow[d,"\simeq"]\\
% \Fun_{\Delta^{\op}}(\Delta_{\{0,1\}}^{\op}, \mathscr V)\arrow[r]&\Fun((\Delta_{\{0,1\}}^{\op})_{[1]}, \mathscr V)\arrow[r]&\Fun({(\Delta_{\{0,1\}}^{\op})_{[1]}}, \mathscr V)
% \arrow[u,"\ev_{(0,1)}"]\end{tikzcd}\]
The $\infty$-categories featuring in the diagram are the following:
\begin{itemize}[leftmargin=*]
\item $(\vcat{\cV})_{*,*}$ is the bipointed model structure, obtained as an instance of \cite{HirschhornOvercategories} applied to the model structure from \cref{EnrichedCatMS}.
\item $\vcat{\cV}^{\{0,1\}}$ is the model category of $\cV$-categories with set of objects $\{0,1\}$, considered in \cite[Lem.\ 3.20]{HaugsengRectification}.
\item $\cV^{\{0,1\}\times\{0,1\}}$ is the category of functors endowed with the injective model structure.
    \item $\Alg_{\Delta_{\{0,1\}}^{\op}}(\mathscr V)$ is an instance of \cite[\textsection 1.2]{GH} with the non-symmetric $\infty$-operad $\Delta_{\{0,1\}}^{\op}$ from \cite[Def.\ 2.8]{HaugsengRectification}. 
    \item $\Alg_{(\Delta_{\{0,1\}}^{\op})_{\mathrm{triv}}}(\mathscr V)$ is an instance of \cite[\textsection 1.2]{GH} with the non-symmetric $\infty$-operad $(\Delta_{\{0,1\}}^{\op})_{\mathrm{triv}}$ from \cite[Def.~3.4.1]{GH}.
  \item $(\Delta_{\{0,1\}}^{\op})_{[1]}$ is the fiber at $[1]$,  which is an object of $\Delta^{\op}$ of the map $\Delta_{\{0,1\}}^{\op}\to\Delta^{\op}$.
\end{itemize}
The functors of $\infty$-categories featuring in the diagram are the following:
\begin{itemize}[leftmargin=*, label=$\diamond$]
\item The functor $\ev_{(0,1)}$ is given by evaluation at the object $(0,1)\in \{0,1\}\times \{0,1\}$.
    \item The functor $V$ is from \cite[Proof of Prop.\ 5.2]{GH}.
    \item The functor $\tau_{\Delta_{\{0,1\}}^{\op}}^*$ is the one considered in \cite[\textsection A.4,~\textsection3.4]{GH}.
    \item The functor $\eta$ is constructed on the level of model categories in \cite[Proof of Prop.\ 5.2]{HaugsengRectification}, and the functor induced at the level of $\infty$-categories is further described in \cite[Def.\ 4.3.1, Prop.\ 5.4.4]{GH}.
\end{itemize}
We address the commutativity of each of the labeled regions as follows.
 \begin{itemize}[leftmargin=*, label=$\star$]
    \item The fact that the region (\ref{regioneta}) commutes is addressed as \cite[Proof of Prop.\ 5.2]{HaugsengRectification}.
    \item The fact that the region (\ref{regionvtau}) commutes is addressed as a combination of \cite[Lem.\ 3.20]{HaugsengRectification}, \cite[\S3.4]{GH} and \cite[\S A.4]{GH}.
    \item The fact that the region (\ref{regiontau}) commutes is addressed in \cite[\textsection A.4,~A.5]{GH}.

 \end{itemize}
This concludes the proof.
\end{proof}

\begin{prop}
The category $\vcat{\sset}$ has
the left Bousfield localization $\cL_{\Sigma^2\Lambda}\vcat{\sset_{(\infty,1)}}$ of the model structure $\vcat{\sset_{(\infty,1)}}$ with respect to the set $\Sigma^2\Lambda$ of maps of the form
\begin{equation}
    \label{suspensionofsphere2}
    \Sigma^2\Delta[0]\hookrightarrow\Sigma^2 S^k, \quad\text{for $k>0$}.
\end{equation}
\end{prop}

\begin{proof}
The Bousfield localization exists because the model category $\sset_{(\infty,1)}$ is combinatorial and left proper by \cite[Prop.~A.3.2.4]{htt}.
\end{proof}

So there is a Quillen reflection pair

\[
\mathrm{Id}\colon\vcat{\sset_{(\infty,1)}}\rightleftarrows \cL_{\Sigma^2\Lambda}\vcat{\sset_{(\infty,1)}}\colon \mathrm{Id}.
\]

To prove the desired result, we will show that the nerve functor induces a right Quillen equivalence ${\bf N}_*\colon\twocat\to \cL_{\Sigma^2\Lambda}\vcat{\sset_{(\infty,1)}}$. First, we prove the following.

\begin{prop} \label{nervestarrightQuillen}
The nerve functor ${\bf N}_*\colon\twocat\to \cL_{\Sigma^2\Lambda}\vcat{\sset_{(\infty,1)}}$ defines a right Quillen embedding.
\end{prop}

\begin{proof}
By \cite[Prop.\ 3.3.18]{hirschhorn} and \cref{EnrichedNerves}(1), it is sufficient to show that $c_*$ sends all elementary maps from \eqref{suspensionofsphere2} to biequivalences of $2$-categories.

Let $k>0$. The functor $c_*$ sends the map
\[\Sigma^2\Delta[0]\hookrightarrow\Sigma^2 S^k\]
to the map
\[c_*\Sigma^2\Delta[0]\hookrightarrow c_*\Sigma^2 S^k,\]
which is the map
\[\Sigma c\Sigma\Delta[0]\hookrightarrow \Sigma c\Sigma S^k,\]
which is the map
\[\Sigma^2\pi_0\Delta[0]\hookrightarrow \Sigma^2\pi_0S^k,\]
which is the identity at $\Sigma^2[0]$. This concludes the proof that the desired functor is right Quillen.

The fact that it is a right Quillen embedding follows directly from \cref{EnrichedNerves}(1) as the derived counits of $\bfN_*\colon \twocat\to \vcat{\sset_{(\infty,1)}}$ and $\bfN_*\colon \twocat\to \cL_{\Sigma^2\Lambda}\vcat{\sset_{(\infty,1)}}$ coincide at a fibrant object in $\cL_{\Sigma^2 \Lambda} \vcat{\sset_{(\infty,1)}}$.
\end{proof}

\begin{prop}
\label{suspension1localization2}
The suspension functor $\Sigma\colon\cL_{\Sigma\Lambda}\sset_{(\infty,1)}\to(\cL_{\Sigma^2\Lambda}{(\vcat{\sset_{(\infty,1)}}}))_{*,*}$ is a left Quillen functor.
\end{prop}

\begin{proof}
As an instance of \cite[Theorem 3.3.20]{hirschhorn} combined with the fact that every object is cofibration in $\sset_{(\infty,1)}$, we know that \[\Sigma\colon\cL_{\Sigma\Lambda}\sset_{(\infty,1)}\to\cL_{\Sigma^2\Lambda}({(\vcat{\sset_{(\infty,1)}})}_{*,*})\]is a left Quillen functor. Further, since left Bousfield localizations commute with taking bipointed model structures, the model structures 
\[\cL_{\Sigma^2\Lambda}({(\vcat{\sset_{(\infty,1)}})}_{*,*})=(\cL_{\Sigma^2\Lambda}{(\vcat{\sset_{(\infty,1)}})})_{*,*}\]
are equal. This concludes the proof.
\end{proof}

\begin{lem}
\label{c*homotopical}
The functor $c_*\colon\vcat{\sset_{(\infty,1)}}\to\twocat$ is homotopical.
\end{lem}

\begin{rmk} \label{DwyerKaneq}
A map $f\colon\cQ\to\cQ'$ is a weak equivalence in $\vcat{\sset_{(\infty,1)}}$ if and only if the following are satisfied.
\begin{enumerate}[leftmargin=*]
    \item The map $f$ is \emph{essentially surjective up to equivalence}; namely it induces an essentially surjective functor
    \[\tau_*f\colon \tau_*\cQ\to \tau_*\cQ',\]
    where $\tau_*\colon \vcat{\sset}\to \cat$ is the base-change functor along Joyal's functor $\tau\colon \sset\to \set$ from \cite[\textsection1]{JoyalVolumeII} given by the composite
    \[ \sset\xrightarrow{c}\cat\xrightarrow{\core}\cG pd\xrightarrow{\pi_0}\set. \]
\item The map $f$ is \emph{a local weak equivalence}; namely it induces a weak equivalence in $\sset_{(\infty,1)}$
\[f\colon\cQ(x,y)\to\cQ'(f(x),f(y))\]
for any objects $x$ and $y$ in $\cQ$.
\end{enumerate}
\end{rmk}

\begin{proof}[Proof of \cref{c*homotopical}]
Given a weak equivalence $f\colon \cQ\to \cQ'$ in $\vcat{\sset_{(\infty,1)}}$, we have a weak equivalence in $\sset_{(\infty,1)}$ 
\[\cQ(x,y)\to \cQ'(x,y),\]
for any objects $x$ and $y$ in $\cQ$, by \cref{DwyerKaneq}(2). Then, since $c$ is homotopical, there is an induced equivalence of categories
\[(c_*\cQ)(x,y)=c\cQ(x,y)\to c\cQ'(x,y)=(c_*\cQ')(x,y).\]
Moreover, by \cref{DwyerKaneq}(1) the functor
\[ (\pi_0)_*(\core)_* c_*\cQ=\tau_*\cQ\to \tau_* \cQ'=(\pi_0)_*(\core)_* c_*\cQ' \]
is essentially surjective on objects. Hence we obtain that the $2$-functor
\[c_*\cQ\to c_*\cQ'\]
is a weak equivalence in $\twocat$, as desired.
\end{proof}

\begin{thm}
\label{QE2}
The nerve functor ${\bf N}_*\colon\twocat\to \cL_{\Sigma^2\Lambda}\vcat{\sset_{(\infty,1)}}$ defines a right Quillen equivalence.
\end{thm}

\begin{proof}
By \cref{nervestarrightQuillen}, it remains to prove that the component of the derived unit at every object $\cQ$ in $\cL_{\Sigma^2\Lambda}\vcat{\sset_{(\infty,1)}}$ is a weak equivalence. We do this by first proving it in the case of $\cQ$ being fibrant in $\cL_{\Sigma^2\Lambda}\vcat{\sset_{(\infty,1)}}$, and then treat the general case.

Assume that $\cQ$ is fibrant in $\cL_{\Sigma^2\Lambda}\vcat{\sset_{(\infty,1)}}$. For any vertices $x$ and~$y$ the tuple $(\cQ,x,y)$ is fibrant in $(\cL_{\Sigma^2\Lambda}\vcat{\sset_{(\infty,1)}})_{*,*}$ so $\cQ(x,y)$ is fibrant in $\cL_{\Sigma\Lambda}\sset_{(\infty,1)}$ by \cref{suspension1localization2}.

By \cref{QE1}, the (derived) unit at $\cQ(x,y)$ is a weak equivalence in $\cL_\Lambda \sset_{(\infty,1)}$
\begin{align*}
\cQ(x,y)  & \simeq {\bf N}(c(\cQ(x,y))&\\
   &\cong  {\bf N}_*((c_*\cQ)(x,y))&\text{\cref{rem:lem12}(1)}\\
  & \cong ({\bf N}_*c_*\cQ)(x,y).&\text{\cref{rem:lem12}(2)}
\end{align*}
between fibrant objects.
Hence, it is already a weak equivalence in $\sset_{(\infty,1)}$.

This weak equivalence
\[\cQ(x,y)\to({\bf N}_*c_*\cQ)(x,y)\]
is precisely the one obtained by taking $\Hom$ of the (derived) unit of $(\cQ,x,y)$.
This means that the (derived) unit of $\cQ$
\[\cQ\to {\bf N}_*c_*\cQ\]
is locally a weak equivalence in $\sset_{(\infty,1)}$, as well as a bijection on objects.
By \cref{DwyerKaneq}
we deduce that the (derived) unit is then a weak equivalence in $\vcat{\sset_{(\infty,1)}}$, so in particular in the localization $\cL_{\Sigma^2\Lambda}\vcat{\sset_{(\infty,1)}}$
as desired.

Now if $\cQ$ is more generally any (cofibrant) simplicial set, we consider a fibrant replacement $\cQ^{\mathrm{fib}}$ in  $\cL_{\Sigma^2\Lambda}\vcat{\sset_{(\infty,1)}}$ and the following naturality diagram.
\begin{tz}
\node[](1) {$\cQ$}; 
\node[right of=1,xshift=.7cm](2) {${\bf N}_*c_*\cQ$}; 
\node[below of=1](3) {$\cQ^\mathrm{fib}$}; 
\node[below of=2](4) {${\bf N}_*c_*(\cQ^\mathrm{fib})$};

\draw[->] (1) to (2);
\draw[->] (1) to (3);
\draw[->] (2) to (4);
\draw[->] (3) to (4); 
\end{tz}
% \[\begin{tikzcd}
%     X\arrow[d]\arrow[r]&{\bf N}c X\arrow[d]\\
%   X^{\mathrm{fib}}\arrow[r]&{\bf N}((c(X^{\mathrm{fib}}))).\\ 
% \end{tikzcd}\]
% \[\begin{tikzcd}
%     X\arrow[d]\arrow[r]&N((c X)^{\mathrm{fib}})\arrow[d]\\
%   X^{\mathrm{fib}}\arrow[r]&N((c(X^{\mathrm{fib}}))^{\mathrm{fib}})\\ 
% \end{tikzcd}\]
Here, the left vertical map is a weak equivalence in $\cL_{\Sigma^2\Lambda}\vcat{\sset_{(\infty,1)}}$ by construction, the right vertical map is a weak equivalence because both $\bf N_*$ and $c_*$ are homotopical by \cref{c*homotopical,nervestarrightQuillen}, and the bottom horizontal arrow is a weak equivalence by the case that we already treated. It follows by $2$-out-of-$3$ that the top horizontal map, which is the (derived) unit of $\cQ$, is a weak equivalence, as desired.
\end{proof}

\begin{rmk}
\label{equivalencewithGH2}
By \cite[Rmk 0.0.4]{LurieGoodwillie}, we know that the underlying $\infty$-category of the model structure $\vcat{\sset_{(\infty,1)}}$ models the established $\infty$-category $\mathscr Cat_{(\infty,2)}$ of $(\infty,2)$-ca\-te\-gories, so there exists an equivalence of $\infty$-categories
\begin{equation}
    \label{correctinfinity2categories}
  [\vcat{\sset_{(\infty,1)}}]_\infty\simeq \mathscr Cat_{(\infty,2)}.
\end{equation}
Any such equivalence can be used to construct a specific
equivalence of $\infty$-categories
\begin{align}\label{correctinfinity2categories2}
    [\twocat]_\infty&\simeq[\cL_{\Sigma^2\Lambda}\vcat{\sset_{(\infty,1)}}]_\infty&\text{\cref{QE2}} \nonumber\\
    &\simeq\mathscr L_{\Sigma^2\Lambda}[\vcat{\sset_{(\infty,1)})}]_\infty&\text{\cref{localization compatible}}\\
    &\simeq\mathscr L_{\Sigma^2\Lambda}\mathscr Cat_{(\infty,2)}&\text{\eqref{correctinfinity2categories},\ \cref{suspensioncomparison}}\nonumber\\
     &\simeq\mathscr Cat_2&\text{\cite[Lem.~6.1.6(1)]{GH}}\nonumber
\end{align}
     between the established $\infty$-category of $2$-categories $\mathscr Cat_2$ and the underlying $\infty$-category of the model structure $\twocat$ on $2$-categories. Via the chosen identifications \eqref{correctinfinity2categories} and \eqref{correctinfinity2categories2},
we see that the functor $[{\mathbf N}_*]_\infty\colon[\twocat]_\infty\to[\vcat{{\sset_{(\infty,1)}}}]_\infty$ and the canonical inclusion $\mathscr Cat_2\hookrightarrow\mathscr Cat_{(\infty,2)}$ from \cite[Lem.~6.1.6(v)]{GH} -- for $n=2$ -- are equivalent. Indeed, this is witnessed by the following diagram of $\infty$-categories
\begin{tz}
\node[](1) {$[\twocat]_\infty$}; 
\node[below of=1](1') {$[\vcat{\sset_{(\infty,1)}}]_\infty$}; 
\node[right of=1,xshift=1.21cm](2) {$[\cL_{\Sigma^2\Lambda}\vcat{\sset_{(\infty,1)}}]_\infty$}; 
\node[below of=2](2') {$[\vcat{\sset_{(\infty,1)}}]_\infty$}; 
\node[right of=2,xshift=2.16cm](3) {$\mathscr L_{\Sigma^2\Lambda}[\vcat{\sset_{(\infty,1)}}]_\infty$};
\node[below of=3](3') {$[\vcat{\sset_{(\infty,1)}}]_\infty$}; 
\node[right of=3,xshift=1.76cm](4) {$\mathscr L_{\Sigma^2\Lambda} \mathscr Cat_{(\infty,1)}$};
\node[below of=4](4') {$\mathscr Cat_{(\infty,2)}$}; 
\node[right of=4,xshift=.61cm,yshift=1pt](5) {$\mathscr Cat_2$}; 
\node[below of=5,yshift=-1pt](5') {$\mathscr Cat_{(\infty,2)}$}; 

\draw[->] ($(1.east)+(0,1pt)$) to node[above,la]{$\simeq$} ($(2.west)+(0,1pt)$);
\draw[->] ($(2.east)+(0,1pt)$) to node[above,la]{$\simeq$} ($(3.west)+(0,1pt)$);
\draw[->] ($(3.east)+(0,1pt)$) to node[above,la]{$\simeq$} ($(4.west)+(0,1pt)$);
\draw[->] ($(4.east)+(0,1pt)$) to node[above,la]{$\simeq$} (5);

\draw[d] (1') to (2');
\draw[d] (2') to (3');
\draw[->] (3') to node[below,la]{$\simeq$} (4');
\draw[d] (4') to (5');

\draw[right hook->] (1) to node[left,la]{$[\bfN_*]_\infty$} (1');
\draw[right hook->] (2) to (2');
\draw[right hook->] (3) to (3');
\draw[right hook->] (4) to (4');
\draw[right hook->] (5) to (5');
\end{tz}
% \[\begin{tikzcd}
% \twocat_\infty\arrow[r,"\simeq"]\arrow[d,"(\mathbf{N}_*)_\infty" swap]&(\cL_{\Sigma^2\Lambda}(\vcat{\sset_{(\infty,1)}}))_\infty\arrow[d,""]\arrow[r,"\simeq"]&\mathscr L_{\Sigma^2\Lambda}(\vcat{(\sset_{(\infty,1)})_\infty})\arrow[d]\arrow[r,"\simeq"]&\mathscr L_{\Sigma^2\Lambda}(\mathscr Cat_{(\infty,2)})\arrow[d,""]\arrow[r,"\simeq"]&\mathscr Cat_2\arrow[d,""]\\
% (\sset_{(\infty,1)})_\infty\arrow[r,equal]&(\sset_{(\infty,2)})_\infty\arrow[r,equal]&\vcat{(\sset_{(\infty,1)}})_\infty\arrow[r,"\simeq" swap]&\mathscr Cat_{(\infty,2)}\arrow[r,equal]&\mathscr Cat_{(\infty,2)}
% \end{tikzcd}\]
% \[\begin{tikzcd}
% [\twocat]_\infty\arrow[r,"(\mathrm{\mathbf N_*)_\infty}"]\arrow[d,"\simeq"]&\sset_{(\infty,1)}\arrow[d,equal]\\
% (\cL(\vcat{\sset_{(\infty,1)})}))_\infty\arrow[r]\arrow[d,"\simeq"]&\sset_{(\infty,1)}\arrow[d,equal]\\
% \mathscr L((\vcat{\sset_{(\infty,1)})}))_\infty\arrow[r]\arrow[d,"\simeq"]&\sset_{(\infty,1)}\arrow[d,"\simeq"]\\
% \mathscr L\mathscr Cat_\infty\arrow[r]\arrow[d,"\simeq"]&\mathscr Cat_\infty\arrow[d,equal]\\
% \mathscr Cat_{(2,2)}\arrow[r]\arrow[d,"\simeq"]&\mathscr Cat_{(\infty,2)}\arrow[d,equal]\\
% \mathscr Cat_{2}\arrow[r]&\mathscr Cat_{(\infty,2)}\\
% \end{tikzcd}\]
which commutes up to equivalence, using \cref{localization compatible,suspensioncomparison} and \cite[Def.~6.1.7(v)]{GH}.
\end{rmk}

%%%
%%%

\appendix

\section{The nerve comparison lemma}

\label{NerveComparisonAppendix}

To assert the commutativity at the level of $\infty$-categories of each of the regions in the diagram from \cref{{MainTheorem}}, we will make use of the following lemma.

\begin{lem}[Nerve comparison lemma]
\label{NerveRecognition}
Let $\cM$ and $\cM'$ be two model categories. Suppose we are given the following:
\begin{itemize}[leftmargin=*]
    \item a left Quillen functor, resp.~right Quillen functor, $F\colon\cM\to\cM'$;
    \item a homotopical functor $H\colon\twocat\to\cM$ that takes values in the subcategory of cofibrant, resp.~fibrant, objects in $\cM$;
    \item a homotopical functor $H'\colon\twocat\to\cM'$; and
    \item a natural weak equivalence $FH\stackrel{\simeq}{\Longrightarrow}H'$.
\end{itemize}
Then, the diagram of categories on the left
\begin{tz}
\node[](1) {$\twocat$}; 
\node[below of=1,xshift=-1.7cm,yshift=.3cm](3) {$\cM$}; 
\node[below of=1,xshift=1.7cm,yshift=.3cm](4) {$\cM'$}; 

\draw[->] (1) to node[above,la,xshift=-15pt]{$H$} (3);
\draw[->] (1) to node[above,la,xshift=15pt]{$H'$} (4); 
\draw[->] (3) to node[below,la]{$F$} (4); 

\node[right of=1,xshift=5cm](1) {$[\twocat]_\infty$}; 
\node[below of=1,xshift=-1.7cm,yshift=.3cm](3) {$[\cM]_\infty$}; 
\node at ($(4)!0.5!(3)+(0,.6cm)$) {$\rightsquigarrow$}; 
\node[below of=1,xshift=1.7cm,yshift=.3cm](4) {$[\cM']_\infty$}; 

\draw[->] (1) to node[above,la,xshift=-15pt]{$[H]_\infty$} (3);
\draw[->] (1) to node[above,la,xshift=15pt]{$[H']_\infty$} (4); 
\draw[->] (3) to node[below,la]{$[F]_\infty$} (4);
\end{tz}
% \[
% \begin{tikzcd}
% &\twocat\ar[dr, "H'"]\ar[dl, "H"swap] & \\
% \cM\ar[rr, "F"swap] & & \cM'
% \end{tikzcd}
% \quad\rightsquigarrow\quad
% \begin{tikzcd}
% &\twocat_{\infty}\ar[dr, "H_\infty'"]\ar[dl, "H_\infty"swap] & \\
% \cM_{\infty}\ar[rr, "F_\infty"swap] & & [\cM']_\infty
% \end{tikzcd}
% \]
induces a diagram of $\infty$-categories that commutes up to equivalence.\footnote{Meaning that the two functors are equivalent in the $\infty$-category of functors from $[\twocat]_{\infty}\to[\cM']_{\infty}$.}
\end{lem}

\begin{rmk}
The second (resp.~third) condition of \cref{NerveRecognition} is automatically satisfied when $H\colon \twocat\to \cM$ (resp.~$H'\colon \twocat\to \cM'$) is right Quillen.
\end{rmk}

We choose to work with the following model of $[\cM]_\infty$ for a model category $\cM$, regarded as a relative category $(\cM,\cW)$ when equipped with its class of weak equivalences $\cW$.

Following e.g.~\cite[Const.\ 15.1]{BarwickSchommerPries}, given a relative category $(\cC,\cW)$, the underlying $\infty$-category is \[ [\cC]_\infty\coloneqq\mathfrak{N}((L_H(\cC, \cW))^{\mathrm{fib}}).\]
Here, the functor $\mathfrak N\colon\vcat{\sset_{(\infty,0)}}\to\sset_{(\infty,1)}$ denotes the homotopy coherent nerve functor defined by \cite{CordierHoCoherent} and which is a right Quillen functor by \cite[Thm 2.2.5.1]{htt}, while $(-)^{\mathrm{fib}}\colon\vcat{\sset_{(\infty,0)}}\to\vcat{\sset_{(\infty,0)}}$ denotes any functorial fibrant replacement in the Bergner model structure $\vcat{\sset_{(\infty,0)}}$ from \cite[Thm 3.2.4, Ex.\ 3.2.23]{htt}; for instance, one could take $(\mathrm{Ex}^{\infty})_*\colon\vcat{\sset_{(\infty,0)}}\to\vcat{\sset_{(\infty,0)}}$. 

The following fact is essentially discussed in \cite[\S A.3.1]{MazelGeeQuillenAdj}, following \cite[Prop.\ 3.3, 3.5]{DwyerKanCalculating}.

\begin{prop}
\label{equivalentHoFunctors}
Let $G,G'\colon (\cC, \cW) \to (\cC', \cW')$ be homotopical functors of relative categories, and let $\alpha\colon G\stackrel{\simeq}{\Longrightarrow}G'$ a natural weak equivalence. Then $G$ and $G'$ induce equivalent functors of quasi-categories
\[ [G]_\infty\simeq [G']_\infty\colon [\cC]_\infty=\mathfrak{N}((L_H(\cC, \cW))^\mathrm{fib}) \to \mathfrak{N}((L_H(\cC', \cW'))^\mathrm{fib})=[\cC']_\infty.\]
\end{prop}

We can now prove the lemma.

\begin{proof}[Proof of \cref{NerveRecognition}] 
The lemma follows from \cref{equivalentHoFunctors} by taking $G=FH$ and $G'=H'$. Indeed, we have equivalences of functors
\[[H']_\infty\simeq [F H]_\infty\simeq [F]_\infty\circ [H]_\infty,\]
which concludes the proof.
\end{proof}

\section{Complements on the Rezk nerve of categories} \label{AppendixRezk}

We collect in this appendix a series of elementary properties of the Rezk nerve that we did not find in the literature. We denote by $\widetilde{[k]}$ the contractible groupoid with $k+1$ objects. 

\begin{const}[{\cite[\S3.5]{rezk}}]
Let $\cC$ be a category. The \emph{Rezk nerve} $\NRezk\cC$ is the simplicial space given for any $j,k\ge0$ by
\[ \NRezk_{j,k}\cC\coloneqq \cat([j]\times \widetilde{[k]}, \cC).\]
The assignment extends to a functor $\NRezk\colon \cat\to\spsh{\Delta}$.
\end{const}

Recall from \cite[Rmk 5.6]{rezk} that the Rezk nerve has a left adjoint $c^R\colon\spsh{\Delta}\to\cat$.

\begin{lem}
\label{RezkInducesMonoidal}
The left adjoint $c^R\colon\spsh{\Delta}\to\cat$ preserves finite products.
\end{lem}

\begin{proof}
Since both $\cat$ and $\spsh{\Delta}$ are cartesian closed, products commute with colimits, hence it suffices to prove that for any $j,k,j',k'\geq 0$ we have an isomorphism of bisimplicial sets
\[c^R(\Delta[j,k]\times \Delta[j',k'])\cong c^R(\Delta[j,k])\times c^R(\Delta[j',k']).\]
We will prove that both sides are isomorphic to $ [j]\times [j'] \times \widetilde{([k]\times [k'])}$. For the right-hand side, we have
\begin{align*}
 c^R(\Delta[j,k])\times c^R(\Delta[j',k'])&\cong [j]\times \widetilde{[k]}\times [j']\times \widetilde{[k']}\\
       &\cong [j]\times [j'] \times \widetilde{[k]}\times \widetilde{[k']}\\
      &\cong [j]\times [j'] \times \widetilde{[k]\times [k']}.
\end{align*}
For the left-hand side, we need the following observations.
\begin{enumerate}[leftmargin=*]
    \item For all $j,k\geq 0$ there is an isomorphism of bisimplicial sets
        \[\Delta[j,k] \cong \Delta[j,0]\times \Delta[0,k].\]
    \item The functor $\widetilde{(-)}\colon \cat\to \cG pd$ is left adjoint to the inclusion functor $\cG pd\hookrightarrow \cat$; in particular, the functor $\widetilde{(-)}$ preserves colimits. 
    \item The left adjoint $c\colon\sset\to\cat$ of the ordinary nerve functor preserves colimits, and it also preserves finite products by \cite[Prop.\ B.0.15]{JoyalVolumeII},
    there attributed to Gabriel--Zisman. Then, for any $j,j'\geq 0$ we obtain an isomorphism of categories
    \begin{align*}
       [j]\times [j']  & \cong c\Delta[j]\times c\Delta[j'] \cong c(\Delta[j]\times \Delta[j'])\\
         & \cong c(\underset{\Delta\downarrow\Delta[j]\times \Delta[j']}{\colim} \Delta[a]) \cong \underset{\Delta\downarrow\Delta[j]\times \Delta[j']}{\colim} c\Delta[a]\\
     & \cong  \underset{\Delta\downarrow\Delta[j]\times \Delta[j']}{\colim} [a].
    \end{align*}
\end{enumerate}
We then have the following isomorphisms of categories
\begin{align*}
  c^R(\Delta[j,k]\times \Delta[j',k'])   & \cong c^R(\Delta[j,0]\times\Delta[j',0]\times \Delta[0,k]\times\Delta[0,k'])&\text{Obs.~(1)}\\
     & \cong c^R(\underset{\Delta\downarrow\Delta[j]\times \Delta[j']}{\colim} \Delta[a,0] \times \underset{\Delta\downarrow\Delta[k]\times \Delta[k']}{\colim} \Delta[0,b])&\\
     & \cong c^R( \underset{\Delta\downarrow\Delta[j]\times \Delta[j']}{\colim}\,\underset{\Delta\downarrow\Delta[k]\times \Delta[k']}{\colim}  (\Delta[a,0] \times  \Delta[0,b]))&\\
     & \cong \underset{\Delta\downarrow\Delta[j]\times \Delta[j']}{\colim}\,\underset{\Delta\downarrow\Delta[k]\times \Delta[k']}{\colim} c^R( \Delta[a,0] \times  \Delta[0,b])&\text{$c^R$ left adjoint}\\
     &\cong \underset{\Delta\downarrow\Delta[j]\times \Delta[j']}{\colim}\,\underset{\Delta\downarrow\Delta[k]\times \Delta[k']}{\colim} c^R( \Delta[a,b])&\text{Obs.~(1)}\\
     &\cong \underset{\Delta\downarrow\Delta[j]\times \Delta[j']}{\colim}\,\underset{\Delta\downarrow\Delta[k]\times \Delta[k']}{\colim} [a]\times \widetilde{[b]}&\\
      &\cong \underset{\Delta\downarrow\Delta[j]\times \Delta[j']}{\colim} [a]\times \underset{\Delta\downarrow\Delta[k]\times \Delta[k']}{\colim} \widetilde{[b]}&\\
      &\cong \underset{\Delta\downarrow\Delta[j]\times \Delta[j']}{\colim} [a]\times\widetilde{(\underset{\Delta\downarrow\Delta[k]\times \Delta[k']}{\colim} [b])}&\text{Obs.~(2)}\\
      &\cong [j]\times [j'] \times \widetilde{[k]\times [k']},&\text{Obs.~(3)}
\end{align*}
as desired.
\end{proof}

\begin{prop}
\label{RezkRightQuillen}
The Rezk nerve $\NRezk\colon\cat\to\spsh{\Delta}_{(\infty,1)}$ is a right Quillen embedding, and in particular a right Quillen and homotopical functor.
\end{prop}

\begin{proof}
We argue that the functor $\NRezk\colon\cat\to\spsh{\Delta}_{(\infty,1)}$ can be understood as the composite of the ordinary nerve ${\bf N}\colon\cat\to\sset_{(\infty,1)}$ of categories into simplicial sets, which is easily seen to be a right Quillen embedding and the functor $t^!\colon \sset_{(\infty,1)}\to\spsh{\Delta}_{(\infty,1)}$ from \cite[\S4]{JT}, which is shown to be a right Quillen equivalence. It will then follow that $\NRezk$ is a right Quillen embedding.

In order to prove the claim, we observe that for any category $\cC$ and $j,k\ge0$ there is a natural bijection
\begin{align*}
     (t^!{\bf N}\cC)_{j,k}&\cong \spsh{\Delta}(\Delta[j,k],t^!{\bf N}\cC)  \cong \sset(t_!\Delta[j,k],{\bf N}\cC)\\
     &\cong \sset(\Delta[j]\times {\bf N}\widetilde{[k]},{\bf N}\cC)\cong \sset({\bf N}([j]\times\widetilde{[k]}),{\bf N}\cC)\\
     &\cong \cat([j]\times\widetilde{[k]},\cC)\cong \NRezk_{j,k}\cC,
     \end{align*}
     as desired.
\end{proof}

%%%
%%%

\section{Complements on the bisimplicial nerve of \texorpdfstring{$2$}{2}-categories}

The homotopically correct nerve of $2$-categories into $2$-quasi-categories is based on the notion of \emph{normal pseudofunctor}, also referred to as \emph{normalized} or \emph{strictly unital pseudofunctor} or \emph{homomorphism}, or \emph{weak functor}.
Roughly speaking, a normal pseudofunctor is a map between $2$-categories that preserves identities strictly and preserves composition up to coherent
isomorphism. We now recall the main aspects of the definitions, referring the reader to other sources, see e.g. B\'enabou \cite[Rmk~4.2]{Benabou}, Street \cite[Ex.~9.7]{street} or Johnson--Yau \cite[Def.~4.1]{JY}, for a more detailed treatment.

Given a $2$-category $\cA$, we denote by $\Ob\cA$, $\Mor\cA$, and $2\Mor\cA$ the sets of objects, $1$-morphisms, and $2$-morphisms in $\cA$, respectively. We denote by $s$, $t$, $i$, and $c$ the source, target, identity, and composition maps for $1$-morphisms, and by $s$, $t$, $i$, $c_h$, and $c_v$ the source, target, identity, horizontal composition, and vertical composition maps for $2$-morphisms. 

We denote by $\Triangle\cA\coloneqq \twocat(\Osim{2},\cA)$, the set of $2$-isomorphisms in $\cA$ of the form 
\begin{tz}
\node[](1) {$x$};
\node[above right of=1,xshift=.3cm](2) {$y$};
\node[below right of=2,xshift=.3cm](3) {$z$};
\punctuation{3}{,}[1pt];
\draw[->] (1) to node[above,la,xshift=-2pt]{$f$} (2);
\draw[->] (2) to node[above,la,xshift=2pt]{$g$} (3);
\draw[->] (1) to node[below,la]{$h$} (3);
\coordinate(a) at ($(1)!0.5!(3)$);
\cell[la,right][n][0.37]{a}{2}{$\cong$};
\end{tz}
which comes with three maps $d_0,d_1,d_2\colon \Triangle\cA\to \Mor\cA$ picking each of the boundary of the $2$-isomorphisms, and two maps $s_0,s_1\colon \Mor\cA\to \Triangle\cA$ sending a $1$-morphism to its identity $2$-morphism in the two usual ways.

Finally, we denote by $2\Iso\cA$, the set of $2$-isomorphisms in $\cA$. Note that there is a map $e\colon \Triangle\cA\to 2\Iso\cA$, which extracts the $2$-isomorphism component, e.g. it sends the above picture to the corresponding $2$-isomorphism $h\cong gf$.

\begin{defn}
\label{normalpseudofunctor}
A \emph{normal pseudofunctor} 
$F\colon\cA\to\cB$ between two $2$-categories $\cA$ and $\cB$ consists of the following data
\begin{enumerate}[leftmargin=*]
\addtocounter{enumi}{-1}
    \item an assignment on objects, namely a function
    $F_0\colon\Ob\cA\to\Ob\cB$;
\item an assignment on $1$-morphisms, namely a function
    $F_1\colon\Mor\cA\to\Mor\cB$;
\item an assignment on $2$-morphisms, namely a function
    $F_2\colon2\Mor\cA\to2\Mor\cB$;
\item a \emph{compositor}
of $F$, namely a function
    $\widetilde F\colon \Mor\cA\times_{\Ob \cA}\Mor\cA\to \Triangle\cB$;
\end{enumerate}
with the requirement that the following axioms be satisfied.
\begin{enumerate}[leftmargin=*, label=(\alph*)]
    \item The assignments of $F$ on objects, $1$- and $2$-morphisms commute with source, target, and identities: 
    \begin{tz}
    \node[](1) {$\Ob\cA$}; 
    \node[below of=1](2) {$\Ob\cB$}; 
    \node[right of=1,xshift=1.5cm](1') {$\Mor\cA$}; 
    \node[below of=1'](2') {$\Mor\cB$}; 
    \node[right of=1',xshift=1.5cm](1'') {$2\Mor\cA$}; 
    \node[below of=1''](2'') {$2\Mor\cB$}; 
    
    \draw[->] ($(1'.west)+(0,5pt)$) to node[above,la]{$s$} ($(1.east)+(0,5pt)$);
    \draw[->] ($(1'.west)-(0,5pt)$) to node[below,la]{$t$} ($(1.east)-(0,5pt)$);
    \draw[->] (1) to node[over,la]{$i$} (1');
    \draw[->] ($(2'.west)+(0,5pt)$) to node[above,la]{$s$} ($(2.east)+(0,5pt)$);
    \draw[->] ($(2'.west)-(0,5pt)$) to node[below,la]{$t$} ($(2.east)-(0,5pt)$);
    \draw[->] (2) to node[over,la]{$i$} (2');
    
    \draw[->] ($(1''.west)+(0,5pt)$) to node[above,la]{$s$} ($(1'.east)+(0,5pt)$);
    \draw[->] ($(1''.west)-(0,5pt)$) to node[below,la]{$t$} ($(1'.east)-(0,5pt)$);
    \draw[->] (1') to node[over,la]{$i$} (1'');
    \draw[->] ($(2''.west)+(0,5pt)$) to node[above,la]{$s$} ($(2'.east)+(0,5pt)$);
    \draw[->] ($(2''.west)-(0,5pt)$) to node[below,la]{$t$} ($(2'.east)-(0,5pt)$);
    \draw[->] (2') to node[over,la]{$i$} (2'');
    
    \draw[->] (1) to node[left,la]{$F_0$} (2);
    \draw[->] (1') to node[left,la]{$F_1$} (2');
    \draw[->] (1'') to node[right,la]{$F_2$} (2'');
    \end{tz}
    This gives that the images under $F$ of a $1$-morphism $f\colon x\to y$ and a $2$-morphism $\alpha\colon f\Rightarrow g$ are of the form ${Ff\colon Fx\to Fy}$ and $F\alpha\colon Ff\Rightarrow Fg$, respectively, and that $F(\id_x)=\id_{Fx}$ and $F(\id_f)=\id_{Ff}$ for any object $x$ and any $1$-morphism $f$.
    \item The boundaries of $\widetilde{F}$ is determined by the following commutative diagram:
    \begin{tz}
    \node[](1) {$\Mor\cA$}; 
    \node[below of=1](2) {$\Mor\cB$}; 
    \node[right of=1,xshift=2cm,yshift=-4pt](1') {$\Mor\cA\ttimes{\Ob\cA}\Mor\cA$}; 
    \node[below of=1',yshift=4pt](2') {$\Triangle\cB$};
    \node[right of=1',xshift=3cm](1'') {$\Mor\cA\ttimes{\Ob\cA}\Mor\cA$};
    \node[below of=1''](2'') {$\Mor\cB\ttimes{\Ob\cB} \Mor\cB$}; 
    
    \draw[->] ($(1'.west)+(0,4pt)$) to node[above,la]{$c$} (1);
    \draw[->] (2') to node[below,la]{$d_1$} (2);
    
    \draw[d] ($(1'.east)+(0,4pt)$) to ($(1''.west)+(0,4pt)$);
    \draw[->] ($(1''.south)+(0,2pt)$) to node[right,la]{$F_1\times F_1$} (2'');
    \draw[->] (2') to node[below,la]{$(d_2,d_0)$} ($(2''.west)+(0,4pt)$);
    
    \draw[->] (1) to node[left,la]{$F_1$} (2);
    \draw[->] ($(1'.south)+(0,2pt)$) to node[left,la]{$\widetilde F$} (2');
    \end{tz}
    When evaluated at an element $(f\colon x\to y, g\colon y\to z)$ this gives a $2$-isomorphism $\widetilde F_{f,g}$ of the following form.
    \begin{tz}
\node[](1) {$Fx$};
\node[above right of=1,xshift=.3cm](2) {$Fy$};
\node[below right of=2,xshift=.3cm](3) {$Fz$};
\draw[->] (1) to node[above,la,xshift=-2pt]{$Ff$} (2);
\draw[->] (2) to node[above,la,xshift=2pt]{$Fg$} (3);
\draw[->] (1) to node[below,la]{$F(gf)$} (3);
\coordinate(a) at ($(1)!0.5!(3)$);
\cell[la,right][n][0.37]{a}{2}{$\widetilde F_{f,g}$};
\end{tz}
    \item The compositor $\widetilde{F}$ is compatible with identities in the sense that the following diagram commutes:
    \begin{tz}
    \node[](1) {$\Mor\cA$}; 
    \node[below of=1](2) {$\Mor\cB$}; 
    \node[right of=1,xshift=2cm,yshift=-4pt](1') {$\Mor\cA\ttimes{\Ob\cA}\Mor\cA$}; 
    \node[below of=1',yshift=4pt](2') {$\Triangle\cB$}; 
    \node[right of=1',xshift=2cm,yshift=4pt](1'') {$\Mor\cA$}; 
    \node[below of=1''](2'') {$\Mor\cB$}; 
    
    \draw[->] (1) to node[above,la]{$(\id,it)$} ($(1'.west)+(0,4pt)$);
    \draw[->] (2) to node[below,la]{$s_1$} (2');
    
    \draw[->] (1'') to node[above,la]{$(is,\id)$} ($(1'.east)+(0,4pt)$);
    \draw[->] (2'') to node[below,la]{$s_0$} (2');
    
    \draw[->] (1) to node[left,la]{$F_1$} (2);
    \draw[->] ($(1'.south)+(0,2pt)$) to node[left,la]{$\widetilde{F}$} (2');
    \draw[->] (1'') to node[right,la]{$F_1$} (2'');
    \end{tz}
    When evaluated at an element $f\colon x\to y$, this gives that the $2$-isomorphisms $\widetilde F_{f,\id_y}$ and $\widetilde F_{\id_x,f}$ are both the identity $2$-morphism at $f$.
    \item The assignment of $F$ on $2$-morphisms commutes with vertical composition of $2$-morphisms: 
    \begin{tz}
    \node[](1') {$2\Mor\cA\ttimes{\Mor\cA}2\Mor\cA$}; 
    \node[right of=1',xshift=4cm](2') {$2\Mor\cB\ttimes{\Mor\cB}2\Mor\cB$}; 
    \node[below of=1',yshift=4pt](1'') {$2\Mor\cA$}; 
    \node[below of=2',yshift=4pt](2'') {$2\Mor\cB$};
    
    \draw[->] ($(1'.south)+(0,2pt)$) to node[left,la]{$c_v$} (1'');
    \draw[->] ($(2'.south)+(0,2pt)$) to node[right,la]{$c_v$} (2'');
    
    \draw[->] ($(1'.east)+(0,4pt)$) to node[above,la]{$F_2\times F_2$} ($(2'.west)+(0,4pt)$);
    \draw[->] (1'') to node[below,la]{$F_2$} (2'');
    \end{tz}
    When evaluated at an element $(\alpha\colon f\Rightarrow g\colon x\to y, \beta\colon g\Rightarrow h\colon x\to y)$, this gives that $F(\beta \alpha)=(F\beta)(F\alpha)$. 
    \item The compositor $\widetilde{F}$ is $2$-natural in the sense that the following diagram commutes:
    \end{enumerate}
    \begin{tz}
    \node[](1) {$2\Mor\cA\ttimes{\Ob\cA}2\Mor\cA$}; 
    \node[right of=1,xshift=6cm,yshift=-2.5pt](2) {$\Triangle\cB\ttimes{\Mor\cB\ttimes{\Ob\cB}\Mor\cB}(2\Mor\cB\ttimes{\Ob\cB}2\Mor\cB)$}; 
    \node[below of=1](1') {$2\Mor\cA\ttimes{\Ob\cA}2\Mor\cA$}; 
    \node[below of=2,yshift=6.5pt](2') {$2\Mor\cB$};
    
    \node[below of=1'](1'') {$2\Mor\cA\ttimes{\Mor\cA} (\Mor\cA\ttimes{\Ob\cA}\Mor\cA)$}; 
    \node[below of=2',yshift=-4pt](2'') {$2\Mor\cB\ttimes{\Mor\cB} \Triangle\cB$}; 
    
    \draw[d] (1') to ($(1.south)+(0,2pt)$);
    \draw[->] ($(2.south)+(0,8pt)$) to node[right,la]{$c_v (\id\times c_h)$} (2');
    
    \draw[->] ($(1'.south)+(0,2pt)$) to node[left,la]{$(c_h,t\times t)$} (1'');
    \draw[->] (2'') to node[right,la]{$c_v$} (2');
    
    \draw[->] ($(1.east)+(0,4pt)$) to node[above,la]{$(\widetilde{F}(s\times s),F_2\times F_2)$} ($(2.west)+(0,6.5pt)$);
    \draw[->] ($(1''.east)+(0,4pt)$) to node[below,la]{$F_2\times\widetilde F$} ($(2''.west)+(0,4pt)$);
    \end{tz}
    \hfill\begin{minipage}{0.948\textwidth}
    When evaluated at an element $(\alpha\colon f\Rightarrow f'\colon x\to y,\beta\colon g\Rightarrow g'\colon y\to z)$, this gives the following pasting equality. 
    \begin{tzbigbig}
\node[](1) {$Fx$};
\node[above of=1](2) {$Fy$}; 
\node[right of=2](3) {$Fz$}; 

\draw[->,bend left] (1) to node(a)[left,la]{$Ff'$} (2);
\draw[->,bend right] (1) to node(b)[over,la]{$Ff$} (2);
\draw[->,bend right=27] (2) to node(d)[over,la]{$Fg$} (3);
\draw[->,bend left=27] (2) to node(c)[above,la]{$Fg'$} (3);
\draw[->,bend right=20] (1) to node(e)[right,la,pos=0.4,yshift=-2pt]{$F(gf)$} (3);
\coordinate(f) at ($(e)+(0,-8pt)$);
\cell[above,la][n][0.42]{b}{a}{$F\alpha$};
\cell[right,la,xshift=-2pt][n][0.42]{d}{c}{$F\beta$};
\cell[above,la,xshift=11pt,yshift=-4pt][n][0.43]{f}{2}{$\widetilde{F}_{f,g}$};

\node[right of=1,xshift=2cm](1) {$Fx$};
\node at ($(1)!0.5!(3)-(.4cm,0)$) {$=$};
\node[above of=1](2) {$Fy$}; 
\node[right of=2](3) {$Fz$}; 

\draw[->,bend left=15] (1) to node[left,la]{$Ff'$} (2);
\draw[->,bend left=15] (2) to node[above,la]{$Fg'$} (3);
\draw[->,bend right=25] (1) to node(a)[right,la,pos=0.4,yshift=-2pt]{$F(gf)$} (3);
\draw[->,bend left=25] (1) to node(b)[over,la]{$F(g'f')$} (3);
\coordinate(c) at ($(a)+(0,-8pt)$);
\cell[left,la,xshift=8pt,yshift=-10pt][n][0.65]{c}{b}{$F(\beta\alpha)$};
\cell[above,la,xshift=14pt,yshift=-4pt][n][0.5]{b}{2}{$\widetilde{F}_{f',g'}$};
\end{tzbigbig}
\end{minipage}
    \begin{enumerate}[leftmargin=*, label=(\alph*)]\addtocounter{enumi}{5}
        \item The compositor $\widetilde{F}$ is compatible with composition of $1$-morphisms in the sense that the following diagram commutes:
    \begin{tz}
    \node[](1) {$(\Mor\cA\ttimes{\Ob\cA}\Mor\cA)\ttimes{\Mor\cA}(\Mor\cA\ttimes{\Ob\cA}\Mor\cA)$}; 
    \node[right of=1,xshift=5.5cm](2) {$\Triangle\cB\ttimes{\Mor\cB}\Triangle\cB$};
    \node[below of=1](1') {$\Mor\cA\ttimes{\Ob\cA}\Mor\cA\ttimes{\Ob\cA}\Mor\cA$}; 
    \node[below of=2,yshift=4pt](2') {$2\Mor\cB$}; 
    \node[below of=1'](1'') {$(\Mor\cA\ttimes{\Ob\cA}\Mor\cA)\ttimes{\Mor\cA}(\Mor\cA\ttimes{\Ob\cA}\Mor\cA)$};
    \node[below of=2',yshift=-4pt](2'') {$\Triangle\cB\ttimes{\Mor\cB}\Triangle\cB$}; 
    
    \draw[->] (1') to node[left,la]{$(\id\times c,!\times\id\times \id)$} node[right,la]{$\cong$} ($(1.south)+(0,2pt)$);
    \draw[->] ($(2.south)+(0,2pt)$) to node[right,la]{$\varphi$} (2');
    
    \draw[->] ($(1'.south)+(0,6pt)$) to node[left,la]{$(\id\times\id\times !, c\times \id)$} node[right,la]{$\cong$} (1'');
    \draw[->] (2'') to node[right,la]{$\psi$} (2');
    
    \draw[->] ($(1.east)+(0,4pt)$) to node[above,la]{$\widetilde F\times\widetilde F$} ($(2.west)+(0,4pt)$);
    \draw[->] ($(1''.east)+(0,4pt)$) to node[below,la]{$\widetilde F\times\widetilde F$} ($(2''.west)+(0,4pt)$);
    \end{tz}
    where $\varphi$ and $\psi$ compute the total composite of the pasting diagrams. When evaluated at an element $(f\colon x\to y, g\colon y\to z, h\colon z\to w)$, this gives the following pasting equality.
\begin{tzbig}
\node[](1) {$Fx$};
\node[above of=1](2) {$Fy$};
\node[right of=2](3) {$Fz$};
\node[below of=3](4) {$Fw$};
\draw[->] (1) to node[left,la]{$Ff$} (2);
\draw[->] (2) to node[above,la]{$Fg$} (3);
\draw[->] (3) to node[right,la]{$Fh$} (4);
\draw[->] (1) to node[below,la]{$F(hgf)$} (4);
\draw[->] (1) to node[over,la]{$F(gf)$} (3);
\coordinate(a) at ($(1)!0.5!(3)$);
\cell[la,right,yshift=7pt][n][0.45]{a}{2}{$\widetilde F_{f,g}$};
\coordinate(a) at ($(1)!0.78!(4)$);
\coordinate(b) at ($(a)+(0,1.5cm)$);
\cell[la,left][n][0.33]{a}{b}{$\widetilde F_{gf,h}$};

\node[right of=4](1) {$Fx$};
\node at ($(3)!0.5!(1)$) {$=$};
\node[above of=1](2) {$Fy$};
\node[right of=2](3) {$Fz$};
\node[below of=3](4) {$Fw$};

\draw[->] (1) to node[left,la]{$Ff$} (2);
\draw[->] (2) to node[above,la]{$Fg$} (3);
\draw[->] (3) to node[right,la]{$Fh$} (4);
\draw[->] (1) to node[below,la]{$F(hgf)$} (4);
\draw[->] (2) to node[over,la]{$F(hg)$} (4);
\coordinate(a) at ($(2)!0.5!(4)$);
\cell[la,left,yshift=7pt][n][0.45]{a}{3}{$\widetilde F_{g,h}$};
\coordinate(a) at ($(1)!0.22!(4)$);
\coordinate(b) at ($(a)+(0,1.5cm)$);
\cell[la,right][n][0.33]{a}{b}{$\widetilde F_{f,hg}$};
\end{tzbig}
        \end{enumerate}
\end{defn}

 The following can be deduced from \cite[\S3]{LackPaoli2Nerves} or \cite[\S2.3.3]{Gurski3D}.
 
 \begin{rmk}
 \label{npsCoherence}
If $F\colon \cA \to \cB$ is a normal pseudofunctor, the compositor $\widetilde F$ can be seen as a map into $2\Iso\cB$ by post-composing with $e\colon \Triangle\cB\to 2\Iso\cB$. Then, using (c) and (f), this map can be uniquely extended to a function
\[\Mor\cA\ttimes{\Ob \cA}\ldots\ttimes{\Ob\cA}\Mor\cA\to 2 \Iso\cB,\]
    which sends a sequence of composable $1$-morphisms  $f_1, \ldots, f_k$ in $\cA$ to a $2$-isomorphism in~$\cB$ of the form
 \[
 F(f_k) \circ \ldots \circ F(f_1) \cong F(f_k \circ \ldots \circ f_1).
 \]
 \end{rmk}

 \begin{lem}
 \label{AdditionalFunctoriality}
 Let $\cT$ be a cofibrant $2$-category, namely a $2$-category whose underlying $1$-category $\Ob_*\cT$ is free. Then any normal pseudofunctor $F\colon \cA \to \cB$ induces a function 
\[
F_*\colon\twocat(\cT,\cA)\xrightarrow{}\twocat(\cT,\cB).
\]
which is natural in $\cA$, $\cB$ and cofibrant $\cT$ with respect to strict $2$-functors.
 \end{lem}

 \begin{proof}
Let $G\colon\cT\to\cA$ be a $2$-functor.
Consider the following data:
\begin{enumerate}[leftmargin=*]
\addtocounter{enumi}{-1}
    \item $(F_*G)_0\colon\Ob\cT\to\Ob\cB$ defined as $F_*G(x)\coloneqq F(G(x))$ on an object $x$ in $\cT$;
    \item $(F_*G)_1\colon\Mor\cT\to\Mor\cB$ defined as $F_*G(f)\coloneqq F(G(f))$ on a generating $1$-morphism~$f$ in $\cT$, and extended appropriately to obtain a functor $\Ob_*\cT\to \Ob_*\cA$, taking advantage of the fact that $\Ob_*\cT$ is a free $1$-category;
    \item $(F_*G)_2\colon2\Mor\cT\to2\Mor\cB$
    with $F_*G(\alpha)$ defined on a $2$-cell $\alpha\colon f_k\circ\ldots\circ f_1\Rightarrow g_l\circ \ldots\circ g_1$ in~$\cT$ as the composite
\begin{tz}[node distance=.8cm]
\node[](1) {$F(G(f_k))\circ \ldots \circ F(G(f_1))$};
\node[below of=1](2) {$F(G(f_k)\circ \ldots \circ G(f_1))$};
\node[below of=2](3) {$F(G(f_k\circ \ldots \circ f_1))$};
\node[right of=1,xshift=5.3cm](1') {$F(G(g_l))\circ \ldots \circ F(G(g_1))$};
\node[below of=1'](2') {$F(G(g_l)\circ \ldots \circ (g_1))$};
\node[below of=2'](3') {$F(G(g_l\circ \ldots \circ g_1))$};

\node at ($(1)!0.5!(2)+(1pt,0)$) {\rotatebox{270}{$\cong$}};
\node at ($(2)!0.5!(3)$) {\rotatebox{270}{$=$}};
\node at ($(1')!0.5!(2')+(1pt,0)$) {\rotatebox{270}{$\cong$}};
\node at ($(2')!0.5!(3')$) {\rotatebox{270}{$=$}};

\draw[n] (1) to node[above,la]{$F_*G(\alpha)$} (1'); 
\draw[n] (3) to node[below,la]{$F(G(\alpha))$} (3');
\end{tz}
    which involves the $2$-isomorphisms for $F$ from \cref{npsCoherence}, and the fact that $G$ preserves compositions strictly.
\end{enumerate}
It remains to see that this does indeed define a $2$-functor $F_* G$. It is clear by construction that $F_* G$ preserves compositions of $1$-morphisms. Then, it preserves horizontal compositions of $2$-morphisms by $2$-naturality of $\widetilde{F}$, and vertical compositions of $2$-morphisms since both $F$ and $G$ preserve those strictly. Note that $F_* G$ preserves $1$- and $2$-identities since both $F$ and $G$ preserve them strictly. 

The desired naturality follows from the definitions.
 \end{proof}

Recall from e.g.~\cite[\S11]{rezkTheta} or \cite[\S7.1]{ara} (resp.~\cite[Def.\ 3.1]{BarwickSchommerPries}) that a $2$-category is said to be \emph{rigid} (resp.~\emph{gaunt}) 
if it has no non-identity invertible $1$- and $2$-morphisms. Examples of gaunt $2$-categories to which we apply the following lemma in this paper are the $2$-categories $\theta$ which are objects of $\Theta_2$. 

Throughout this section, we follow the notational convention that
\[\ts{\alpha^*}{\beta^*}\colon \NMos_{i',j',0}\cA\to \NMos_{i,j,0}\cB \]
denotes the simplicial map induced by the simplicial operators $\alpha\colon [i]\to [i']$, $\beta\colon [j]\to [j']$, and $\id\colon [0]\to [0]$. 
 
\begin{prop}
For any gaunt $2$-category $\cA$ and any $2$-category $\cB$ there is a natural function
\[\NMos\colon \twocat_{\nps}(\cA,\cB)\to\spsh{(\Delta\times\Delta)}(\NMos\cA,\NMos\cB).\]
\end{prop}

 \begin{proof}
 This follows directly from \cref{AdditionalFunctoriality} with $\cT=L^{\simeq}\mathbb C\Delta[i,j,k]$, using that $L^{\simeq}\mathbb C$ is a left Quillen functor by \cite[Thm~6.1.1]{MoserNerve} and hence that every $2$-category in its image is cofibrant.
 \end{proof}

\begin{rmk}
For a $2$-category $\cB$, we give explicit relations between the sets $\NMos_{i,j,k}\cB$ for low values of $i,j,k\ge0$ with the structural data of $\cB$.
\begin{itemize}[leftmargin=*]
    \item For $(i,j,k)=(0,0,0)$ there is a bijection \[ \NMos_{0,0,0}\cB\cong\twocat([0],\cB)=\Ob\cB, \] 
    \item For $(i,j,k)=(1,0,0)$ there a bijection \[ \NMos_{1,0,0}\cB\cong\twocat([1],\cB)=\Mor\cB, \]
\item For $(i,j,k)=(1,1,0)$ there is an inclusion 
 \[\NMos_{1,1,0}\cB\leftarrow
 \twocat(\Osim{1}\tensoricon\Osim{1},\cB)\cong\twocat(\Sigma[1],\cB)\cong2\Mor\cB  \]
 induced by the map $L^{\simeq}\bC\Delta[1,1,0]\to L\bC[1,1,0]$ from \cref{LandLsimeq}; note that this inclusion can also be obtained as the pullback
\begin{tz}
\node[](1) {$2\Mor\cB$}; 
\node[below of=1](2) {$\NMos_{1,1,0}\cB$}; 
\node[right of=1,xshift=3cm](3) {$\NMos_{0,0,0}\cB\times \NMos_{0,0,0}\cB$}; 
\node[below of=3](4) {$\NMos_{0,1,0}\cB\times \NMos_{0,1,0}\cB$}; 

\draw[->] (1) to (2);
\draw[->] (1) to (3);
\draw[->] (2) to node[below,la]{$(\ts{d_1}{\id},\ts{d_0}{\id})$} (4);
\draw[->] (3) to node[right,la]{$\ts{s_0}{\id}\times \ts{s_0}{\id}$} (4);

\pullback{1};
\end{tz}
which only makes use of the simplicial structure of $\NMos\cB$.
\item For $(i,j,k)=(2,0,0)$ there is a bijection
 \[\NMos_{2,0,0}\cB\cong\twocat(\Osim{2},\cB)=\Triangle\cB.\]
\end{itemize}
\end{rmk}

\begin{rmk}
If $\cA$ is a gaunt $2$-category, the following relations hold.
\begin{itemize}[leftmargin=*]
    \item For $(i,j,k)=(1,1,0)$ there is a bijection
    \[\NMos_{1,1,0}\cA\cong 2\Mor\cA, \]
    \item For $(i,j,k)=(2,0,0)$ there is a bijection
    \[ \NMos_{2,0,0}\cA\cong \Mor\cA\ttimes{\Ob\cA}\Mor\cA, \]
    \item For $(i,j,k)=(2,1,0)$ there is a bijection
    \[ \NMos_{2,1,0}\cA\cong 2\Mor\cA\ttimes{\Ob\cA}2\Mor\cA, \]
    \item For $(i,j,k)=(3,0,0)$ there is a bijection
    \[ \NMos_{3,0,0}\cA\cong \Mor\cA\ttimes{\Ob\cA}\Mor\cA\ttimes{\Ob\cA}\Mor\cA, \]
    \item For $(i,j,k)=(1,2,0)$ there is a bijection
    \[ \NMos_{1,2,0}\cA\cong 2\Mor\cA\ttimes{\Mor\cA}2\Mor\cA. \] 
\end{itemize}
\end{rmk}

\begin{prop}
For any gaunt $2$-category $\cA$ and any $2$-category $\cB$ there is a natural function
\[\Cadj\colon\spsh{(\Delta\times\Delta)}(\NMos\cA,\NMos\cB)\to\twocat_{\nps}(\cA,\cB).\]
\end{prop}

\begin{proof}
Given a map $f\colon \NMos\cA\to \NMos\cB$ in $\spsh{(\Delta\times\Delta)}$, we produce a normal pseudo\-functor $\Cadj f\colon\cA\to\cB$ as follows:
\begin{enumerate}[leftmargin=*]
\setcounter{enumi}{-1}
\item the assignment on objects, $(\Cadj f)_0\colon \Ob\cA\to \Ob\cB$, is given by 
\[ f_{0,0,0}\colon \NMos_{0,0,0}\cA\to \NMos_{0,0,0}\cB; \]
\item the assignment on $1$-morphisms, $(\Cadj f)_1\colon \Mor\cA\to \Mor\cB$, is given by
\[ f_{1,0,0}\colon \NMos_{1,0,0}\cA\to \NMos_{1,0,0}\cB; \]
\item the assignment on $2$-morphisms, $(\Cadj f)_2\colon 2\Mor\cA\to 2\Mor\cB$, is induced by 
\[ f_{1,1,0}\colon \NMos_{1,1,0}\cA\to \NMos_{1,1,0}\cB \]
by requesting that $(\Cadj f)_2$ is the unique map that fits into the following commutative diagram:
\begin{tz}
\node[](1) {$2\Mor\cA$}; 
\node[below of=1](2) {$\NMos_{1,1,0}\cA$}; 
\node[right of=1,xshift=1.3cm](3) {$2\Mor\cB$}; 
\node[below of=3](4) {$\NMos_{1,1,0}\cB$};
\draw[->,dashed] (1) to node[above,la]{$(\Cadj f)_2$} (3); 
\draw[->] (1) to node[left,la]{$\cong$} (2);
\draw[->] (2) to node[below,la]{$f_{1,1,0}$} (4); 
\draw[->] (3) to (4);
\end{tz}
\item the compositor $\widetilde{\Cadj f}\colon \Mor\cA\ttimes{\Ob\cA}\Mor\cA\to \Triangle\cB$ is induced by 
\[ f_{2,0,0}\colon \NMos_{2,0,0}\cA\to \NMos_{2,0,0}\cB \]
by requesting that $\widetilde{\Cadj f}$ is the unique map that fits into the following commutative diagram:
\begin{tz}
\node[](1) {$\Mor\cA\ttimes{\Ob\cA}\Mor\cA$}; 
\node[below of=1](2) {$\NMos_{2,0,0}\cA$}; 
\node[right of=1,xshift=2cm,yshift=3.5pt](3) {$\Triangle\cB$}; 
\node[below of=3,yshift=-3.5pt](4) {$\NMos_{2,0,0}\cB$};
\draw[->,dashed] ($(1.east)+(0,3.5pt)$) to node[above,la]{$\widetilde{\Cadj f}$} (3); 
\draw[->] ($(1.south)+(0,2pt)$) to node[left,la]{$\cong$} (2);
\draw[->] (2) to node[below,la]{$f_{2,0,0}$} (4); 
\draw[->] (3) to node[right,la]{$\cong$} (4);
\end{tz}
\end{enumerate}
We verify that $\Cadj f$ does indeed define a normal pseudofunctor. 
\begin{enumerate}[leftmargin=*, label=(\alph*)]
    \item The compatibility of $\Cadj f$ with source, target, and identities follows from the commutativity of the following diagrams: 
    \begin{tz}
    \node[](1) {$\NMos_{0,0,0}\cA$}; 
    \node[below of=1](2) {$\NMos_{0,0,0}\cB$}; 
    \node[right of=1,xshift=2.5cm](1') {$\NMos_{1,0,0}\cA$}; 
    \node[below of=1'](2') {$\NMos_{1,0,0}\cB$}; 
    \node[right of=1',xshift=2.5cm](1'') {$\NMos_{1,1,0}\cA$}; 
    \node[below of=1''](2'') {$\NMos_{1,1,0}\cB$};
    
    \draw[->] ($(1'.west)+(0,6pt)$) to node[above,la]{$\ts{d_1}{\id}$} ($(1.east)+(0,6pt)$);
    \draw[->] ($(1'.west)-(0,6pt)$) to node[below,la]{$\ts{d_0}{\id}$} ($(1.east)-(0,6pt)$);
    \draw[->] (1) to node[over,la]{$\ts{s_0}{\id}$} (1');
    \draw[->] ($(2'.west)+(0,6pt)$) to node[above,la]{$\ts{d_1}{\id}$} ($(2.east)+(0,6pt)$);
    \draw[->] ($(2'.west)-(0,6pt)$) to node[below,la]{$\ts{d_0}{\id}$} ($(2.east)-(0,6pt)$);
    \draw[->] (2) to node[over,la]{$\ts{s_0}{\id}$} (2');
    
    \draw[->] ($(1''.west)+(0,6pt)$) to node[above,la]{$\ts{\id}{d_1}$} ($(1'.east)+(0,6pt)$);
    \draw[->] ($(1''.west)-(0,6pt)$) to node[below,la]{$\ts{\id}{d_0}$} ($(1'.east)-(0,6pt)$);
    \draw[->] (1') to node[over,la]{$\ts{\id}{s_0}$} (1'');
    \draw[->] ($(2''.west)+(0,6pt)$) to node[above,la]{$\ts{\id}{d_1}$} ($(2'.east)+(0,6pt)$);
    \draw[->] ($(2''.west)-(0,6pt)$) to node[below,la]{$\ts{\id}{d_0}$} ($(2'.east)-(0,6pt)$);
    \draw[->] (2') to node[over,la]{$\ts{\id}{s_0}$} (2'');
    
    \draw[->] (1) to node[left,la]{$f_{0,0,0}$} (2);
    \draw[->] (1') to node[left,la]{$f_{1,0,0}$} (2');
    \draw[->] (1'') to node[right,la]{$f_{1,1,0}$} (2'');
    \end{tz}
    \item The boundaries of $\widetilde{\Cadj f}$ satisfy the required condition because of the commutativity of the following diagram:
    \begin{tz}
    \node[](1) {$\NMos_{1,0,0}\cA$}; 
    \node[below of=1](2) {$\NMos_{1,0,0}\cB$}; 
    \node[right of=1,xshift=2cm](1') {$\NMos_{2,0,0}\cA$}; 
    \node[below of=1'](2') {$\NMos_{2,0,0}\cB$}; 
    \node[right of=1',xshift=3cm,yshift=-4pt](1'') {$\NMos_{1,0,0}\cA\ttimes{\NMos_{0,0,0}\cA} \NMos_{1,0,0}\cA$}; 
    \node[below of=1''](2'') {$\NMos_{1,0,0}\cB\ttimes{\NMos_{0,0,0}\cB} \NMos_{1,0,0}\cB$}; 
    
    \draw[->] (1') to node[above,la]{$\ts{d_1}{\id}$} (1);
    \draw[->] (2') to node[below,la]{$\ts{d_1}{\id}$} (2);
    
    \draw[->] (1') to node[above,la]{$\ts{(d_2,d_0)}{\id}$} node[below,la]{$\cong$} ($(1''.west)+(0,4pt)$);
    \draw[->] (2') to node[below,la]{$\ts{(d_2,d_0)}{\id}$} ($(2''.west)+(0,4pt)$);
    
    \draw[->] (1) to node[left,la]{$f_{1,0,0}$} (2);
    \draw[->] (1') to node[left,la]{$f_{2,0,0}$} (2');
    \draw[->] ($(1''.south)+(0,2pt)$) to node[right,la]{$f_{1,0,0}\times f_{1,0,0}$} ($(2''.north)-(0,4pt)$);
    \end{tz}
    \item The compatibility of $\widetilde{\Cadj f}$ with identities follows from the commutativity of the following diagram:
    \begin{tz}
    \node[](1) {$\NMos_{1,0,0}\cA$}; 
    \node[below of=1](2) {$\NMos_{1,0,0}\cB$}; 
    \node[right of=1,xshift=2cm](1') {$\NMos_{2,0,0}\cA$}; 
    \node[below of=1'](2') {$\NMos_{2,0,0}\cB$}; 
    \node[right of=1',xshift=2cm](1'') {$\NMos_{1,0,0}\cA$}; 
    \node[below of=1''](2'') {$\NMos_{1,0,0}\cB$};
        \draw[->] (1) to node[above,la]{$\ts{s_1}{\id}$} (1');
    \draw[->] (2) to node[below,la]{$\ts{s_1}{\id}$} (2');
    
    \draw[->] (1'') to node[above,la]{$\ts{s_0}{\id}$} (1');
    \draw[->] (2'') to node[below,la]{$\ts{s_0}{\id}$} (2');
    
    \draw[->] (1) to node[left,la]{$f_{1,0,0}$} (2);
    \draw[->] (1') to node[left,la]{$f_{2,0,0}$} (2');
    \draw[->] (1'') to node[right,la]{$f_{1,0,0}$} (2'');
    \end{tz}
    \item The fact that $\Cadj f$ preserves vertical composition of $2$-morphisms strictly follows from the commutativity of the following diagram:
    \begin{tz}
    \node[](1) {$\NMos_{1,1,0}\cA\ttimes{\NMos_{1,0,0}\cA} \NMos_{1,1,0}\cA$}; 
    \node[right of=1,xshift=5cm](2) {$\NMos_{1,1,0}\cB\ttimes{\NMos_{1,0,0}\cB} \NMos_{1,1,0}\cB$}; 
    \node[below of=1](1') {$\NMos_{1,2,0}\cA$}; 
    \node[below of=2](2') {$\NMos_{1,2,0}\cB$}; 
    \node[below of=1'](1'') {$\NMos_{1,1,0}\cA$}; 
    \node[below of=2'](2'') {$\NMos_{1,1,0}\cB$}; 
    
    \draw[->] (1') to node[left,la]{$\ts{\id}{(d_2,d_0)}$} node[right,la]{$\cong$} ($(1.south)+(0,2pt)$);
    \draw[->] (2') to node[right,la]{$\ts{\id}{(d_2,d_0)}$} ($(2.south)+(0,2pt)$);
    
    \draw[->] (1') to node[left,la]{$\ts{\id}{d_1}$} (1'');
    \draw[->] (2') to node[right,la]{$\ts{\id}{d_1}$} (2'');
    
    \draw[->] ($(1.east)+(0,4pt)$) to node[above,la]{$f_{1,1,0}\times f_{1,1,0}$} ($(2.west)+(0,4pt)$);
    \draw[->] (1') to node[above,la]{$f_{1,2,0}$} (2');
    \draw[->] (1'') to node[below,la]{$f_{1,1,0}$} (2'');
    \end{tz}
    \item The $2$-naturality of $\widetilde{\Cadj f}$ follows from the commutativity of the following diagram:
    \end{enumerate}
    \begin{tz}
    \node[](1) {$\NMos_{2,0,0}\cA\ttimes{X}(\NMos_{1,1,0}\cA\ttimes{\NMos_{0,0,0}\cA}\NMos_{1,1,0}\cA)$}; 
    \node[right of=1,xshift=6.2cm](2) {$\NMos_{2,0,0}\cB\ttimes{Y}(\NMos_{1,1,0}\cB\ttimes{\NMos_{0,0,0}\cB}\NMos_{1,1,0}\cB)$}; 
    \node[below of=1](1') {$\NMos_{2,1,0}\cA$}; 
    \node[below of=2](2') {$\NMos_{2,1,0}\cB$}; 
    \node[below of=1',yshift=-4pt](1'') {$\NMos_{1,1,0}\cA\ttimes{\NMos_{1,0,0}\cA} \NMos_{2,0,0}\cA$}; 
    \node[below of=2',yshift=-4pt](2'') {$\NMos_{1,1,0}\cB\ttimes{\NMos_{1,0,0}\cB} \NMos_{2,0,0}\cB$}; 
    
    \draw[->] (1') to node[left,la]{$(\ts{\id}{d_1}, \ts{(d_2,d_0)}{\id})$} node[right,la]{$\cong$} ($(1.south)+(0,8pt)$);
    \draw[->] (2') to node[right,la]{$(\ts{\id}{d_1},\ts{(d_2,d_0)}{\id})$} ($(2.south)+(0,8pt)$);
    
    \draw[->] (1') to node[left,la]{$(\ts{d_1}{\id},\ts{\id}{d_0})$} ($(1''.north)-(0,4pt)$);
    \draw[->] (2') to node[right,la]{$(\ts{d_1}{\id},\ts{\id}{d_0})$} ($(2''.north)-(0,4pt)$);
    
    \draw[->] ($(1.east)+(0,4pt)$) to node[above,la]{$f_{2,0,0}\times$} node[below,la]{$(f_{1,1,0}\times f_{1,1,0})$} ($(2.west)+(0,4pt)$);
    \draw[->] (1') to node[above,la]{$f_{2,1,0}$} (2');
    \draw[->] ($(1''.east)+(0,4pt)$) to node[below,la]{$f_{1,1,0}\times f_{2,0,0}$} ($(2''.west)+(0,4pt)$);
    \end{tz}
    \hfill\begin{minipage}{0.948\textwidth}
    where 
    \[ X=\NMos_{1,0,0}\cA\ttimes{\NMos_{0,0,0}\cA}\NMos_{1,0,0} \quad \text{and}\quad Y=\NMos_{1,0,0}\cB\ttimes{\NMos_{0,0,0}\cB} \NMos_{1,0,0}\cB. \] The fact that we retrieve the diagram of \cref{normalpseudofunctor}(e) comes from the fact that $\NMos_{2,1,0}\cB$ is the following pullback
    \begin{tz}
    \node[](1) {$\NMos_{2,1,0}\cB$}; 
    \node[right of=1,xshift=5cm,yshift=-4.5pt](2) {$\NMos_{2,0,0}\cB\ttimes{Y}(\NMos_{1,1,0}\cB\ttimes{\NMos_{0,0,0}\cB}\NMos_{1,1,0}\cB)$}; 
    \node[below of=1,yshift=-5pt](3) {$\NMos_{1,1,0}\cB\ttimes{\NMos_{1,0,0}\cB} \NMos_{2,0,0}\cB$}; 
    \node[below of=2,yshift=5pt](4) {$2\Mor\cB$}; 
    
    \draw[->] ($(1.east)+(0,.5pt)$) to node[above,la]{$(\ts{\id}{d_1},\ts{(d_2,d_0)}{\id})$} ($(2.west)+(0,5pt)$);
    \draw[->] (1) to node[left,la]{$(\ts{d_1}{\id},\ts{\id}{d_0})$} ($(3.north)-(0,4pt)$);
    \draw[->] ($(2.south)+(0,10pt)$) to node[right,la]{$\Phi$} (4);
    \draw[->] ($(3.east)+(0,5pt)$) to node[below,la]{$\Psi$} (4);
    
    \pullback{1};
    \end{tz}
    where $\Phi$ and $\Psi$ compute the total composite of the pasting diagrams. \end{minipage}
    \begin{enumerate}[leftmargin=*,label=(\alph*)]\addtocounter{enumi}{5}
        \item The compatibility of $\widetilde{\Cadj f}$ with respect to composition of $1$-morphisms follows from the commutativity of the following diagram:
    \begin{tz}
    \node[](1) {$\NMos_{2,0,0}\cA\ttimes{\NMos_{1,0,0}\cA}\NMos_{2,0,0}\cA$}; 
    \node[right of=1,xshift=5cm](2) {$\NMos_{2,0,0}\cB\ttimes{\NMos_{1,0,0}\cB}\NMos_{2,0,0}\cB$}; 
    \node[below of=1](1') {$\NMos_{3,0,0}\cA$}; 
    \node[below of=2](2') {$\NMos_{3,0,0}\cB$}; 
    \node[below of=1',yshift=-4pt](1'') {$\NMos_{2,0,0}\cA\ttimes{\NMos_{1,0,0}\cA} \NMos_{2,0,0}\cA$}; 
    \node[below of=2',yshift=-4pt](2'') {$\NMos_{2,0,0}\cB\ttimes{\NMos_{1,0,0}\cB} \NMos_{2,0,0}\cB$}; 
    
    \draw[->] (1') to node[left,la]{$\ts{(d_2,d_0)}{\id}$} node[right,la]{$\cong$} ($(1.south)+(0,2pt)$);
    \draw[->] (2') to node[right,la]{$\ts{(d_2,d_0)}{\id}$} ($(2.south)+(0,2pt)$);
    
    \draw[->] (1') to node[left,la]{$\ts{(d_3,d_1)}{\id}$} node[right,la]{$\cong$} ($(1''.north)-(0,4pt)$);
    \draw[->] (2') to node[right,la]{$\ts{(d_3,d_1)}{\id}$} ($(2''.north)-(0,4pt)$);
    
    \draw[->] ($(1.east)+(0,4pt)$) to node[above,la]{$f_{2,0,0}\times f_{2,0,0}$} ($(2.west)+(0,4pt)$);
    \draw[->] (1') to node[above,la]{$f_{3,0,0}$} (2');
    \draw[->] ($(1''.east)+(0,4pt)$) to node[below,la]{$f_{2,0,0}\times f_{2,0,0}$} ($(2''.west)+(0,4pt)$);
    \end{tz}
    The fact that we retrieve the diagram of \cref{normalpseudofunctor}(f) comes from the fact that $\NMos_{3,0,0}\cB$ is the following pullback
    \begin{tz}
    \node[](1) {$\NMos_{3,0,0}\cB$}; 
    \node[right of=1,xshift=4cm,yshift=-4.5pt](2) {$\NMos_{2,0,0}\cB\ttimes{\NMos_{1,0,0}\cB}\NMos_{2,0,0}\cB$}; 
    \node[below of=1,yshift=-5pt](3) {$\NMos_{2,0,0}\cB\ttimes{\NMos_{1,0,0}\cB} \NMos_{2,0,0}\cB$}; 
    \node[below of=2,yshift=5pt](4) {$2\Mor\cB$}; 
    
    \draw[->] ($(1.east)+(0,.5pt)$) to node[above,la]{$\ts{(d_2,d_0)}{\id}$} ($(2.west)+(0,5pt)$);
    \draw[->] (1) to node[left,la]{$\ts{(d_3,d_1)}{\id}$} ($(3.north)-(0,4pt)$);
    \draw[->] ($(2.south)+(0,2pt)$) to node[right,la]{$\varphi$} (4);
    \draw[->] ($(3.east)+(0,5pt)$) to node[below,la]{$\psi$} (4);
    
    \pullback{1};
    \end{tz}
    where $\varphi$ and $\psi$ compute the total composite of the pasting diagrams. 
    \end{enumerate}
    The desired naturality follows from the definitions. 
\end{proof}

We will need the following auxiliary fact, asserting a type of fully faithfulness for $\NMos$ when restricted to certain $2$-categories. 

\begin{prop}
\label{bisimplicialnerveFF}
For any gaunt $2$-category $\cA$ and any $2$-category $\cB$ there is a natural bijection
\[\NMos\colon \twocat_{\nps}(\cA,\cB)\cong\spsh{(\Delta\times\Delta)}(\NMos\cA,\NMos\cB).\]
\end{prop}

\begin{proof}
We now argue that given a map $f\colon \NMos\cA\to \NMos\cB$ in $\spsh{(\Delta\times \Delta)}$ we have
\[\NMos (\Cadj f)=f.\]
Since $\NMos \cA$ is $3$-coskeletal\footnote{For a Reedy category $A$, a presheaf $A^{\op}\to \set$ is $k$-coskeletal if it is canonically isomorphic to its $k$-coskeleton, in the sense of \cite[\S3.8]{RVreedy}.} and $\NMos_{i,j,k}\cA=\NMos_{i,j,0}\cA$ for all $i,j,k\geq 0$, it is enough to check that $\NMos_{i,j,0} (\Cadj f)=f_{i,j,0}$ for any $i,j\ge0$ with $i+j\leq2$,
which we see by direct inspection. 

We now argue that
\[\Cadj {(\NMos F)}=F.\]
For this, it is enough to observe that by definition $F$ and $\Cadj {(\NMos F)}$ agree on objects, $1$- and $2$-morphisms, and on the compositors.
\end{proof}

\bibliographystyle{amsalpha}
\bibliography{ref}

\end{document}

%% file: PreambleSep2020.tex
\usepackage{amsmath}
\usepackage{latexsym,amsfonts,amssymb,mathrsfs}
\usepackage{mathdots}
\usepackage{color}
\usepackage[all,2cell,color]{xy}
\usepackage{enumitem} %allows e.g. to change labels in enumerations
\usepackage{bbm}%allows to make mathbb-like numbers
\usepackage{stmaryrd}%allows for double square brackets
\usepackage[margin=1.55in]{geometry}

\usepackage{tikz} %pictures
\usepackage{tikz-cd} %diagrams
\usetikzlibrary{arrows,calc,decorations.markings,shapes}
\pgfdeclarelayer{background}
\pgfsetlayers{background,main}

\usepackage{xparse}
\usepackage{blindtext}

\newenvironment{tz}{\begin{center}\begin{tikzpicture}}{\end{tikzpicture}\end{center}}
\newenvironment{tzbig}{\begin{center}\begin{tikzpicture}[node distance=2cm]}{\end{tikzpicture}\end{center}}
\newenvironment{tzbigbig}{\begin{center}\begin{tikzpicture}[node distance=2.3cm]}{\end{tikzpicture}\end{center}}

\newcommand{\eqarrow}{\mathrel{\begin{tikzpicture}
    [baseline=(current bounding box.south)]
    \node[](a) {};
\node[right of=a,xshift=-.6cm](b) {};
\draw ($(a.east)+(0,2pt)$) sin ($(a)!0.3!(b)+(0,6pt)$) cos ($(a)!0.5!(b)+(0,4pt)$) sin ($(a)!0.7!(b)+(0,2pt)$) cos ($(b.west)+(0,4pt)$);
\draw ($(a.east)$) -- ($(b.west)$);
  \end{tikzpicture}}}
  
\newcommand{\longeqarrow}{\mathrel{\begin{tikzpicture}
    [baseline=(current bounding box.south)]
    \node[](a) {};
\node[right of=a,xshift=-.4cm](b) {};
\draw ($(a.east)+(0,2pt)$) sin ($(a)!0.3!(b)+(0,6pt)$) cos ($(a)!0.5!(b)+(0,4pt)$) sin ($(a)!0.7!(b)+(0,2pt)$) cos ($(b.west)+(0,4pt)$);
\draw ($(a.east)$) -- ($(b.west)$);
  \end{tikzpicture}}}
    
\tikzstyle{d}=[double distance=.3ex]

\usepackage{tikz}
\usetikzlibrary{cd}
\usetikzlibrary{math}
\usetikzlibrary{fit}
\usetikzlibrary{positioning}
\usetikzlibrary{decorations.pathmorphing, arrows.meta}
\tikzset{Rightarrow/.style={double equal sign distance,>={Implies},->},
triple/.style={-,preaction={draw,Rightarrow}},
quadruple/.style={preaction={draw,Rightarrow,shorten >=0pt},shorten >=1pt,-,double,double
distance=0.2pt}}

\tikzset{%
node distance=1.5cm, la/.style={scale=0.8}, rr/.style={xshift=1.5cm},
space/.style={xshift=.5cm}, over/.style={auto=false,fill=white,inner sep=1.5pt, minimum size=0, outer sep=0},
    symbol/.style={%
        draw=none,
        every to/.append style={%
            edge node={node [sloped, allow upside down, auto=false]{$#1$}}},
            
    }, pro/.style={postaction={decorate,decoration={
        markings,
        mark=at position .5 with {\node at (0,0) {$\bullet$};}
      }},
      inner sep=.9ex,
      },
  n/.style={double equal sign distance, -implies}, e/.style={double equal sign distance}, t/.style={double distance=2.5pt, -implies, postaction={draw,-}},
}

\newcommand{\pushout}[1]{\node at ($({#1})-(10pt,-10pt)$) {$\ulcorner$};}
\newcommand{\pullback}[1]{\node at ($({#1})+(10pt,-10pt)$) {$\lrcorner$};}
\NewDocumentCommand{\punctuation}{ m m O{4pt} }{\node at ($(#1.east)-(0,#3)$) {#2};}

\def\cellslide{0.5}
\def\celllength{.2cm}

\NewDocumentCommand{\cell}{ O{} O{n} O{\cellslide} O{\celllength} m m m }{
  \coordinate (mid) at ($({#5})!{#3}!({#6})$);
  \coordinate (start) at ($(mid)!{#4}!({#5})$);
  \coordinate (end) at ($(mid)!{#4}!({#6})$);
  \draw[#2] (start) to node
  [inner sep=6pt,outer sep=0,minimum size=0,#1]{{#7}} (end);
}

\usepackage[T1]{fontenc}
 \usepackage[utf8]{inputenc}
 \usepackage{array}%allows to modify arrays, e.g. the column width
 \usepackage{blkarray} %package for matrices with labels
 \usepackage{mathtools} %package for alignment in matrices
 \usepackage{longtable}%allows to make reasonable tables, in particular with pagebreaks
 \usepackage{multicol} %allows to write in two columns
 \usepackage{hyperref}
\usepackage[capitalise]{cleveref}

\usetikzlibrary{positioning} %should allow to draw more precise arrows
\usetikzlibrary{decorations.pathreplacing}%need this to draw nice curly brackets
\usetikzlibrary{arrows, decorations.markings} %need this for nice arrows in tikz pictures, markings in particular for double arrows

\newcolumntype{D}{>{\hfil$}p{1.2cm}<{$\hfil}}

\tikzset{%
    symbol/.style={%
        draw=none,
        every to/.append style={%
            edge node={node [sloped, allow upside down, auto=false]{$#1$}}}
    }
}

\tikzset{%
scalearrow/.style n args={3}{
  decoration={
    markings,
    mark=at position (1-#1)/2*\pgfdecoratedpathlength
      with {\coordinate (#2);},
    mark=at position (1+#1)/2*\pgfdecoratedpathlength
      with {\coordinate (#3);},
    },
  postaction=decorate,
  }
}

\theoremstyle{plain}   % This is the default, anyway
\newtheorem{thm}{Theorem}[section] % numbered theorem
\makeatletter\let\c@thm\c@thm\makeatother
\newtheorem{cor}{Corollary}[section]
\makeatletter\let\c@cor\c@thm\makeatother
\newtheorem{lem}{Lemma}[section]
\makeatletter\let\c@lem\c@thm\makeatother
\newtheorem{prop}{Proposition}[section]
\makeatletter\let\c@prop\c@thm\makeatother

\makeatletter\let\c@claim\c@thm\makeatother

\makeatletter\let\c@conjecture\c@thm\makeatother

\makeatletter\let\c@wconjecture\c@thm\makeatother

\newtheorem{thmalph}{Theorem}%[thmcount]
\theoremstyle{definition}

\newtheorem{defn}{Definition}[section]
\makeatletter\let\c@defn\c@thm\makeatother
\newtheorem{const}{Construction}[section]
\makeatletter\let\c@const\c@thm\makeatother
\newtheorem{notn}{Notation}[section]
\makeatletter\let\c@notn\c@thm\makeatother

\makeatletter\let\c@convention\c@thm\makeatother

\theoremstyle{remark}

\newtheorem{rmk}{Remark}[section]
\makeatletter\let\c@rmk\c@thm\makeatother

\makeatletter\let\c@ex\c@thm\makeatother

\makeatletter\let\c@observation\c@thm\makeatother

\makeatletter\let\c@warning\c@thm\makeatother

\makeatletter\let\c@digression\c@thm\makeatother

\makeatletter\let\c@answ\c@thm\makeatother

\makeatletter\let\c@answ\c@thm\makeatother

\makeatletter\let\c@aside\c@thm\makeatother

\makeatletter\let\c@question\c@thm\makeatother

\makeatletter
\let\c@equation\c@thm
\numberwithin{equation}{section}
\makeatother

\crefname{lem}{Lemma}{Lemmas}
\crefname{thm}{Theorem}{Theorems}
\crefname{defn}{Definition}{Definitions}
\crefname{notn}{Notation}{Notations}
\crefname{const}{Construction}{Constructions}
\crefname{prop}{Proposition}{Propositions}
\crefname{rmk}{Remark}{Remarks}
\crefname{cor}{Corollary}{Corollaries}
\crefname{equation}{Display}{Displays}
\crefname{ex}{Example}{Examples}
\crefname{thmalph}{Theorem}{Theorems}
\crefname{answ}{Answer}{Answers}
\crefname{question}{Question}{Questions}

\newcounter{diagram}
\newenvironment{diagram}{\setcounter{diagram}{\value{thm}}\refstepcounter{thm}\refstepcounter{diagram}\begin{center}\normalfont{(\thediagram)}\hfill\begin{tikzpicture}[baseline=(current bounding box.center)]}
    {\end{tikzpicture}\hfill\text{ }\end{center}}

\crefname{diagram}{Diagram}{Diagrams}
\numberwithin{diagram}{section} % this was subsection

\makeatletter
\newcommand{\@bbify}[1]{
  \ifcsname b#1\endcsname
  \message{WARNING: Overwriting b#1 with blackboard letter!}
  \fi
  \expandafter\edef\csname b#1\endcsname
  {\noexpand\ensuremath{\noexpand\mathbb #1}\noexpand\xspace}}
\newcommand{\@calify}[1]{
  \ifcsname c#1\endcsname
  \message{WARNING: Overwriting c#1 with calligraphic letter!}
  \fi 
  \expandafter\edef\csname c#1\endcsname
  {\noexpand\ensuremath{\noexpand\mathcal #1}\noexpand\xspace}}
\newcommand{\@bfify}[1]{
  \ifcsname bf#1\endcsname
  \message{WARNING: Overwriting c#1 with bold letter!}
  \fi
  \expandafter\edef\csname bf#1\endcsname
  {\noexpand\ensuremath{\noexpand\mathbf #1}\noexpand\xspace}}
\newcounter{@letter}\stepcounter{@letter}
\loop\@bbify{\Alph{@letter}}\@calify{\Alph{@letter}}\@bfify{\Alph{@letter}}
\ifnum\the@letter<26\stepcounter{@letter}\repeat
\makeatother
\newcommand{\Cadj}{\gamma}

\newcommand{\cat}{\cC\!\mathit{at}}
\newcommand{\dblcat}{\cD\!\mathit{bl}\cC\!\mathit{at}}
\newcommand{\set}{\cS\!\mathit{et}}
\newcommand{\sset}{\mathit{s}\set}

\newcommand{\psh}[1]{\set^{#1^{\op}}}
\newcommand{\spsh}[1]{\sset^{#1^{\op}}}
\newcommand{\msset}{m\sset}
 \newcommand{\twocat}{2\cat}
 \newcommand{\scat}{s\cat}

 \newcommand{\vcat}[1]{\cat_{#1}}
 
\newcommand{\Fun}{\mathrm{Fun}}

\DeclareMathOperator{\Hom}{Hom}
\DeclareMathOperator{\Map}{Map} 
\DeclareMathOperator{\Mor}{Mor}
\DeclareMathOperator{\Triangle}{Comp}
\DeclareMathOperator{\Iso}{Iso}
\DeclareMathOperator{\Ob}{Ob}

\newcommand{\Alg}{\mathrm{Alg}}

\newcommand{\iconI}[2]{[#1, #2]_{\mathrm{ic}}}

\DeclareMathOperator{\core}{core}
\newcommand{\inttrunc}[1]{\tau_{\leq n}^{\text{i}}}

\DeclareMathOperator{\colim}{colim}

\DeclareMathOperator{\id}{id}

\newcommand{\aamalg}[1]{\underset{#1}{{\amalg}}} %vo: for pushouts
\newcommand{\ttimes}[1]{\underset{#1}{\times}} %for pullbacks
\DeclareMathOperator{\op}{op}

\DeclareMathOperator{\ev}{ev}
\DeclareMathOperator{\cosk}{cosk}

\newcommand{\nps}{\mathrm{nps}}

\newcommand{\NMos}{{\bf N}^{\Delta\times\Delta}}
\newcommand{\NGHL}{{\bf N}^{sc}}
\newcommand{\NOR}{{\bf N}^{t\Delta}}

\newcommand{\NLein}{{\bf N}^{\Theta_2}}

\newcommand{\Nduskin}{{\bf N}^{D}}
\newcommand{\NRezk}{{\bf N}^{R}}

\newcommand{\BRequi}{\mathfrak R}
\newcommand{\schoco}{\mathfrak{N}^{sc}}

\newcommand{\Osim}[1]{\cO_2^\sim[{#1}]}
\newcommand{\Otilde}[1]{\widetilde{\cO_2[{#1}]}}
\newcommand{\DelQuot}{\overline{\Delta}}

\newcommand{\ts}[2]{[{#1},{#2}]}

\newlist{rome}{enumerate}{7}
\setlist[rome]{label=(\roman*)}

\newcommand{\pseudo}{\otimes_{\mathrm{ps}}}
\newcommand{\tensoricon}{\otimes_{\mathrm{ic}}}
\newcommand{\dblgray}{\otimes_{\mathrm{ps}}^{\mathrm{dbl}}}

\newcommand{\chaos}{\mathrm{ch}}

\makeatletter
\newcommand\notsotiny{\@setfontsize\notsotiny{7}{8}}
\makeatother

\tikzset{twoarrowlonger/.style={double,double distance=1.5pt,
shorten <=5pt,shorten >=6pt,
decoration={markings,mark=at position -4pt with {\arrow[scale=1.75]{>}}},
preaction={decorate}}}

\AtBeginDocument{%
   \def\MR#1{}
}

\newcounter{shorten}

%% file: ComparisonOfNervesMay2022.bbl
\providecommand{\bysame}{\leavevmode\hbox to3em{\hrulefill}\thinspace}
\providecommand{\MR}{\relax\ifhmode\unskip\space\fi MR }
\providecommand{\MRhref}[2]{%
  \href{http://www.ams.org/mathscinet-getitem?mr=#1}{#2}
}
\providecommand{\href}[2]{#2}
\begin{thebibliography}{HORR21}

\bibitem[AC22]{CamarenaSets}
Omar Antol\'{i}n-Camarena, \emph{The nine model category structures on the
  category of sets}, note available at
  \url{https://www.matem.unam.mx/~omar/notes/modelcatsets.html}, retrieved in
  May 2022.

\bibitem[AL20]{AraLucas}
Dimitri Ara and Maxime Lucas, \emph{The folk model category structure on strict
  {$\omega$}-categories is monoidal}, Theory Appl. Categ. \textbf{35} (2020),
  Paper No. 21, 745--808. \MR{4105933}

\bibitem[AM20]{AraMaltsiniotisJoin}
Dimitri Ara and Georges Maltsiniotis, \emph{Joint et tranches pour les
  $\infty$-cat\'{e}gories strictes}, M\'{e}m. Soc. Math. Fr. (N.S.) (2020),
  no.~165, vi+213.

\bibitem[Ara14]{ara}
Dimitri Ara, \emph{Higher quasi-categories vs higher {R}ezk spaces}, J.
  K-Theory \textbf{14} (2014), no.~3, 701--749. \MR{3350089}

\bibitem[Bar05]{BarwickThesis}
Clark Barwick, \emph{{$(\infty, n)$-Cat as a closed model category}}, PhD
  thesis (University of Pennsylvania) available at
  \url{https://repository.upenn.edu/dissertations/AAI3165639}, 2005, retrieved
  in May 2022.

\bibitem[BC03]{BullejosCegarra}
M.~Bullejos and A.~M. Cegarra, \emph{On the geometry of 2-categories and their
  classifying spaces}, $K$-Theory \textbf{29} (2003), no.~3, 211--229.
  \MR{2028502}

\bibitem[B{\'e}n67]{Benabou}
Jean B{\'e}nabou, \emph{Introduction to bicategories}, Reports of the {M}idwest
  {C}ategory {S}eminar, Springer, Berlin, 1967, pp.~1--77. \MR{0220789}

\bibitem[Ber02]{BergerCellular}
Clemens Berger, \emph{A cellular nerve for higher categories}, Adv. Math.
  \textbf{169} (2002), no.~1, 118--175. \MR{1916373}

\bibitem[BM13]{BergerMoerdijkEnriched}
Clemens Berger and Ieke Moerdijk, \emph{On the homotopy theory of enriched
  categories}, Q. J. Math. \textbf{64} (2013), no.~3, 805--846. \MR{3094501}

\bibitem[B{\"{o}}h20]{BoehmGray}
Gabriella B{\"{o}}hm, \emph{The {G}ray monoidal product of double categories},
  Appl. Categ. Structures \textbf{28} (2020), no.~3, 477--515. \MR{4089626}

\bibitem[BR13]{br1}
Julia~E. Bergner and Charles Rezk, \emph{Comparison of models for
  {$(\infty,n)$}-categories, {I}}, Geom. Topol. \textbf{17} (2013), no.~4,
  2163--2202. \MR{3109865}

\bibitem[BR20]{br2}
\bysame, \emph{Comparison of models for {$(\infty, n)$}-categories, {II}}, J.
  Topol. \textbf{13} (2020), no.~4, 1554--1581. \MR{4186138}

\bibitem[BSP21]{BarwickSchommerPries}
Clark Barwick and Christopher Schommer-Pries, \emph{On the unicity of the
  theory of higher categories}, J. Amer. Math. Soc. \textbf{34} (2021), no.~4,
  1011--1058. \MR{4301559}

\bibitem[BV73]{BoardmanVogt}
J.~M. Boardman and R.~M. Vogt, \emph{Homotopy invariant algebraic structures on
  topological spaces}, Lecture Notes in Mathematics, Vol. 347, Springer-Verlag,
  Berlin-New York, 1973. \MR{0420609}

\bibitem[Cam20]{CampbellHoCoherent}
Alexander Campbell, \emph{A homotopy coherent cellular nerve for bicategories},
  Adv. Math. \textbf{368} (2020), 107138, 67. \MR{4088416}

\bibitem[CCG10]{CCG}
Pilar Carrasco, Antonio~M. Cegarra, and Antonio~R. Garz{\'o}n, \emph{Nerves and
  classifying spaces for bicategories}, Algebr. Geom. Topol. \textbf{10}
  (2010), no.~1, 219--274. \MR{2602835 (2011d:18010)}

\bibitem[Cis06]{cisinski}
Denis-Charles Cisinski, \emph{{Les pr\'efaisceaux comme mod\`eles des types
  d'homotopie.}}, Paris: Soci\'et\'e Math\'ematique de France, 2006 (French).

\bibitem[Cis19]{CisinskiBook}
\bysame, \emph{Higher categories and homotopical algebra}, Cambridge Studies in
  Advanced Mathematics, Cambridge University Press, 2019.

\bibitem[CKM20]{CKM}
Tim {Campion}, Chris {Kapulkin}, and Yuki {Maehara}, \emph{A cubical model for
  $(\infty, n)$-categories},
  \href{https://arxiv.org/abs/2005.07603v2}{arXiv:2005.07603v2}, 2020.

\bibitem[Cor82]{CordierHoCoherent}
Jean-Marc Cordier, \emph{Sur la notion de diagramme homotopiquement
  coh\'{e}rent}, Cahiers Topologie G\'{e}om. Diff\'{e}rentielle \textbf{23}
  (1982), no.~1, 93--112, Third Colloquium on Categories, Part VI (Amiens,
  1980). \MR{648798}

\bibitem[Cru09]{CruttwellThesis}
G.S.H. Cruttwell, \emph{Normed spaces and the change of base for enriched
  categories}, PhD Thesis (Dalhousie University) available at
  \url{http://www.reluctantm.com/gcruttw/publications/thesis4.pdf}, 2009,
  retrieved in May 2022.

\bibitem[DK80a]{DwyerKanCalculating}
W.~G. Dwyer and D.~M. Kan, \emph{Calculating simplicial localizations}, J. Pure
  Appl. Algebra \textbf{18} (1980), no.~1, 17--35. \MR{578563}

\bibitem[DK80b]{DwyerKanFunction}
\bysame, \emph{Function complexes in homotopical algebra}, Topology \textbf{19}
  (1980), no.~4, 427--440. \MR{584566}

\bibitem[DKM21]{DKM}
Brandon Doherty, Chris Kapulkin, and Yuki Maehara, \emph{Equivalence of cubical
  and simplicial approaches to $(\infty,n)$-categories},
  \href{https://arxiv.org/abs/2106.09428v2}{arXiv:2106.09428v2}, 2021.

\bibitem[Dus02]{duskin}
John~W. Duskin, \emph{Simplicial matrices and the nerves of weak
  {$n$}-categories. {I}. {N}erves of bicategories}, Theory Appl. Categ.
  \textbf{9} (2001/02), 198--308, CT2000 Conference (Como). \MR{1897816
  (2003f:18005)}

\bibitem[EK66]{EilenbergKelly}
Samuel Eilenberg and G.~Max Kelly, \emph{Closed categories}, Proc. {C}onf.
  {C}ategorical {A}lgebra ({L}a {J}olla, {C}alif., 1965), Springer, New York,
  1966, pp.~421--562. \MR{0225841}

\bibitem[GH15]{GH}
David Gepner and Rune Haugseng, \emph{Enriched {$\infty$}-categories via
  non-symmetric {$\infty$}-operads}, Adv. Math. \textbf{279} (2015), 575--716.
  \MR{3345192}

\bibitem[GHL22]{GHL}
Andrea Gagna, Yonatan Harpaz, and Edoardo Lanari, \emph{On the equivalence of
  all models for $(\infty,2)$-categories}, Journal of the London Mathematical
  Society (2022).

\bibitem[Gra74]{GrayFormal}
John~W. Gray, \emph{Formal category theory: adjointness for {$2$}-categories},
  Lecture Notes in Mathematics, Vol. 391, Springer-Verlag, Berlin-New York,
  1974. \MR{0371990}

\bibitem[Gur13]{Gurski3D}
Nick Gurski, \emph{Coherence in three-dimensional category theory}, Cambridge
  Tracts in Mathematics, vol. 201, Cambridge University Press, Cambridge, 2013.
  \MR{3076451}

\bibitem[Hau13]{HaugsengThesis}
Rune Haugseng, \emph{Weakly {E}nriched {H}igher {C}ategories}, PhD Thesis
  (Massachusetts Institute of Technology) available at
  \url{https://math.mit.edu/~hrm/thesis/haugseng-thesis.pdf}, 2013, retrieved
  in May 2022.

\bibitem[Hau15]{HaugsengRectification}
\bysame, \emph{Rectification of enriched {$\infty$}-categories}, Algebr. Geom.
  Topol. \textbf{15} (2015), no.~4, 1931--1982. \MR{3402334}

\bibitem[Hau21]{HaugsengLax}
\bysame, \emph{On lax transformations, adjunctions, and monads in
  {$(\infty,2)$}-categories}, High. Struct. \textbf{5} (2021), no.~1, 244--281.
  \MR{4367222}

\bibitem[Hir03]{hirschhorn}
Philip~S. Hirschhorn, \emph{Model categories and their localizations},
  Mathematical Surveys and Monographs, vol.~99, American Mathematical Society,
  Providence, RI, 2003. \MR{1944041}

\bibitem[Hir21]{HirschhornOvercategories}
\bysame, \emph{Overcategories and undercategories of cofibrantly generated
  model categories}, J. Homotopy Relat. Struct. \textbf{16} (2021), no.~4,
  753--768. \MR{4343079}

\bibitem[HNP19]{HarpazNuitenPrasmaCohomology}
Yonatan Harpaz, Joost Nuiten, and Matan Prasma, \emph{Quillen cohomology of
  {$(\infty,2)$}-categories}, High. Struct. \textbf{3} (2019), no.~1, 17--66.
  \MR{3939045}

\bibitem[HORR21]{HORR}
Philip Hackney, Viktoriya Ozornova, Emily Riehl, and Martina Rovelli, \emph{An
  $(\infty,2)$-categorical pasting theorem},
  \href{https://arxiv.org/abs/2106.03660v3}{arXiv:2106.03660v3}, 2021.

\bibitem[Hov99]{hovey}
Mark Hovey, \emph{Model categories}, Mathematical Surveys and Monographs,
  vol.~63, American Mathematical Society, Providence, RI, 1999. \MR{1650134
  (99h:55031)}

\bibitem[JFS17]{JFSMorita}
Theo Johnson-Freyd and Claudia Scheimbauer, \emph{({O}p)lax natural
  transformations, twisted quantum field theories, and ``even higher'' {M}orita
  categories}, Adv. Math. \textbf{307} (2017), 147--223. \MR{3590516}

\bibitem[Joy97]{JoyalDisks}
Andr{\'e} Joyal, \emph{{Disks, Duality and $\Theta$-categories}}, preprint
  available at \url{https://ncatlab.org/nlab/files/JoyalThetaCategories.pdf},
  1997, retrieved in May 2022.

\bibitem[Joy08a]{joyalnotes}
\bysame, \emph{Notes on quasi-categories}, preprint available at
  \url{https://www.math.uchicago.edu/~may/IMA/Joyal.pdf}, 2008, retrieved in
  May 2022.

\bibitem[Joy08b]{JoyalVolumeII}
\bysame, \emph{The theory of quasi-categories and its applications}, preprint
  available at
  \url{http://mat.uab.cat/~kock/crm/hocat/advanced-course/Quadern45-2.pdf},
  2008, retrieved in May 2022.

\bibitem[JT07]{JT}
Andr\'e Joyal and Myles Tierney, \emph{Quasi-categories vs {S}egal spaces},
  Categories in algebra, geometry and mathematical physics, Contemp. Math.,
  vol. 431, Amer. Math. Soc., Providence, RI, 2007, pp.~277--326. \MR{2342834}

\bibitem[JY19]{JY}
Niles Johnson and Donald Yau, \emph{A bicategorical pasting theorem},
  \href{https://arxiv.org/abs/1910.01220v1}{arXiv:1910.01220v1}, 2019.

\bibitem[JY21]{JohnsonYau}
\bysame, \emph{2-dimensional categories}, Oxford University Press, Oxford,
  2021. \MR{4261588}

\bibitem[Lac02]{lack1}
Stephen Lack, \emph{A {Q}uillen model structure for 2-categories}, $K$-Theory
  \textbf{26} (2002), no.~2, 171--205. \MR{1931220}

\bibitem[Lac04]{lack2}
\bysame, \emph{A {Q}uillen model structure for bicategories}, $K$-Theory
  \textbf{33} (2004), no.~3, 185--197. \MR{2138540}

\bibitem[Lac10]{LackIcons}
\bysame, \emph{Icons}, Appl. Categ. Structures \textbf{18} (2010), no.~3,
  289--307. \MR{2640216}

\bibitem[Lei02]{LeinsterSurvey}
Tom Leinster, \emph{A survey of definitions of {$n$}-category}, Theory Appl.
  Categ. \textbf{10} (2002), 1--70. \MR{1883478}

\bibitem[LP08]{LackPaoli2Nerves}
Stephen Lack and Simona Paoli, \emph{2-nerves for bicategories}, $K$-Theory
  \textbf{38} (2008), no.~2, 153--175. \MR{2366560}

\bibitem[Lur09a]{htt}
Jacob Lurie, \emph{Higher topos theory}, Annals of Mathematics Studies, vol.
  170, Princeton University Press, Princeton, NJ, 2009. \MR{2522659}

\bibitem[Lur09b]{LurieGoodwillie}
\bysame, \emph{$(\infty, 2)$-categories and the {G}oodwillie calculus {I}},
  \href{https://arxiv.org/abs/0905.0462v2}{arXiv:0905.0462v2}, 2009.

\bibitem[Lur18]{LurieHA}
\bysame, \emph{Higher algebra}, preprint available at
  \url{http://www.math.harvard.edu/~lurie/papers/HA.pdf}, retrieved in May
  2022, 2018.

\bibitem[MG16]{MazelGeeQuillenAdj}
Aaron Mazel-Gee, \emph{Quillen adjunctions induce adjunctions of
  quasicategories}, New York J. Math. \textbf{22} (2016), 57--93. \MR{3484677}

\bibitem[Mos20]{MoserNerve}
Lyne Moser, \emph{A double $(\infty,1)$-categorical nerve for double
  categories}, \href{https://arxiv.org/abs/2007.01848v4}{arXiv:2007.01848v4},
  2020.

\bibitem[MSV22]{MSV1}
Lyne Moser, Maru Sarazola, and Paula Verdugo, \emph{A $2\mathrm{Cat}$-inspired
  model structure for double categories}, Cah. Topol. G\'{e}om. Diff\'{e}r.
  Cat\'{e}g. \textbf{LXIII} (2022), no.~2, 184--236.

\bibitem[MZ01]{MakkaiZawadowski}
Mihaly Makkai and Marek Zawadowski, \emph{Duality for simple
  {$\omega$}-categories and disks}, Theory Appl. Categ. \textbf{8} (2001),
  114--243. \MR{1825434}

\bibitem[OR20a]{ORfundamentalpushouts}
Viktoriya Ozornova and Martina Rovelli, \emph{Fundamental pushouts of
  $n$-complicial sets},
  \href{https://arxiv.org/abs/2005.05844v1}{arXiv:2005.05844v1}, 2020.

\bibitem[OR20b]{or}
\bysame, \emph{Model structures for ({$\infty$},n)--categories on
  (pre)stratified simplicial sets and prestratified simplicial spaces}, Algebr.
  Geom. Topol. \textbf{20} (2020), no.~3, 1543--1600. \MR{4105558}

\bibitem[OR21]{Nerves2Cat}
\bysame, \emph{Nerves of 2-categories and 2-categorification of
  $(\infty,2)$-categories}, Advances in Mathematics \textbf{391} (2021),
  107948.

\bibitem[Qui67]{QuillenHA}
Daniel~G. Quillen, \emph{Homotopical algebra}, Lecture Notes in Mathematics,
  No. 43, Springer-Verlag, Berlin-New York, 1967. \MR{0223432}

\bibitem[Rez96]{RezkCat}
Charles Rezk, \emph{A model category for categories}, note available at
  \url{https://faculty.math.illinois.edu/~rezk/cat-ho.dvi}, 1996, retrieved in
  May 2022.

\bibitem[Rez01]{rezk}
\bysame, \emph{A model for the homotopy theory of homotopy theory}, Trans.
  Amer. Math. Soc. \textbf{353} (2001), no.~3, 973--1007 (electronic).
  \MR{1804411 (2002a:55020)}

\bibitem[Rez10]{rezkTheta}
\bysame, \emph{A {C}artesian presentation of weak {$n$}-categories}, Geom.
  Topol. \textbf{14} (2010), no.~1, 521--571. \MR{2578310}

\bibitem[Rie17]{RiehlCTC}
Emily Riehl, \emph{Category theory in context}, Courier Dover Publications,
  2017.

\bibitem[Rie18]{EmilyNotes}
\bysame, \emph{Complicial sets, an overture}, 2016 {MATRIX} annals, MATRIX Book
  Ser., vol.~1, Springer, Cham, 2018, pp.~49--76. \MR{3792516}

\bibitem[RV14]{RVreedy}
Emily Riehl and Dominic Verity, \emph{The theory and practice of {R}eedy
  categories}, Theory and Applications of Categories \textbf{29} (2014), no.~9,
  256--301.

\bibitem[RV22]{RiehlVerityBook}
\bysame, \emph{Elements of {$\infty$}-category theory}, Cambridge Studies in
  Advanced Mathematics, vol. 194, Cambridge University Press, Cambridge, 2022.
  \MR{4354541}

\bibitem[Str82]{StreetHandwritten}
Ross Street, \emph{Higher-dimensional nerves}, hand-written notes available at
  \url{http://maths.mq.edu.au/~street/HDN1982.pdf}, 1982, retrieved in May
  2022.

\bibitem[Str87]{StreetOrientedSimplexes}
\bysame, \emph{The algebra of oriented simplexes}, J. Pure Appl. Algebra
  \textbf{49} (1987), no.~3, 283--335. \MR{920944}

\bibitem[Str96]{street}
\bysame, \emph{Categorical structures}, Handbook of algebra, {V}ol.\ 1, Handb.
  Algebr., vol.~1, Elsevier/North-Holland, Amsterdam, 1996, pp.~529--577.
  \MR{1421811 (97j:18007)}

\bibitem[Ver08a]{VerityComplicialAMS}
Dominic Verity, \emph{Complicial sets characterising the simplicial nerves of
  strict {$\omega$}-categories}, Mem. Amer. Math. Soc. \textbf{193} (2008),
  no.~905, xvi+184. \MR{2399898}

\bibitem[Ver08b]{VerityComplicialI}
\bysame, \emph{Weak complicial sets. {I}. {B}asic homotopy theory}, Adv. Math.
  \textbf{219} (2008), no.~4, 1081--1149. \MR{2450607}

\bibitem[Ver17]{VeritySlides}
\bysame, \emph{A complicial compendium}, slides available at
  \url{https://www.cirm-math.fr/ProgWeebly/Renc1773/Verity.pdf}, 2017,
  retrieved in May 2022.

\end{thebibliography}
